\documentclass[12pt,twoside,english]{article}
\usepackage{ifthen}                % if then else Abfragen
\newboolean{boldeqnitalic}
\setboolean{boldeqnitalic}{true}

\usepackage{caption}
\usepackage{umlaute}
\usepackage[intlimits] {amsmath}   % muss for defjs# stehen
\usepackage{defjs2}                % eigene Definitionen
\usepackage{epsfig}
\usepackage{psfig}
\usepackage{graphicx}
\usepackage{psfrag}
\usepackage[usenames]{color}
\usepackage{colortbl}
\usepackage{ifthen}
\usepackage{fancyhdr}
\usepackage{rotating}
\usepackage{multirow}
\usepackage{booktabs}              % Dicke Tabellenlinien
\usepackage{amssymb}
\usepackage{amsbsy}
\usepackage{bm}                    % fette Schrift in Formeln
\usepackage{bbm}
\usepackage{babel}
\usepackage{theorem}
\usepackage{upgreek}  % gerade (roman) griechische Buchstaben z.B. \upalpha
\usepackage{pifont}   % zusaetzliche Schriften, z.B. eingekreiste Zahlen

\newcommand{\Sum}{\sum\limits}
\newcommand {\parziell}[2]{\frac{\partial#1}{\partial#2}}
\newcommand {\dparziell}[2]{\dfrac{\partial#1}{\partial#2}}
\newcommand{\1}{\mathbf{1}}
\newcommand{\Real}{\mathbb{R}}
\newcommand{\dK}{\hat{K'}}
\newcommand{\ddK}{\hat{K''}}

\definecolor{dunkelgrau}{rgb}{0.8,0.8,0.8}
\definecolor{hellgrau}{rgb}{0.90,0.90,0.90} %... slightly darker
\usepackage[numbers]{natbib}
\fancypagestyle{plain}{%
  \fancyhead{}%
  \fancyfoot[c]{\sffamily\thepage}%
}
\makeatletter                           % Leere Fuellseiten
\def\cleardoublepage{\clearpage\if@twoside \ifodd\c@page\else
  \hbox{}
  \vspace*{\fill}
  \thispagestyle{empty}
  \newpage
  \if@twocolumn\hbox{}\newpage\fi\fi\fi}
\makeatother
\begin{document}
\unitlength1.0cm
\frenchspacing
\thispagestyle{empty}

\vspace{-3mm}
\begin{center}
{\bf \large
    Is there an order-barrier $p\leq2$ for time integration \\
    in computational elasto-plasticity?
}
\end{center}

\vspace{4mm}
\ce{Bernhard Eidel${}^{1,\ast}$, Charlotte Kuhn${}^2$}

\vspace{4mm}
\ce{${}^1$ Chair of Computational Mechanics, Mechanical Engineering %Fakult\"at f\"ur Ingenieurwissenschaften
}
\ce{Universit\"at Siegen, 57076 Siegen, Paul-Bonatz-Str. 9-11,
Germany}
\ce{${}^{\ast}$\small e-mail: bernhard.eidel@uni-siegen.de,
    phone: +49 271 740-2224,
    fax: +49 271 740-2436}

\vspace*{2mm}

\ce{${}^2$ Lehrstuhl f\"ur Technische Mechanik, Fakult\"at f\"ur Maschinenbau und Verfahrenstechnik}
\ce{TU Kaiserslautern, 67653 Kaiserslautern, Germany}

%\vspace{2mm}
%\begin{center}
%{This paper is dedicated to the memory of Juan Carlos Sim\'{o} (* 1952)}
%\end{center}

\vspace{4mm}
\begin{center}
{\bf \large Abstract}
\bigskip

{\footnotesize
\begin{minipage}{14.5cm}
\noindent
This paper is devoted to the question, whether there is an order barrier $p\leq2$ for time integration in computational 
elasto-plasticity. In the analysis we use an implicit Runge-Kutta (RK) method of order $p=3$ for integrating 
the evolution equations of plastic flow within a nonlinear finite element framework. We show that two novel algorithmic 
conditions are necessary to overcome the order barrier, (i) total strains must have the same order in time as the time 
integrator itself, (ii) accurate initial data must be calculated via detecting the elastic-plastic switching point (SP) 
in the predictor step. Condition (i) is for a \emph{consistent} coupling of the global boundary value problem (BVP) with the 
local initial value problems (IVP) via displacements/strains. Condition (ii) generates consistent initial data of the IVPs. 
The third condition, which is not algorithmic but physical in nature, is that (iii) the total strain path in time must be 
smooth such that condition (i) can be fulfilled at all. This requirement is met by materials showing a sufficiently smooth 
elastic-plastic transition in the stress-strain curve. We propose effective means to fulfil conditions (i) and (ii). We 
show in finite element simulations that, if condition (iii) is additionally met, the present method yields the full, theoretical 
convergence order 3 thus overcoming the barrier $p\leq 2$ for the first time. The observed speed-up for a 3rd order RK method is considerable 
compared with Backward Euler. 
\end{minipage}
}
\end{center}

{\bf Keywords:}
Finite elements; Elasto-plasticity; Time integration; Differential-algebraic equations; Runge-Kutta; Order-barrier
%\hfill vers.\,\today  

%elasto-plasticity; time-integration; Runge-Kutta methods; finite element method; order reduction; 
%space-time coupling 
%----------------------------------------------------------------------------------------------------------------

%=========================================================================
\section{Introduction}
\label{sec-intro-DAE-OrdReduc}

Time integration is the numerical workhorse of computational plasticity. It critically determines the accuracy and efficiency 
of finite element calculations. Backward Euler has been the very standard for time integration, see \cite{Simo98} and the references therein. Second order methods $(p=2)$ have been used for plastic rate equations \cite{OrtizPopov85, GoSi91, PapadopoulosTaylor94, Wieners, Ellsiepen, EckertBaaserGross, esc2, ArtioliAuricchioVeiga2006} but remarkably, the true convergence order was measured only in very few cases. The formulation of even higher-order $p \geq 3$ time-integration schemes for computational plasticity without order reduction is a longstanding, unsolved problem. Order reduction in this context means that higher-order Runge-Kutta (RK) methods with $p \geq 3$ fail to achieve the full, theoretical convergence order for integrating plastic evolution equations. Instead, the integrators fall back to convergence order 2. But why so? And why order 2? 

Suchlike order reduction is reported by \cite{Ellsiepen} for Prandtl-Reuss plasticity (i.e. von-Mises 
plasticity with zero hardening), where the 3rd order, stiffly-accurate, diagonally implicit RK method (SDIRK) of 
\cite{Alexander} was used. As a conclusion in \cite{Ellsiepen}, third order methods are seen as not competitive 
for \emph{numerical hard problems} like elasto-plasticity without hardening. However, 
even considerable hardening cannot solve the problem of order reduction for $p\geq3$ was the observation of 
\cite{HartmannBier2008}. The same authors observed that attempts to regularize the solution of elasto-plasticity 
by means of a Perzyna-type viscosity-model do fail. Instead, {\color{black}\emph{''only if an un-physically large viscosity 
is chosen, the expected order can be achieved}''}, \cite{HartmannBier2008}, p.34. Other authors observed order reduction 
from order 3 to 2 for the application of a third-order RK scheme of Radau IIa to von-Mises plasticity, see \cite{BueSi}. 
%and \cite{Buettner-Disse}, 

Beyond these results there is -- to the best knowledge of the authors -- not any publication, where time-integration for 
elasto-plasticity in finite element simulations achieves order 3 or higher. This suggests, that there exists an 
\emph{order-barrier} $p\leq2$. In view of the relevance of time integration in computational inelasticity, this order 
barrier would be severe. It is an issue, that even the reason is not well understood. The notion of \emph{low regularity} 
of elasto-plastic problems as an explanation is rather vague and a deeper analysis is missed. At best, precise conditions 
for the failure or success of 3rd order methods should be formulated.

A famous example of an order-barrier describing the failure of integration algorithms is \emph{Dahlquist's 2nd barrier}. 
It states that A-stable multistep methods (also called BDF for backward difference formulas) have order of 
convergence at most 2. Interestingly, this order barrier is the reason, why we will not attack the barrier $p=2$ 
using multistep methods. In some cases, A-stability might be not necessary, but for stiff problems as for inelasticity, 
A-stability had turned out to be even not sufficient such that stricter and additional stability notions were introduced
like L-stability and S-stability, the latter by \cite{ProtheroRobinson}. For a nice description of 
the history of stiff problems, see \cite{HaWaRadau}. Therefore, RK methods as one-step, multistage integrators are used 
here to go for order $p=3$.

Elasto-plastic continuum constitutive laws with internal variables typically exhibit the format of
differential algebraic equations (DAE), where the evolution of plastic flow is described by ordinary differential
equations (ODE) which are controlled by the yield condition as an algebraic constraint. Since RK methods have shown 
the full, theoretical convergence order beyond order 2 for various DAEs, a true order barrier as a strict, proven upper bound
for particular RK methods can be excluded. {\color{black}Hence, the \emph{order barrier} that is analyzed in the present work
has its source in the model characteristics and is a conditional one as it will turn out.}

Elasto-plasticity from the algorithmic point of view is to solve an initial boundary value problem (IBVP), where typically 
the semi-discrete (i.e. space-discrete) DAEs are solved on the local level of quadrature points. This means that the solution 
of the IVPs in elasto-plasticity is not isolated, but coupled to the solution of a BVP. From this perspective of elasto-plasticity 
as a \emph{space-time coupled problem} we take a new fresh look onto the order barrier issue, cf. \cite{Eidel-Habil}. 

In summary, the problem of order reduction in computational elasto-plasticity is still not solved.
Despite considerable research efforts, there is no solution available and even the reason is not well
understood. The present paper attacks this problem and aims at convergence order 3. In more detail, 
our main aim is three-fold:

\begin{itemize}
\item[(i)] To put things into perspective, the algorithmic structure of inelastic stress analysis within the finite element method as a space-time coupled problem
           is briefly recalled. {\color{black} Two different concepts are described.} In Sec. \ref{sec:ElastoplasticModel}, the constitutive model of a {\color{black} 
           geometrically nonlinear, small-strain von-Mises plasticity}\footnote{\color{black}Since the model is based on the geometrically nonlinear Green-Lagrangean strain 
           tensor, which is objective opposed to the infinitesimal strain tensor, the present approach is even valid for finite rotations, where strains remain small 
           as it is the case for e.g. shell-structures.} is summarized which serves the purpose of a prototype model for applying implicit Runge-Kutta methods within 
           nonlinear finite element calculations. For the two-stage, third order Radau IIa method, which is a stiffly accurate RK-version, the solution algorithm is 
           reformulated for the chosen elasto-plasticity model.

\item[(ii)] The main thrust is to identify the reasons for order reduction and, based on that, to propose simple and effective means to overcome order reduction in 
            cases where it is possible at all. Two aspects are crucial:
            \begin{enumerate}
            \item {\color{black} The first aspect is that the numerical solution of the plastic IVP needs consistent ''initial data''. As it will turn out 
                  in the following, these ''initial data'' is the time and the corresponding total strain at the elastic-plastic transition\footnote{We use 
                  the term ''initial data'' instead of ''initial values'', since plastic strain as the proper initial values of the evolution equations do not
                  belong to these ''initial data''.}. 
                  This is a new aspect that has not been considered in standard predictor-corrector schemes so far. A possible explanation is that a 
                  missing detection of the accurate elastic-plastic transition is inconsequential for the convergence order of Backward-Euler as a 
                  linear method. 
                  \emph{But how does SP detection affect the accuracy and convergence order of higher order methods with $p\geq2$?}}
%                  We propose a simple method for the detection of the switching-point within 
%                  the elastic predictor step. As it will turn out, switching point detection is a necessary but not sufficient remedy against order reduction.
            \item The second aspect refers to the characteristic of RK methods that they are multi-stage methods. As it will turn out in the present context of 
                  inelasticity, the input of total strains is required at the RK stages in the considered time interval. Since nodal displacements are 
                  the solution of the BVP, they are the link of the global BVP to the local IVPs of elasto-plasticity. Hence, the perspective of a \emph{space-time coupling} 
                  comes into play, which can be stated more precisely: 
                  \emph{Which is the correct, i.e. consistent format for the total strains serving as input in the RK stages of the local IVPs?}                  
%                  Here we propose a higher-order approximation of the strain path and, complementary, show that the approximation error in displacement/strain 
%                  interpolation is passed to the order of convergence of time integration. Thus, higher-order polynomials can reduce the local approximation error 
%                  made in strain interpolation, which has been proved for viscoelastic material models in \cite{EidelKuhn2010}. \\
%                  A critical issue closely related to elasto-plasticity in contrast to viscoelasticity is that the strain path in time may not be smooth but instead, 
%                  typically exhibits kinks indicating the onset of plastic flow in the considered material point or in its neighborhood. Suchlike low regularity 
%                  of the strain path which is obviously related to the above mentioned switching point cannot be easily tackled by higher order approximations. It is 
%                  a property of the very continuous problem of elasto-plasticity and not only related to the time-discrete counterpart.
            \end{enumerate}
The analysis of the above aspects and questions will lead to novel consistency conditions and to the development of adequate algorithmic means to meet these conditions. 

\item[(iii)] The third thrust of the paper is an assessment of the proposed methodology in representative finite element simulations, where the plastic 
             evolution equations are integrated by a fully implicit, 2-stage (3-stage), 3rd (5th) order RK method 
             of Radau IIa class. In particular, we carefully analyze the influence of SP detection and a nonlinear approximation of the strain path 
             onto the convergence order. Moreover, the efficiency of the RK methods will be measured in terms of the speed-up compared with 
             Backward-Euler.   
\end{itemize}

%----------------------------------------------------------------------------------------------------------------
%\vfill
%\newpage
%=====================================================================================
\section{Elasto-plasticity in the finite element method}
\label{sec:ElastoplasticModel}
%=====================================================================================

%=========================================================================
\subsection{General algorithmic structure}
\label{subsec:Plastic-FEM-StandardStructure}
%-------------------------------------------------------------------------

% (\emph{BE-2010}:)
The finite element method is typically used for implicit, quasi-static simulations
of solids and structures undergoing inelastic material behavior, which are modeled
within the framework of continuum constitutive models with internal variables. 
\\[1mm]
{\bf \color{black} Classical Concept.} In standard (commercial/research) finite element codes the numerical solution 
exhibits mostly the same algorithmic structure: The variational form of the balance of momentum 
is discretized in space by finite elements, which leads to a set of nonlinear algebraic equations
with displacements as unknowns in the nodes of the triangulation. These equations are typically 
solved using Newton's method, which requires a consistent linearization for the sake of quadratic 
convergence in the asymptotic range. The displacements as the solution from the \emph{global} level
are passed over to and serve as input (in the format of a total deformation tensor as derived from displacements)
in the local IVPs of elasto-plastic flow living on the \emph{local} level of quadrature (typically Gauss) points. 
Here the time integration algorithms, which are frequently called \emph{stress-update} algorithms, integrate
the rate equations. One global-local or space-time Newton-iteration circuit is completed, if the updated quantities 
from the quadrature level, i.e. stress and the algorithmically consistent tangent, are passed over to the global level 
of equilibrium iterations. This solution scheme will be called the \emph{partitioned ansatz} and is the algorithmic 
backbone of the present work. {\color{black} Note, that here even for the case of e.g. a two-stage mid-point rule the
global solution of the global algebraic equations is required only once.}

{\color{black} {\bf An Alternative: DAE/MLNA.} An alternative is the understanding of inelasticity within finite elements as a 
global DAE system, where displacements are the algebraic unknowns and the internal variables are the differential variables. 
The corresponding solution algorithm is referred to Multi-Level-Newton-Algorithm (MLNA).
Of course, even in this approach, the DAE-problem of elasto-plasticity on the local level of Gauss points still exists.   
A description of this alternative concept is described in detail in \cite{EllsiepenHartmann2002} and follow-up papers of the second author.  
%\\[1mm]
%{\bf WAS IST DAS???} MLNA stands for Multi-Level-Newton-Algorithm is 
\\[1mm]
In the context of the present work, one aspect of the DAE-approach is of particular relevance; the set of global algebraic 
equations for equilibrium is solved not only at the end of the time interval, but additionally at each and every of the altogether
$s$ RK stages, see p.687 in \cite{EllsiepenHartmann2002}, Tab.~3. in \cite{Hartmann2002}, Tab.~3.~\cite{HartmannBier2008}. 
Since the rest of the algorithmic structure is the same as in the classical concept, the numerical effort increases by a factor 
of $s$.} 

\subsection{Von-Mises constitutive model}
\label{sec:VonMisesPlasticity}

%------------------------------------------------------------------------------------------------------------- 
%\fbox{\parbox{15.3cm}{

\begin{Table}[htbp]
\begin{tabular}{ p{15.5cm}}
\hline\\
\end{tabular}
\begin{eqnarray}
\label{Elastoplasticity-DecompStrain}
 \mbox{decomposition of strain}        & &  \bE = \bE^{\mbox{\scriptsize e}} + \bE^{\mbox{\scriptsize p}}         \\[1mm]
\label{FreeEnergy} 
 \mbox{decomposition of free energy}   & &  \Psi = W(\bE^{\mbox{\scriptsize e}}) + \hat{K}(\alpha).               \\[1mm]
\label{Elasticity law}
 \mbox{elasticity law}                 & &  \bS = \partial_{\bE^{\mbox{\scriptsize e}}}\Psi(\bE^{\mbox{\scriptsize e}}, \alpha)\, , 
                                            \bS^{\mbox{\scriptsize D}} = \mbox{dev}(\bS)                          \\[1mm]
\label{Elastoplasticity-YieldCondition}
 \mbox{yield function}                 & &  f(\bS,\alpha) =  || \bS^{\mbox{\scriptsize D}} ||-\sqrt{\frac{2}{3}}\left(\sigma_Y 
                                            + \hat{K'}(\alpha)\right)                                              \\[1mm]
\label{Elastoplasticity-FlowRules1}
\mbox{associated flow rules }          & &  \dot{\bE^{\mbox{\scriptsize p}}} = \gamma\, \frac{\partial f}{\partial \bS} = \gamma\, \frac{\bS^{\mbox{\scriptsize D}}}{||\bS^{\mbox{\scriptsize D}}||} := \gamma\, \bn               \\
\label{Elastoplasticity-FlowRules2}
                                       & &  \dot{\alpha} = \sqrt{\frac{2}{3}}\,\gamma                                    \\
\label{IsotropicHardening}
 \mbox{isotropic hardening}            & & \hat{K'}(\alpha) = H\alpha + (\sigma_\infty - \sigma_Y)(1-e^{-\delta\alpha})  \\[1mm]
\label{KKT-conditions}
 \mbox{Karush-Kuhn-Tucker conditions}  & &  \gamma\geq 0,\;\;\;f\leq0, \;\;\; \gamma f = 0  \label{KKT}                  \\[2mm]
\label{consistency-condition}
 \mbox{consistency condition}         & &  \gamma \dot{f} = 0 \quad \text{if} \quad f=0    \label{consistency-condition}
\end{eqnarray}
\begin{tabular}{ p{15.5cm}}
\hline\\[-4mm]
\end{tabular}
\caption{Elasto-plasticity model of von-Mises.\label{tab:vonMisesModel}}
\end{Table}
 
Equations \eqref{Elastoplasticity-DecompStrain}--\eqref{consistency-condition} summarize {\color{black} a small-strain} von-Mises plasticity model
in a geometrical nonlinear setting. This model will be used in the following to showcase the numerical solution of the corresponding evolution 
equations of DAE-type by RK methods within the partitioned ansatz. For the modeling of elasto-plasticity in the framework of 
continuum constitutive laws with internal variables we refer to the books of \cite{Lubliner1990, Maugin1992, Haupt2000, Lubarda2002, Bertram2012}.
For a sound introduction into the algorithmic treatment of elasto-plasticity we refer to \cite{Simo98}.

The model is based on an additive decomposition of the Green-Lagrangean strain tensor $\bE=1/2 (\bF^T \bF - \mathbf{1})$ into an elastic part 
$\bE^{\mbox{\scriptsize e}}$ and a plastic part $\bE^{\mbox{\scriptsize p}}$ according to \eqref{Elastoplasticity-DecompStrain}. Here, $\bF$ 
is the deformation gradient and $\1$ denotes the second order unity tensor with components $\delta_{ij}$. It is assumed that the free energy $\Psi$ 
in \eqref{FreeEnergy} can be additively decomposed into a hyperelastic part $W(\bE^{\mbox{\scriptsize e}})$ and a plastic part $\hat K(\alpha)$
for isotropic hardening with the equivalent plastic strain $\alpha$. In the hyperelasticity law, $\bS$ is the Second Piola-Kirchhoff stress 
tensor, which is energetically conjugate to $\bE^{\mbox{\scriptsize e}}$. The yield function $f$ is the criterion for the distinction of elastic
states from plastic states. It enters the associated flow rule for $\bE^{\mbox{\scriptsize p}}$ according to \eqref{Elastoplasticity-FlowRules1}, 
which is accompanied by the evolution equation for $\alpha$ according to \eqref{Elastoplasticity-FlowRules2}. The Lagrange-multiplier $\gamma$ 
follows from the consistency condition \eqref{consistency-condition}. For the isotropic hardening $\hat K'(\alpha)$ a format is assumed ($(\bullet)'$ 
denotes a derivative with respect to the argument, here: $\alpha$) that is well-suited for adaption to experiments by 
the composition of a linear part and a saturation part, see \eqref{IsotropicHardening}. The hardening modulus $H$, the initial yield stress $\sigma_{Y}$, 
the saturation stress $\sigma_{\infty}$ as well as $\delta$ are material parameters. Elastic-plastic loading and unloading are governed by the Karush-Kuhn-Tucker 
optimality conditions according to \eqref{KKT-conditions}.   

In the decompositions
\begin{equation}
\bE^{\mbox{\scriptsize D}} =\bE - \frac{1}{3}\tr(\bE)\mathbf{1}\,, \qquad 
\bS^{\mbox{\scriptsize D}} =\bS - p \, \mathbf{1} \quad \mbox{with} \quad p={1}/{3}\tr(\bS) \, ,
\end{equation}
$\bE^{\mbox{\scriptsize D}}$ $(\bS^{\mbox{\scriptsize D}})$ are the strain (stress) deviator, 
$\tr (\bE)=E_{ii}$ is the trace operator and $p$ is the hydrostatic pressure.

\subsection{Elasto-plasticity is a DAE}
\label{sec:ElPl-Index2-DAE}

The von-Mises elasto-plasticity model is an instance of the broad class of constitutive models of the type

\begin{equation}
%\left.
\begin{array}{rcl}
     \bS  &=& \bh(\bE, \by)                                    \\
\dot{\by} &=& \bbf(\bE, \by) \, , \quad \by(t_0)=\by_0         \\
       0  &=& f(\bS)                                       
\end{array} %\left.%\right\} %i = 1,...,s \, .
\label{eq:Elasticity-and-DAE}
\end{equation}

where $\bS = \bh(\bE, \by)$ represents an elasticity relation and $\by \in \Real^{n_z}$ denotes the set of internal variables describing 
plastic material behavior. Here, for $\by=\{\bE^{\mbox{\scriptsize p}},\alpha\}$ it holds
$n_z=6+1=7$. % holds in the 3-d case. 
The evolution equations of plastic flow, \eqref{eq:Elasticity-and-DAE}$_2$, are subject to the algebraic constraint of the yield
condition \eqref{eq:Elasticity-and-DAE}$_3$ thus forming a set of DAEs.
For plastic loading of von-Mises plasticity \eqref{eq:Elasticity-and-DAE}$_{2,3}$ read

\begin{equation}
\label{eq:vMises-DAE} 
%\left.
\begin{array}{rcl}
\dot{\bE}^{\mbox{\scriptsize p}}   &=&  \gamma \dfrac{\bE^{\mbox{\scriptsize D}}-\bE^{\mbox{\scriptsize p}}}{||\bE^{\mbox{\scriptsize D}}-\bE^{\mbox{\scriptsize p}} ||}\,,  \\
\dot{\alpha}  &=&  \gamma \sqrt{\dfrac{2}{3}}\,,                                              \\
           0  &=&  2\mu||\bE^{\mbox{\scriptsize D}}-\bE^{\mbox{\scriptsize p}} ||-\sqrt{\dfrac{2}{3}}\left(\sigma_Y+\dK(\alpha)\right)   \,,
\end{array} %\right\} % i = 1,...,s \, .
\end{equation}

where we replace for notational convenience the stress deviator $\bS^{\mbox{\scriptsize D}}$ by $2\mu(\bE^{\mbox{\scriptsize D}} -\bE^{\mbox{\scriptsize p}})$, 
$\mu$ is the shear-modulus. We call $\bE^p$ the \emph{differential variable}, $\gamma$ is the \emph{algebraic variable}.

In the present work we will use RK methods of Radau IIa class, since they are well suited for stiff problems, and 
inelasticity is a typical instance for stiff problems. For Radau IIa the theoretical convergence orders of both the 
differential variables and the algebraic variables of the elasto-plastic DAE system are 
%\eqref{Ind2-DAE-1}--\eqref{Ind2-DAE-3}
listed in the left of Tab. \ref{tab:ConvergenceOrder} (from Table VII.4.1 in \cite{HaNoWaI}) and are contrasted 
against the reduced order as it was observed in finite element simulations in \cite{BueSi}.
\begin{Table}[htbp]
\renewcommand{\arraystretch}{1.5}
\centering
\begin{tabular}{l c c }
\hline
variable type                               & theoretical order &  reduced order   \\
\hline differential variable  $\by$         & $2s-1$            &       $2$        \\
       algebraic variable     $\gamma$      & $s$               &       $2$        \\
       \hline
\end{tabular}
\newline 
\caption{Global, theoretical vs. reduced convergence order of $s$-stage ($s\geq2$) Radau IIa methods 
for the DAE case of rate-independent elasto-plasticity.\label{tab:ConvergenceOrder}}
\end{Table}

\subsection{Stiffly accurate RK methods for the integration of plastic flow}

The solution of IVPs by implicit RK (IRK) methods is a standard task of numerical mathematics. 
For a comprehensive overview we refer to the monographs \cite{HaLuRo, HaNoWaI, HaWaII}. 

The solution of the IVP of plastic flow by RK methods in $s$ stages for a time step 
$\Delta t=t_{n+1}-t_n$ is obtained as follows. The nonlinear system of equations

\begin{equation}\label{step-1}
\left.
\begin{array}{rcl}
\bY_{ni} &=& \by_n + \Delta t\Sum_{j=1}^s a_{ij} \bbf(\bY_{nj},\Gamma_{nj})  \\ %  \;\;\;\; i = 1,...,s\\
0&=&f(\bY_{ni})
\end{array}
\right\} \quad i = 1,...,s \, 
\end{equation}

is solved for $\bY_{ni}$ and $\Gamma_{ni}$ $(i = 1,...,s)$, which are the \textit{stage solutions} of $\by$ and $\gamma$.
In \eqref{step-1}, the weighting factors $a_{ij}$ are the coefficients of the \textit{Runge-Kutta-Matrix} 
$A = (a_{ij})_{i,j=1,...,s}$. The position of the time \emph{stages} $t_{ni} = t_n + c_i \Delta t_{n}$ are
prescribed by the coefficients $c_i, {i=1,...,s}$. The Butcher arrays for the Radau IIa schemes up to 
order 3 are listed in Appendix A.   

Radau IIa schemes are an instance of \textit{stiffly accurate} RK methods. They obey the relation $a_{si} = b_i$ for  
all $i = 1,...,s$. As a consequence, it holds for the differential variable $\by$ %and for $\gamma$
\begin{equation}
\label{Yn+1UpdatdStiffMethods}
\by_{n+1} = \by_n+ \Delta t \, \Sum_{j+1}^s b_j\dot{\bY}_{nj} =  \by_n + \Delta t\Sum_{j=1}^s a_{sj} \dot{\bY}_{nj}= \bY_{ns}\, .  
\end{equation}
This characteristic, that the approximate solution $\by_{n+1}$ coincides with the last stage solution $\bY_{ns}$ is favorable, 
since then, the stress tensor $\bS_{n+1}$ satisfies the yield condition because of \eqref{step-1}$_2$, $f(\by_{n+1})=0$. 
{\color{black}Moreover, as \eqref{Yn+1UpdatdStiffMethods} implies that $y_{n+1} = Y_{ns}$, $y_{n+1}$ immediately follows from the solution 
of \eqref{step-1} and the stage dervivatives do not need to be computed.}

Next, we consider the application of stiffly accurate RK-methods to Eqs.~\eqref{eq:vMises-DAE}$_{1-3}$ 
of von-Mises plasticity with isotropic hardening within the partitioned framework. Here, we introduce 
$\Delta\gamma :=\gamma\Delta t$. With the stage values $\{\bE^{\mbox{\scriptsize p}}_{ni},\Lambda_{ni}\}$ for 
$\{\bE^{\mbox{\scriptsize p}},\alpha\}$, $\Delta\Gamma_{ni}$ for $\Delta\gamma$ and $t_{ni} = t_n + c_i\Delta t$ for $t$, 
the nonlinear set of equations for the stage values reads as
\begin{equation}
\label{nonlinEqSys}
\left.
\begin{array}{rcl}
\bE^{\mbox{\scriptsize p}}_{ni} &=& \bE^{\mbox{\scriptsize p}}_n + \Sum_{j=1}^s a_{ij} \Delta\Gamma_{nj} \dfrac{\bE^{\mbox{\scriptsize D}}(t_{nj}) - \bE^{\mbox{\scriptsize p}}_{nj}}{||\bE^{\mbox{\scriptsize D}}(t_{nj}) - \bE^{\mbox{\scriptsize p}}_{nj} ||}\\
\Lambda_{ni} &=& \alpha_n + \Sum_{j=1}^s a_{ij} \sqrt{\dfrac{2}{3}} \Delta\Gamma_{nj}\\
0&=&2\mu ||\bE^{\mbox{\scriptsize D}}(t_{ni}) - \bE^{\mbox{\scriptsize p}}_{ni} || - \sqrt{\dfrac{2}{3}}\left(\sigma_Y + \dK(\Lambda_{ni})\right)
\end{array}\right\} i = 1,...,s \, .
\end{equation}
%------------------------------------------------------------------------------------------------------------------------------
\section{Reasons for and remedies against order reduction}
\label{sec:Reasons4OrderReduction}

In this section we identify and analyze two critical sources for order reduction. Based on the findings
we propose solutions to overcome the loss of the theoretical convergence order. 
More specifically, we examine the influence of
\begin{enumerate}
\item Accurate, i.e. \emph{consistent} strain interpolation in time, where not only the nonlinearity of the strain path in time
      comes into play, but moreover its \emph{smoothness}.
\item \emph{Consistent} initial data for the IVP via (SP) detection between elastic state and plastic state.
\end{enumerate}

\subsection{Consistent approximation of strain in time}
\label{subsec:DehnInterpol}

\begin{Figure}[htbp]
   \begin{minipage}{15cm}
   \centering
   %      \scalebox{.75}{\input{bild1202.pstex_t}}
   %      \scalebox{.75}{\input{bild1202.pstex_t}}
          \includegraphics[width=.70\textwidth, angle=0, clip=]{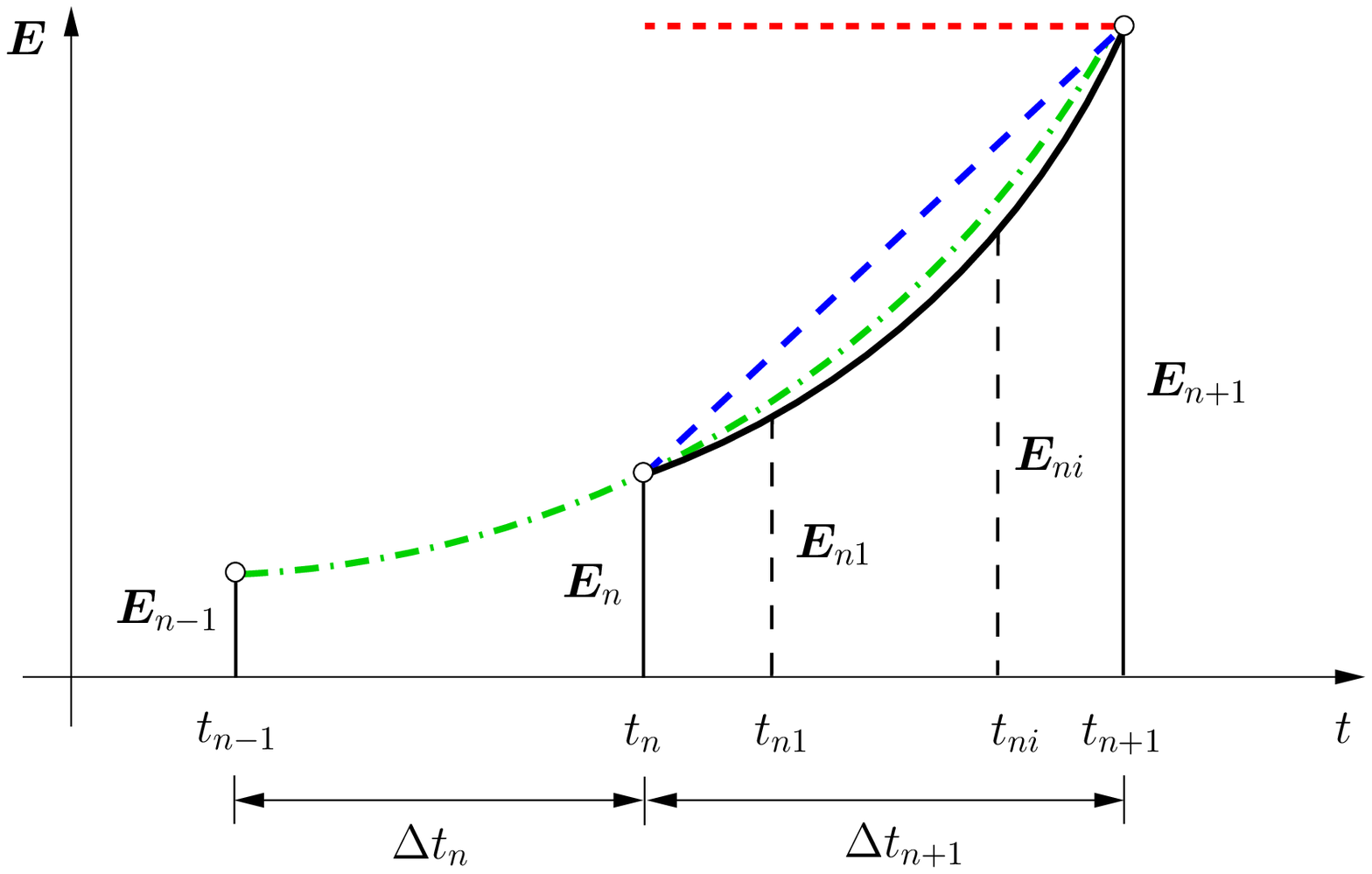}
   \end{minipage}
\caption{
The exact, unknown strain path (black, solid line) in the time interval $t \in [t_n, t_{n+1}]$ and its approximations by constant 
deformation $\bE=\bE_{n+1}$ (red, dotted line), by linear (blue, dashed line) interpolation based on $\bE_n$ and $\bE_{n+1}$ data, and 
by quadratic (green, dot and dash line) interpolation based on $\bE_n$, $\bE_{n+1}$ and $\bE_{n-1}$ data.\label{fig:Interpol-LinVsQuadraticVsExact}
}
\end{Figure}

Equations \eqref{nonlinEqSys}$_{1,3}$ indicate, that RK schemes require total strain values at their stages $\bE(t_{nj})$. 
Moreover, Fig.~\ref{fig:Interpol-LinVsQuadraticVsExact} hints at the fact, that the strain path in time of a material 
point in an inelastic solid is typically nonlinear. If the strain path is nonlinear and smooth, strain representation 
by higher-order polynomials yield a better approximation as e.g. linear interpolation can do. Keeping the total strain tensor 
constant, which is the assumption of BE within a classical predictor-corrector scheme, is a very rough approximation 
of the true strain path. 

But which is the correct or \emph{consistent} format to approximate total strains $\bE(t_{nj})$ by $\bE_{nj}$? 

\subsubsection{Linear strain approximation: order $p=2$.} \quad
Introducing a new time scale $\tilde{t} := t-t_n$ in the time
interval $t \in [t_n, t_{n+1}]$, the linear interpolation polynomial,
supported by $(0,\bE_n)$ and $(\Delta t, \bE_{n+1})$, reads
\begin{equation}
\bp_1(\tilde{t}) = \bE_n + \frac{\tilde{t}}{\Delta t}(\bE_{n+1} -\bE_n) \, ,
\end{equation}
which leads to the approximation of $\bE$ at time $t = t_{ni} = t_n + c_i \Delta t$
\begin{equation}
\bE(t = t_{ni}) \approx \bp_1(\tilde{t} = c_i\Delta t) = \bE_n + c_i (\bE_{n+1}-\bE_n)\,.
\label{linapproxE}
\end{equation}
Under the assumption that the strains in time interval $t \in [t_n, t_{n+1}]$ and $\tilde{t}\in [0,\Delta t]$ are 
twice continuously differentiable, the interpolation error in the strains is 
%\begin{equation}
%\label{eq:error-LinInterpolation}
$|\bE(\tilde{t}) - \bp_1(\tilde{t})| \leq {M_2}/{2} (\Delta t)^2$
%\end{equation}
with $M_2 = \mbox{max}\lbrace|\bE''(\tilde{t})| \,: \,\tilde{t}\in [0,\Delta t]\rbrace$. %, \cite{Bronstein}.
%The approximation of the stage values via linear interpolation was applied e.g. in \cite{Buettner-Disse}.

\subsubsection{Quadratic strain approximation: order $p=3$.} \quad
\label{subsubsec:QuadStrainInterpol}
Given that $\bE(\tilde{t})$ is three times continuously
differentiable in $[-\Delta t, \Delta t]$, it holds for the
interpolation error
%\begin{equation}
$|\bE(\tilde{t}) - \bp_2(\tilde{t})| \leq {M_3}/{6} (\Delta t)^3$
%\end{equation}
with $M_3 = \mbox{max}\lbrace|\bE'''(\tilde{t})| \,: \,\tilde{t}\in [-\Delta t,\Delta t]\rbrace$. %, \cite{Bronstein}.
The interpolation polynomial $\bp_2$ reads
\begin{equation}
\label{poly2}
\bp_2(\tilde{t}) = \frac{1}{2(\Delta t)^2} \left(\bE_{n-1}-2\bE_n +\bE_{n+1}\right)\, \tilde{t}^2 + \frac{1}{2 \Delta t} \left(\bE_{n+1} - \bE_{n-1}\right) \, \tilde{t} + \bE_n \, .
\end{equation}
Then, stage values of $\bE$ at $t=t_{ni}=t_n + c_i\Delta t$ and $\tilde{t}= c_i\Delta t$ read
\begin{equation}
\bE(t = t_{ni}) \approx \bp_2(\tilde{t} = c_i\Delta t) = \frac{c_i}{2}(c_i-1)\,\bE_{n-1} + (1-c_i^2)\,\bE_n + \frac{c_i}{2}(c_i+1)\,\bE_{n+1}\,.
\label{quadapproxE}
\end{equation}

The fact that the error for linear interpolation is of the order $\mathcal{O}(\Delta t^2)$ suggests that this 
low-order approximation is a candidate to cause order reduction in elasto-plastic stress computations when higher-order 
methods, $p \geq 3$, are used. In \cite{EidelKuhn2010} it was shown, that the convergence order in viscoelastic stress 
computations directly depends on the approximation error in strain interpolation.
More precisely, it was shown for a 3rd order RK method, that quadratic interpolation based on $t_{n-1}$, $t_{n}$ and $t_{n+1}$ data
enables convergence order 3 in viscoelastic time integration. Linear interpolation, however, leads to a reduced order of 2.

%-----------------------------------------------------------------------------------------------------------------------
\subsection{Consistent initial data via elasto-plastic switching point detection}
\label{subsec:SwitchingPoint}

\begin{Figure}[htbp]
   \begin{minipage}{15cm}
   \centering
  %     \scalebox{.70}{\input{bild1200a.pstex_t}}
      \includegraphics[width=.45\textwidth, angle=0, clip=]{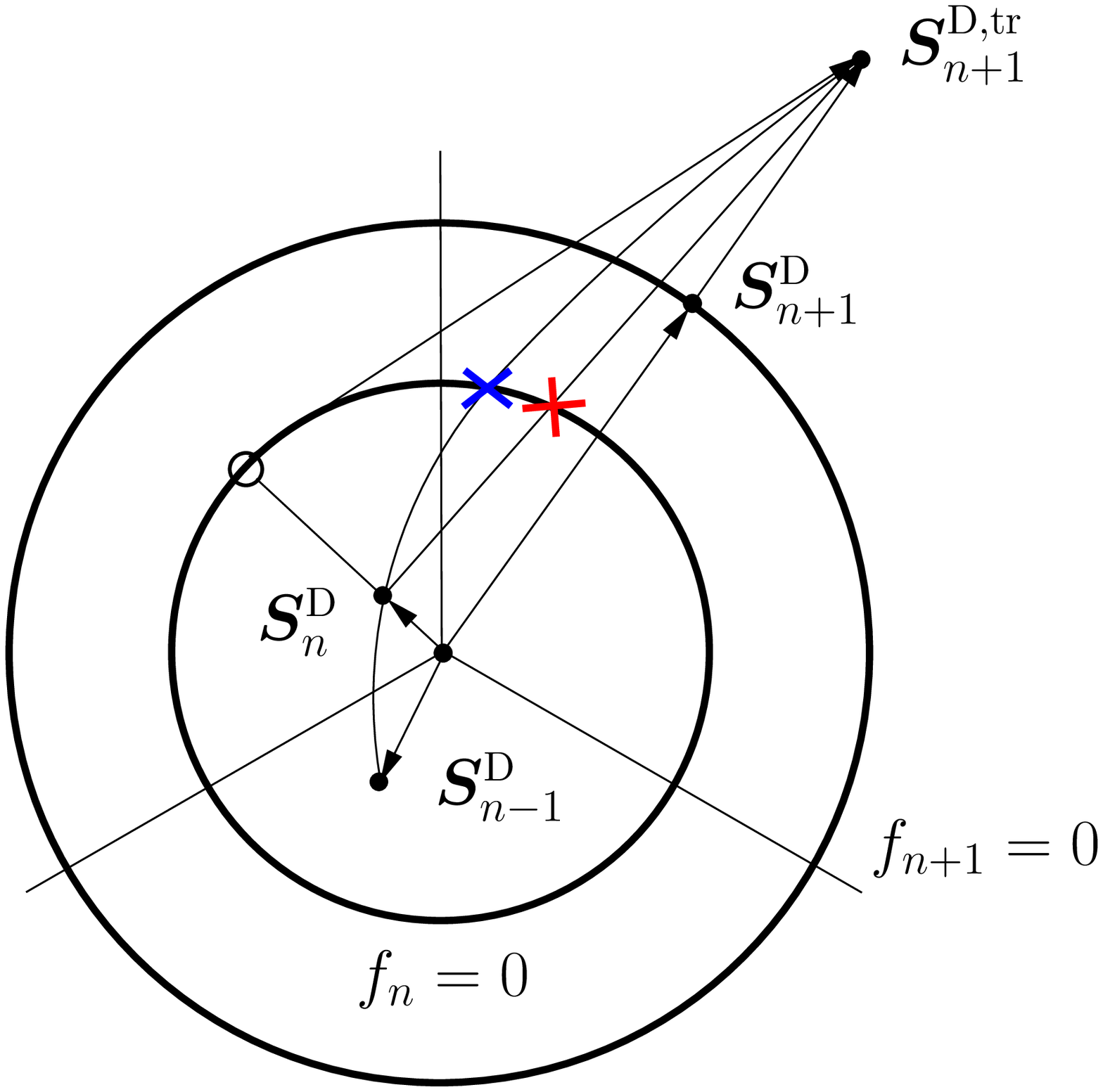}
   \end{minipage}
\caption{Switching point (SP) detection for von-Mises elasto-plasticity with isotropic hardening in the predictor step of the return-mapping scheme.
         Three versions are displayed: for linear strain interpolation {\bf \color{red}$+$}, for quadratic strain interpolation 
         {\color{blue}$\times$}, and for linear extrapolation {\bf \color{red}$\circ$}. \label{fig:SPD}}
\end{Figure}

The transition from an elastic state to a plastic state is characterized by the SP at time $t_{\mbox{\scriptsize SP}}$.
Since the SP defines initial values for IVPs of time integration in elasto-plasticity, its determination is crucial for the accuracy 
of inelastic stress calculations, in particular for higher-order methods. Linear Backward Euler and BDF2, a second order backward-difference formula, 
\cite{PapadopoulosTaylor94}, \cite{Simo98}, achieve full convergence order without SP detection.
\\[2mm] 
In \cite{Marques84} and \cite{Bicanic89}, the SP was termed \emph{contact point} and considered for final stress 
calculations, so without the use of a predictor-corrector scheme based on an operator split. In \cite{esc2} the SP 
was used for a novel exponential type integrator to decompose the strain increment into its proper elastic and plastic parts. 
Therefore, this usage of the SP comes close to the purpose it serves in the present analysis. 
\\[2mm]
During an elastic state, $f<0$, the consistency parameter $\gamma$ is zero and the plastic variables $\bE^{\mbox{\scriptsize p}}$ and 
$\alpha$ remain constant. In a plastic phase $f=0$, $\gamma > 0$, the evolution equations are active. The strain path in time $\bE(t)$ 
of a material point is continuous after the onset of plastic flow. However, it is not necessarily continuously differentiable at 
$t=t_{\mbox{\scriptsize SP}}$. Instead, it often exhibits a kink at the SP indicating a characteristic loss of stiffness 
when plastic yielding commences.
 
The trial step value of the yield condition indicates, whether there is a SP within the current time interval $(t_n,t_{n+1}]$. 
A SP exists, if it holds $f_n^{\mbox{\scriptsize tr}} < 0$ and $f_{n+1}^{\mbox{\scriptsize tr}}\geq 0$. 
 
If switching occurs, the true time step size of inelastic flow is reduced to $\Delta t=t_{n+1}-t_{\mbox{\scriptsize SP}}$ and is therefore 
smaller than the nominal time step size $\Delta t = t_{n+1}-t_n$. Since plastic variables are constant during the elastic 
state, it holds
\begin{equation}
\bE^{\mbox{\scriptsize p}}_{\mbox{\scriptsize SP}} = \bE^{\mbox{\scriptsize p}}_n\,, \qquad \alpha_{\mbox{\scriptsize SP}} = \alpha_n\,.
\end{equation}
The SP and the corresponding total strain $\bE(t_{\mbox{\scriptsize SP}})$ have to be calculated.
In the following, we consider three different strategies to determine the SP consistent with the method of 
strain interpolation. The first one is a simple ansatz for the approximate calculation of the SP which is based 
on linear interpolation. Following our ideas of higher-order interpolation of strain we adapt SP detection 
to the case of quadratic interpolation from Sec.~\ref{subsubsec:QuadStrainInterpol}. The third method is based on linear extrapolation.
All versions are visualized in Fig.~\ref{fig:SPD}.
%---------------------------------------------------------------------------
\subsubsection{Switching point detection for linear interpolation.}
%\label{SPlin}
%
If the strains within the time interval $[t_n, t_{n+1}]$ are calculated via linear interpolation according to \eqref{linapproxE}, 
the SP can be approximately calculated as described in \cite{Marques84} and similarly in \cite{Bicanic89} and 
\cite{esc2}. It is assumed that deformation in the considered material point is purely elastic in the first subinterval 
$[t_n, t_{n+x})$ with $x\in (0,1)$ and that plastic flow is restricted to $[t_{n+x},t_{n+1}]$. For the calculation of 
parameter $x$, the scalar equation
\begin{equation}
2\mu ||(\bE^{\mbox{\scriptsize D}}_n+x(\bE^{\mbox{\scriptsize D}}_{n+1}-\bE^{\mbox{\scriptsize D}}_n))-\bE^{\mbox{\scriptsize p}}_n || - \sqrt{\frac{2}{3}}\left(\sigma_Y + \dK(\alpha_n)\right) = 0
\end{equation}
has to be solved. Then, the SP is at $t_{\mbox{\scriptsize SP}}\approx t_n + x\Delta t$ with the corresponding total
strain $\bE(t_{\mbox{\scriptsize SP}})\approx \bE_n + x(\bE_{n+1}-\bE_n)$. Consequently, the purely elastic, first step is calculated 
using the reduced step size $(1-x)\Delta t$.

%---------------------------------------------------------------------------
\subsubsection{Switching point detection for quadratic interpolation.}
%\label{SPquad}
If total strains in the interval $[t_n, t_{n+1}]$ are approximated by quadratic interpolation according to \eqref{quadapproxE}, the corresponding 
equation to determine  $x\in (0,1)$ reads as 
\begin{equation}
2\mu || [\frac{x}{2}(x-1)\,\bE^{\mbox{\scriptsize D}}_{n-1} + (1-x^2)\,\bE^{\mbox{\scriptsize D}}_n + \frac{x}{2}(x+1)\bE^{\mbox{\scriptsize D}}_{n+1}] -\bE^{\mbox{\scriptsize p}}_n ||- \sqrt{\frac{2}{3}}\left(\sigma_Y + \dK(\alpha_n)\right) = 0\,.
\end{equation}
Again it holds $t_{\mbox{\scriptsize SP}}\approx t_n + x\Delta t$ and time integration starts at $t_{\mbox{\scriptsize SP}}$ along with a reduced time step size $(1-x)\Delta t$. For the total strains at
the SP it holds
\begin{equation}
\bE(t_{\mbox{\scriptsize SP}}) \approx \frac{x}{2}(x-1)\,\bE_{n-1} + (1-x^2)\,\bE_n + \frac{x}{2}(x+1)\,\bE_{n+1}.
\end{equation}

%---------------------------------------------------------------------------
\subsubsection{Switching point detection via linear extrapolation.}
%\label{exSP}
If $\bE(t)$ is not differentiable in $t_{\mbox{\scriptsize SP}}\in (t_n, t_{n+1}]$, an approximation of the true strain path
via linear interpolation is inaccurate and quadratic interpolation is no improvement. Then, a piecewise bi-linear approximation
 is a better choice to account for the kink in the strain path. For the approximation of the strains in $[t_n,t_{n+x})$ (i.e. before the
 SP $t_{\mbox{\scriptsize SP}}=t_{n+x}$) the strain values at $t_{n-1}$ and at $t_n$ are linearly extrapolated:
\begin{equation}
\bE(t)\approx \bE_n + \frac{t-t_n}{\Delta t}(\bE_n-\bE_{n-1})\,.
\end{equation}
The equation to determine parameter $x$ with $t_{\mbox{\scriptsize SP}}\approx t_n + x\Delta t$ then reads
\begin{equation}
2\mu ||(\bE^{\mbox{\scriptsize D}}_n+x(\bE^{\mbox{\scriptsize D}}_n-\bE^{\mbox{\scriptsize D}}_{n-1})) -\bE^{\mbox{\scriptsize p}}_n || - \sqrt{\frac{2}{3}}\left(\sigma_Y + \dK(\alpha_n)\right) = 0\,,
\end{equation}
and the strains at the SP are approximated by
\begin{equation}
\bE(t_{\mbox{\scriptsize SP}})\approx \bE_{\mbox{\scriptsize SP}}=\bE_n + x(\bE_n - \bE_{n-1}).
\end{equation}

If, for the approximation of strains in the plastic phase $[t_{\mbox{\scriptsize SP}},t_{n+1}]$ only values later than the kink
shall be used, then only data at $\bE_{\mbox{\scriptsize SP}}$ and $\bE_{n+1}$ are available, thus 

\begin{equation}
\bE(t) \approx  \bE_{\mbox{\scriptsize SP}} + \frac{t-(t_n+x\Delta t)}{(1-x)\Delta t}(\bE^{\mbox{\scriptsize}}_{n+1} - \bE_{\mbox{\scriptsize SP}}) \,.
\end{equation}

Since $\bE_{\mbox{\scriptsize SP}}$ cannot be determined exactly, it is not reasonable to use SP data for interpolation in the consecutive
time step. Then, the SP is within the time interval $[t_{n-1},t_n]$. If no values earlier than the switching
point shall be used, linear interpolation is the method of choice. In the consecutive steps, quadratic interpolation is to be used.

\begin{Figure}[htbp]
\center
   \begin{minipage}{15cm}
   \centering
    \includegraphics[width=0.60\textwidth, angle=0, clip=]{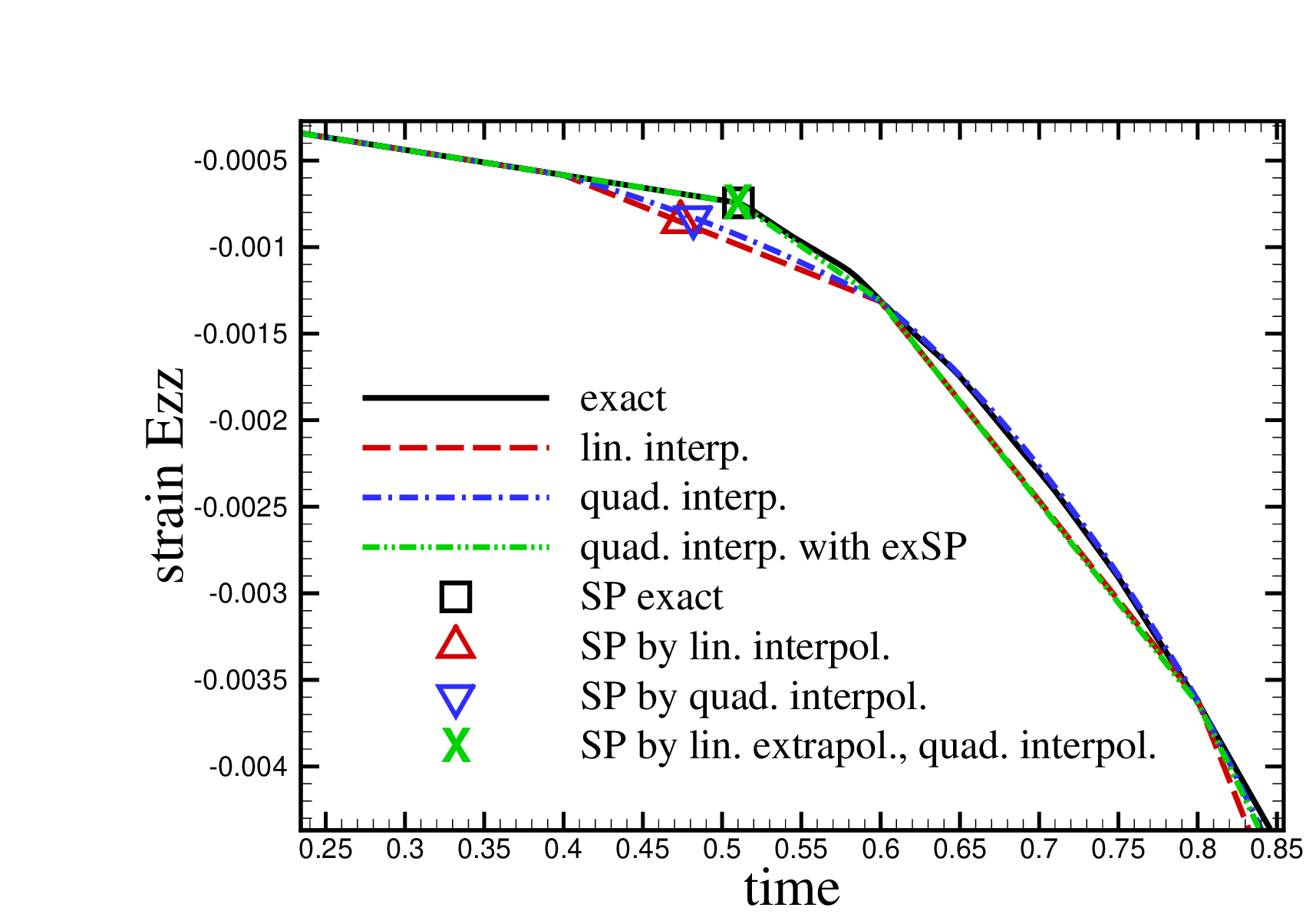}
   \end{minipage}
\caption{Strain component $E_{zz}$ in time: comparison of different methods for SP detection along with strain interpolation }
\label{fig:VergleichSP}
\end{Figure}
 
In Fig.~\ref{fig:VergleichSP} the three different variants for strain interpolation along with SP detection are compared. The considered 
example is a homogeneous, biaxial tension-compression test. In the diagram, the strain component in the direction of compression, here $E_{zz}$, 
is displayed. It can be seen that the SP detection via linear extrapolation is quite accurate, if the strain path before the SP is 
linear.

\subsection{Solution}
\label{subsec:Solution}
Now, having determined the stage values of $\bE^{\mbox{\scriptsize D}}_{ni}$ by (linear or quadratic) interpolation, the set of nonlinear 
equations \eqref{nonlinEqSys} is reformulated according to 
\begin{equation}
\label{Solution}
\left.
\begin{array}{rcl}
\bE^{\mbox{\scriptsize p}}_{ni} &=& \bE^{\mbox{\scriptsize p}}_n + \Sum_{j=1}^s a_{ij} \Delta\Gamma_{nj} \dfrac{\bE^{\mbox{\scriptsize D}}_{nj} - \bE^{\mbox{\scriptsize p}}_{nj}}{||\bE^{\mbox{\scriptsize D}}_{nj} - \bE^{\mbox{\scriptsize p}}_{nj} ||} \\
\Lambda_{ni} &=& \alpha_n + \Sum_{j=1}^s a_{ij} \sqrt{\dfrac{2}{3}} \Delta\Gamma_{nj}                                         \\
0&=&2\mu ||\bE^{\mbox{\scriptsize D}}_{ni} - \bE^{\mbox{\scriptsize p}}_{ni} || - \sqrt{\dfrac{2}{3}}\left(\sigma_Y + \dK(\Lambda_{ni})\right)
\end{array}\right\} i = 1,...,s \, .
\end{equation}
Equation \eqref{Solution}$_2$ is linear and can directly be inserted into \eqref{Solution}$_3$. Then, the nonlinear set
of equations \eqref{Solution}$_1$ and \eqref{Solution}$_3$ have to be solved for the stage values $(\bE^{\mbox{\scriptsize p}}_{ni},\Delta\Gamma_{ni})$, $(i = 1,...,s)$.
Based on this, $\Lambda_{ni}$ is calculated according to \eqref{Solution}$_2$. The number of unknowns in the nonlinear set of equations amounts to $(6+1)s=7s$.
For the class of stiffly accurate RK methods, the update at $t_{n+1}$ coincides with the last stage solution $s$, 
hence 
%\begin{equation}
$\bE^{\mbox{\scriptsize p}}_{n+1} = \bE^{\mbox{\scriptsize p}}_{ns}, 
\alpha_{n+1} = \Lambda_{ns},   
\Delta\gamma_{n+1} = \Delta\Gamma_{ns}$.
%\end{equation}
The values of $\bE^{\mbox{\scriptsize p}}_{n+1}$ and $\alpha_{n+1}$ are saved as history variables for the following time step, whereas $\Delta\gamma_{n+1}$ 
is not required for the next time step.

\section{Augmented consistency conditions}
\label{sec:ConsistentAlgorithmicTangent-El}
 
\subsection{Consistent tangent moduli}
\label{subsec:ConsistentAlgorithmicTangent-El}
In order to obtain quadratic convergence order of the Newton algorithm, the tangent moduli must be algorithmically 
consistent, i.e. consistent with the time integration method, 
$\mathbb{C}_{n+1}^{\mbox{\scriptsize ep}} = \partial{\bS_{n+1}}/{\partial \bE_{n+1}}$. 
In the context of the present work, the moduli newly depend on the chosen approximation of strain interpolation. 
With the 4th order elasticity tensor $\mathbb{C}$, the consistent tangent can be written as
\begin{equation}
\mathbb{C}_{n+1}^{\mbox{\scriptsize ep}} = \parziell{\bS_{n+1}}{\bE_{n+1}}= \mathbb{C}- 2\mu\parziell{\bE^{\mbox{\scriptsize p}}_{n+1}}{\bE_{n+1}} \, .
%\label{eq:Elpl-Tangent-General}
\end{equation}
Differentiation of \eqref{Solution}$_1$ and \eqref{Solution}$_3$ with respect to $\bE_{n+1}$ yields a linear set of equations for calculating
$\dfrac{\partial\bE^{\mbox{\scriptsize p}}_{n+1}}{\partial\bE_{n+1}} = \dfrac{\partial\bE^{\mbox{\scriptsize p}}_{ns}}{\partial\bE_{n+1}}$.
The explicit derivations are carried out in Appendix B. 

The concept of the algorithmic tangent, that is consistent with the time integration algorithm, was introduced by \cite{Nagtegaal1982}
%\footnote{For a discussion on the application of the Newton-Raphson method in nonlinear finite element analyses including earlier works we refer to \cite{Hartmann2005}.} 
for elasto-plasticity with linear isotropic hardening and generalized to other hardening laws and non-associative elasto-plasticity 
by \cite{SimoTaylor}. Since then, it is a very standard in computational inelasticity and textbook knowledge. 
\cite{Miehe1996} presents a general scheme to calculate algorithmic (consistent) tangent moduli by numerical means in terms of difference formulae.
Since this material-independent, simple method overcomes the necessity of algebraic expressions (along with tedious and error-prone analytical manipulations)
it is a milestone in computational inelasticity. 

Since the algorithmic tangent as well as stress are the result of time integration on quadrature-point level, and which are then passed over to 
the global Newton-Raphson iterations, they can be understood as a \emph{local-global} -- or equally -- \emph{time-space} coupling-agency, 
see \cite{Eidel-Habil}. This coupling in terms of $\mathbb{C}_{n+1}^{\mbox{\scriptsize ep}}$ and $\bS_{n+1}$ along with its particular 
direction of information passing is visualized in Fig.~\ref{fig:Spatio-Temporal-Coupling-Augmented}. 

\subsection{Novel conditions}
\label{subsec:NovelConsistency}

\begin{Figure}[htbp]
%\fbox{\parbox{15.5cm}{\hspace*{2mm}\parbox{15.4cm}{
 \begin{minipage}{15.4cm}
\fbox{\parbox{15.5cm}{\hspace*{2mm}\parbox{15.4cm}{
\centering{\textsc{Space-Time Coupling in Elasto-Plastic Finite Element Algorithms}}
\\[4mm]
    \begin{minipage}{2.9cm}
       \vspace*{14mm}
       variational form \\
       of the BVP  \\
       -- \emph{global level} --
       \vspace*{12mm}\\
       {\bf \color{blue}space-time}  \\
       {\bf \color{blue}coupling} \\
       \vspace*{2cm} \\
       IVP as a DAE \\
       -- \emph{local level} --
    \end{minipage}
    \begin{minipage}{8.0cm}
  %     \scalebox{.9}{\input{bild1202a.pstex_t}}
    \includegraphics[width=1.\textwidth, angle=0, clip=]{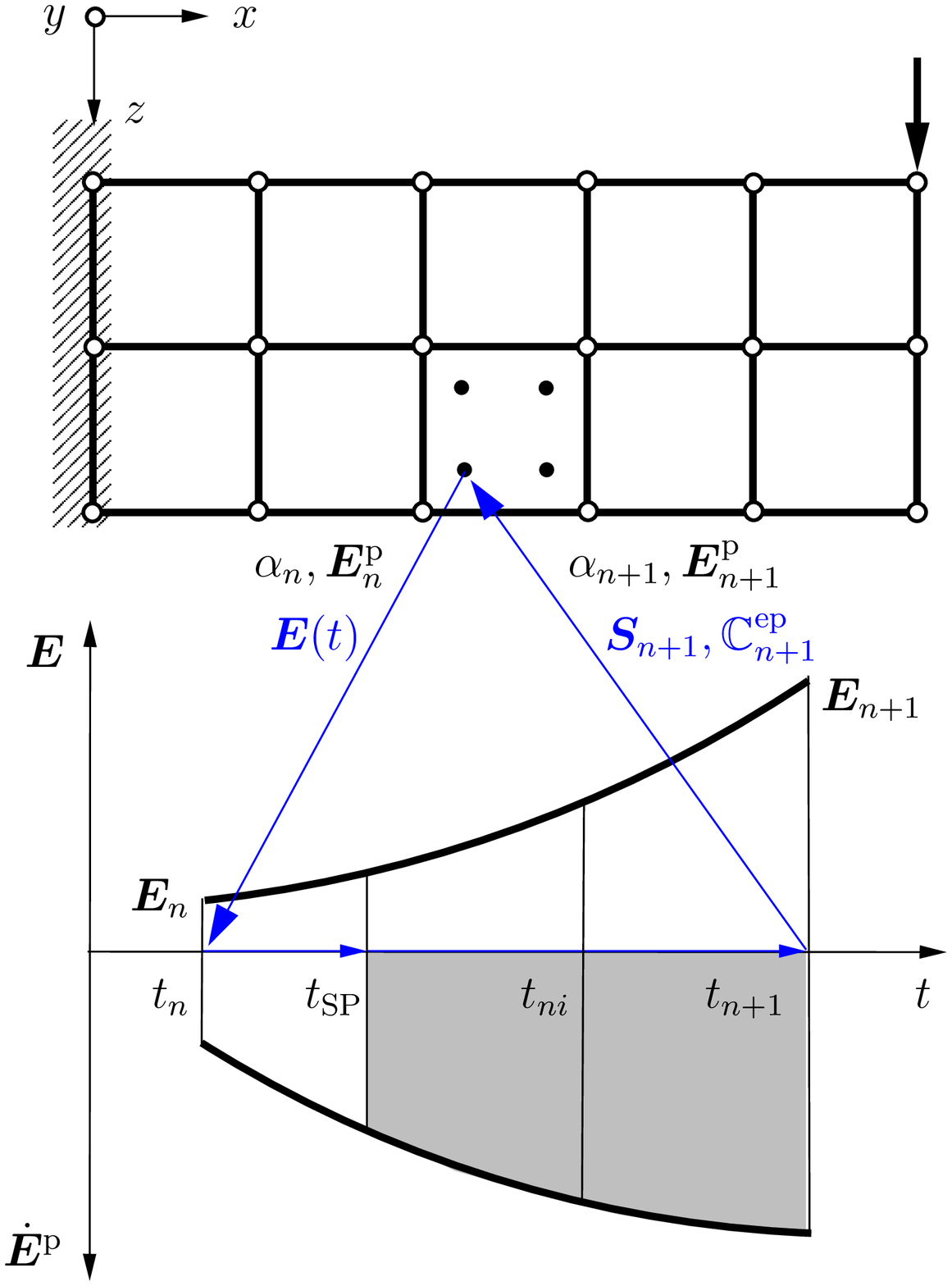}
    \end{minipage}
    \begin{minipage}{3.5cm}
       \vspace*{4.0cm}
       {\bf {\color{red}consistency} \\
       of coupling:}   \\
%       \vspace*{0cm} \\
       \qquad $${\color{red}\bullet} \quad \mathbb{C}_{n+1}^{\mbox{\scriptsize ep}}= \dfrac{\partial \bS_{n+1}}{\partial \bE_{n+1}}$$
       \vspace*{3mm}\\
       \fbox{\parbox{3.4cm}{\hspace*{0mm}\parbox{3.3cm}{
       %\vspace*{-3mm}
       ${\color{red}\bullet}$ \,\, $q=p$  \quad \mbox{for} \\[4mm]
       \hspace*{3mm} \begin{array}[t]{l}
                     \bE = \mathcal{O}(\Delta t^q) \\
                     \bE^{\mbox{\scriptsize p}} = \mathcal{O}(\Delta t^p)
        \end{array}
       }}}
       \vspace*{3mm}\\
       \fbox{\parbox{3.4cm}{\hspace*{0mm}\parbox{3.3cm}{
       %\vspace*{-3mm}
       ${\color{red}\bullet}$ \,\, $f_{n+1}^{\mbox{\scriptsize tr}}(t_{\mbox{\scriptsize SP}})=0$  \\[4mm]
       %\begin{array}[t]{rl}
       $t_{\mbox{\scriptsize SP}} \in (t_n, t_{n+1}]$ for\\[2mm]
       $f^{\mbox{\scriptsize tr}}_n<0, \,\,f^{\mbox{\scriptsize tr}}_{n+1}\geq0$
        %\end{array}
       }}}
    \end{minipage}
  \vspace*{-2mm}
  \begin{eqnarray*}
     {\bE^{\mbox{\scriptsize p}}_{n+1}} &=& \bE^{\mbox{\scriptsize p}}_n + \int_{t_n/t_{\mbox{\tiny SP}}}^{t_{n+1}} \dot{\bE^{\mbox{\scriptsize p}}}(\bE, \bE^{\mbox{\scriptsize p}}) \, \mbox{d}\tau \\
      \alpha_{n+1} &=& \alpha_n + \int_{t_n/t_{\mbox{\tiny SP}}}^{t_{n+1}} \dot{\alpha} \, \mbox{d}\tau \\
                0  &=& f(\bS_{n+1})
  \end{eqnarray*}
  }}}
  \end{minipage}
\caption{Augmented scheme for the consistent coupling of spatio-temporal discretizations in elasto-plastic finite element algorithms for higher order 
         methods $p\geq 2$, \cite{Eidel-Habil}.\label{fig:Spatio-Temporal-Coupling-Augmented}}
%}}}
\end{Figure}

If we follow this coupling perspective in the opposite direction, namely from the global-local or space-time coupling, an aspect 
shows up, which has been overseen so far, but which is crucial for the order barrier issue. The coupling agency in the opposite 
direction is the total strain tensor, here $\bE(t)$, see Fig.~\ref{fig:Spatio-Temporal-Coupling-Augmented}. The approximation of 
total strains by interpolation polynomials was addressed in Sec.~\ref{subsec:DehnInterpol}. The correct, i.e. consistent format 
of interpolation follows from a consistency analysis of the space-time coupling. The DAE 
\begin{eqnarray}
     \dot{\bE}^{\mbox{\scriptsize p}} &=&  \bbf\left(\bE(t), \bE^{\mbox{\scriptsize p}}(t)\right) \label{eq:DAE-DependencyFromE} \\
              0  &=& f(\bS_{n+1}) 
\end{eqnarray}
describing plastic flow is solved by a time integration algorithm of the nominal order $p$, 
\begin{equation}
\label{eq:Order-of-Ep-is-p}
     \bE^{\mbox{\scriptsize p}} =  \mathcal{O}(\Delta t^{p}) \,.
\end{equation}
Since the time-dependent, total strain tensor $\bE(t)$ serves as argument in
\eqref{eq:DAE-DependencyFromE}, its (polynomial) approximation order in time 
will influence the order of time integration. Remarkably, this dependency, though obvious, 
has been neglected in previous works with the exception of \cite{EidelKuhn2010}.
An interpolation polynomial of degree $q-1$ exhibits the approximation order $q$ 
\begin{equation}
\label{eq:Order-of-E-is-q}
     \bE =  \mathcal{O}(\Delta t^q) \, ,
\end{equation}
which implies that $q$ is due to \eqref{eq:DAE-DependencyFromE} an upper bound for the convergence order in time integration. 
Combining \eqref{eq:DAE-DependencyFromE} with \eqref{eq:Order-of-Ep-is-p} and with \eqref{eq:Order-of-E-is-q}, it follows, cf. \cite{Eidel-Habil}

\begin{minipage}{15.5cm}
\fbox{\parbox{15.4cm}{\hspace*{1mm}\parbox{15.0cm}{

\begin{equation}
     \bE^{\mbox{\scriptsize p}} =  \mathcal{O}(\Delta t^{\mbox{\scriptsize min}\{p,q\}}) \, .
\label{eq:min-pq}
\end{equation}

}}}
\end{minipage}

Two conclusions are obvious. First, to empower time integration to its full order, the order of strain approximation
and the order of time integration must be consistent\footnote{Of course, for $q>p$, full convergence order $p$ will equally 
be achieved.}, $q=p$. Second, a polynomial of lower order $q=p-m$, $m\geq1$ will reduce the consistency order of the differential 
variable $\bE^{\mbox{\scriptsize p}}$ to $p-m$ and, as a consequence, of the stress tensor $\bS$ as well.
\\[2mm]
Summarizing, the standard picture of \emph{consistency} of space-time coupling in elasto-plasticity has been restricted so far 
to the local-to-global or time-to-space coupling, which is fulfilled by the algorithmically consistent tangent for the sake of quadratic
convergence of Newton's method. 
For full convergence order in time integration however, this picture is incomplete and therefore must be augmented. 
Two novel consistency requirements arise, which are necessary conditions for order $p\geq3$. They are highlighted by boxes 
in Fig.~\ref{fig:Spatio-Temporal-Coupling-Augmented}.
 
\begin{enumerate}
\item[(i)] The novel consistency condition $q=p$ ensures a strain approximation consistent to the order of the integrator.  
%\\[2mm]
\item[(ii)] The condition $f_{n+1}^{\mbox{\scriptsize tr}}(t_{\mbox{\scriptsize SP}})=0$ for $f_{n}^{\mbox{\scriptsize tr}}<0$ 
with $t_{\mbox{\scriptsize SP}} \in (t_n, t_{n+1}]$ detects the SP and leads to its corresponding total strain, which improves 
the consistency of the initial data $\bm E_{\mbox{\scriptsize SP}}$. 
\end{enumerate}

\bigskip

{\bf Remarks.} 
\\[1mm]
(i) Note, that Fig.~\ref{fig:Interpol-LinVsQuadraticVsExact} showing the nonlinear strain path of $\bE(t)$ is part of Fig.~\ref{fig:Spatio-Temporal-Coupling-Augmented}.
    Figure~\ref{fig:Spatio-Temporal-Coupling-Augmented} displays $\dot{\bE}^{\mbox{\scriptsize p}}$ in the time interval $\Delta t$ 
    and contains two distinct cases; the case that plastic flow commences at $t_{SP}$ and the case that plastic flow continues at $t_{n}$. 
    For the first case with a SP in $\Delta t$, the rate $\dot{\bE}^{\mbox{\scriptsize p}}$ is zero for $t \in (t_n, t_{SP})$.
\\[2mm]
(ii) Backward-Euler is included as a special case in the \emph{consistency} scheme of Fig.~\ref{fig:Spatio-Temporal-Coupling-Augmented}. 
     As a linear method, $p=1$, it requires for $q=1$ a constant strain approximation, hence $\bE(t):=\bE_{n+1}$ in the considered time interval. 
\\[2mm]
(iii) For RK methods of order $p=2$ (e.g. a diagonally implicit RK method in two stages), a consistent order 
      of strain approximation $q=2$ requires linear interpolation based on total strain data at $t_n$ and $t_{n+1}$. 
      For $p=4$, cubic polynomials have to be used based on $t_{n+1-i}$ data with $i=0,1,2,3$ data.  
      This generalization to arbitrary order $1\leq p \leq 4$ is suggested in \cite{Eidel-Habil} and is analyzed for 
      finite-strain viscoelasticity in \cite{EidelStumpfSchroeder2012}.
\\[2mm]
%\bigskip
%
%{\color{red}{\bf Final Remark.} 
%\\[1mm]
(iv) It is worth to mention that in the RK stages (except of at the end of the time interval) stress and strain only have to satisfy 
the equations \eqref{Solution} describing the numerical solution of the plastic IVP. Hence, at the aforementioned RK stages strain 
may violate compatibility and stress may violate the (weak form of the) balance of momentum. The violation of meaningful physical 
laws at intermediate stages of a solution process is not an exception to the rule but a common property of numerical methods. 
Of course, at the end of the time interval all relevant mechanical equations must be and are in fact fulfilled; strain is compatible, 
stress does fulfil the variational form of the balance of momentum and it does not violate the yield condition. 

\bigskip

Table \ref{tab:PredictorCorrector} displays the consistent algorithmic embedding of RK methods for $p=3$ into a standard predictor-corrector 
algorithm, which is based on an elastic-plastic operator split. The key new contribution is the SP detection as well as the quadratic strain 
interpolation, which is the result of the presented consistency analysis.   

%-------------------------------- r e t u r n   m a p ---------------------------------------------

%\vfill
%\newpage

\begin{Table}[htbp]
\begin{tabular}{ p{15.5cm}}
\hline\\[-4mm]
\end{tabular}
%\begin{center}
%{\bf Predictor-Corrector for Elasto-Plasticity with Runge-Kutta Methods    \\
%     of Order $p=3$ via Consistent Strain Interpolation of Order $q=p$}    \\[2mm]
%\end{center}
%\begin{tabular}{ p{15.5cm}}
%\hline\\[-2mm]
%\end{tabular}

\begin{enumerate}
\item[{\bf (I)}] Provide $t_n$-history data $\{\bE^{\mbox{\scriptsize p}}_n, \alpha_n\}$ and total strain data $\{\bE_{n-1}, \bE_n, \bE_{n+1}\}$. 

\item[{\bf (II)}] Elastic predictor and switching point (SP) detection:
\begin{eqnarray*}
    \bS^{\mbox{\scriptsize D,tr}}_{n+1} = 2 \mu (\bE_{n+1}-\bE^{\mbox{\scriptsize p}}_{n}) \quad & & \quad
    f^{\mbox{\scriptsize tr}}_{n+1} =  || \bS_{n+1}^{\mbox{\scriptsize D,tr}} || - \sqrt{\dfrac{2}{3}} \, \sigma_Y(\alpha_n) 
\end{eqnarray*}
If \hspace*{4mm} $f^{\mbox{\scriptsize tr}}_{n+1}<0$ \, then \, $\bE^{\mbox{\scriptsize p}}_{n+1}=\bE^{\mbox{\scriptsize p}}_{n}$, $\alpha_{n+1}=\alpha_n$, 
$\bS_{n+1} = p\, \1 + \bS^{\mbox{\scriptsize D,tr}}_{n+1}$, $\mathbb{C}_{n+1}^{\mbox{\scriptsize ep}}= \mathbb{C}^{\mbox{\scriptsize}}$, \quad exit.\\[2mm]
elseif \hspace*{0mm} $f^{\mbox{\scriptsize tr}}_{n}<0$ \hspace*{1mm} then calculate $t_{\mbox{\scriptsize SP}}= t_n + x \Delta t$ with $x$ from:
%\hspace*{-15mm}
\begin{eqnarray*}
% \mbox{(i)}     & & 2\mu ||(\bE^{\mbox{\scriptsize D}}_n+x(\bE^{\mbox{\scriptsize D}}_{n+1}-\bE^{\mbox{\scriptsize D}}_n))-\bE^{\mbox{\scriptsize p}}_n || - \sqrt{\frac{2}{3}}\left(\sigma_Y(\alpha_n)\right) = 0                                                   \\
 \mbox{(i)}    & & 2\mu || [\frac{x}{2}(x-1)\,\bE^{\mbox{\scriptsize D}}_{n-1} + (1-x^2)\,\bE^{\mbox{\scriptsize D}}_n + \frac{x}{2}(x+1)\bE^{\mbox{\scriptsize D}}_{n+1}] -\bE^{\mbox{\scriptsize p}}_n ||- \sqrt{\frac{2}{3}} \, \sigma_Y(\alpha_n) = 0  \\[-2mm]
 \mbox{(ii)}   & & 2\mu ||(\bE^{\mbox{\scriptsize D}}_n+x(\bE^{\mbox{\scriptsize D}}_n-\bE^{\mbox{\scriptsize D}}_{n-1})) -\bE^{\mbox{\scriptsize p}}_n || - \sqrt{\frac{2}{3}} \, \sigma_Y(\alpha_n) = 0 
\end{eqnarray*}
\hspace*{10mm} with (i) for quadratic interpolation, (ii) for linear extrapolation.\\[2mm]
endif

\item[{\bf (III)}] Calculate stage values $\bE_{ni}$ by the polynomial of consistent degree 2 ($=p-1$):
\begin{eqnarray*}
\label{StageSolution1}
\mbox{If SP in $\Delta t$} \hspace*{10mm}  \bE_{ni} &=& \bE_{\mbox{\scriptsize SP}} + \frac{t_{ni}-(t_n+x\Delta t)}{(1-x)\Delta t}(\bE^{\mbox{\scriptsize}}_{n+1} - \bE_{\mbox{\scriptsize SP}}) \\
                             && \mbox{with} \quad \bE_{\mbox{\scriptsize SP}}=\bE_n + x(\bE_n - \bE_{n-1})  \\
\mbox{else}   \hspace*{10mm}                \bE_{ni} &=& \frac{c_i}{2}(c_i-1)\,\bE_{n-1} + (1-c_i^2)\,\bE_n + \frac{c_i}{2}(c_i+1)\,\bE_{n+1}\,
\end{eqnarray*}

\item[{\bf (IV)}] Plastic corrector: solve for $\bE^{\mbox{\scriptsize p}}_{ni}$, $i = 1,...,s$: \\[-5mm]
\begin{eqnarray*}
  \bE^{\mbox{\scriptsize p}}_{ni} &=& \bE^{\mbox{\scriptsize p}}_n + \Sum_{j=1}^s a_{ij} \Delta\Gamma_{nj} \dfrac{\bE^{\mbox{\scriptsize D}}_{nj} - \bE^{\mbox{\scriptsize p}}_{nj}}{||\bE^{\mbox{\scriptsize D}}_{nj} - \bE^{\mbox{\scriptsize p}}_{nj} ||} \\[-2mm]
\Lambda_{ni} &=& \alpha_n + \Sum_{j=1}^s a_{ij} \sqrt{\dfrac{2}{3}} \Delta\Gamma_{nj}                                          \\[-4mm]
           0 &=& 2 \mu ||\bE^{\mbox{\scriptsize D}}_{ni} - \bE^{\mbox{\scriptsize p}}_{ni} || - \sqrt{\dfrac{2}{3}} \, \sigma_Y(\Lambda_{ni})
\end{eqnarray*}

\item[{\bf (V)}] Update of internal variables $\{\bE^{\mbox{\scriptsize p}}_{n+1}, \alpha^{\mbox{\scriptsize}}_{n+1}\}$ and of stress $\bS_{n+1}$:\\[-2mm]
\begin{equation*}
\begin{array}{rcl}
        \bE^{\mbox{\scriptsize p}}_{n+1} = \bE^{\mbox{\scriptsize p}}_{ns}\, , \qquad \alpha^{\mbox{\scriptsize}}_{n+1} = \Lambda_{ns}\, ,  
                                            & &  \quad \bS_{n+1} = \mathbb{C}:\bE_{n+1} - 2\mu\bE^{\mbox{\scriptsize p}}_{n+1}          %\\[2mm]
%                               \bS_{n+1} &=& \mathbb{C}:\bE_{n+1} - 2\mu\bE^{\mbox{\scriptsize p}}_{n+1}         
\end{array}
\end{equation*}

\item[{\bf (VI)}] Compute $\mathbb{C}_{n+1}^{\mbox{\scriptsize ep}}$ with \eqref{cepsyst1} and \eqref{cepsyst2}, see Appendix B.   
\begin{equation*}
\mathbb{C}_{n+1}^{\mbox{\scriptsize ep}} = \parziell{\bS_{n+1}}{\bE_{n+1}}= \mathbb{C}^{\mbox{\scriptsize}}- 2\mu\parziell{\bE^{\mbox{\scriptsize p}}_{n+1}}{\bE_{n+1}} 
\end{equation*}
\end{enumerate}
\begin{tabular}{ p{15.5cm}}
\hline\\[-4mm]
\end{tabular}
 \caption{Predictor-corrector scheme for von-Mises elasto-plasticity using RK methods of order $p=3$ via consistent strain interpolation 
          of order $q=p$ along with SP detection. \label{tab:PredictorCorrector}}
\end{Table}
  
\vfill
\newpage
%---------------------------------------------------------------------------------------------------------------------------------------------------------------------

%=========================================================================================
\section{Numerical assessment}
\label{sec:Elastoplastic-Tests}
%=========================================================================
In this section, we assess the effectivity of the two measurements undertaken to
overcome the order barrier of $p\leq2$; (i) the improved, quadratic approximation of the strain path
in time and (ii) the SP detection. %as introduced in Sec.\ref{sec:Reasons4OrderReduction}.
For that aim the von-Mises plasticity model of Sec.~\ref{sec:ElastoplasticModel} along with the
third order (two-stage) Radau IIa scheme ($A_{ij}, c_i, b_j$ see Appendix A) have been implemented into 
an 8-node hexahedron element within FEAP, a general purpose finite element code \cite{ZienkiewiczTaylor}.
The test sets considered are the following: %(1) simple shear, 
(1) biaxial stretch and (2) the radial contraction of an annulus.

For all test sets an accurate reference solution $\bX^{\mbox{\scriptsize ex}}$ for tensor $\bX$ with $\bX \in \{\bE, \bE^{\mbox{\scriptsize p}}, \bS\}$
is calculated in each plastified Gauss-point by \emph{numerical overkill} using a very small time step size where the accuracy of the results is in the
range of machine precision. Based on this reference solution, a relative, global error for finite time step sizes is calculated
according to
 
\begin{equation}
\label{RelativeTotalError-ElPl}
e(\bX) = 
\dfrac{1}{ N_{\mbox{\scriptsize el}} \cdot N_{\mbox{\scriptsize gauss}}} 
\sum_{i=1}^{N_{\mbox{\scriptsize el}}} 
\sum_{j=1}^{N_{\mbox{\scriptsize gauss}}}
%\sum_{k=1}^{Nstress}
\frac{||\bX^{(ij)}(\Delta t)-\bX^{(ij)\mbox{\scriptsize ex}}||}{||\bX^{(ij)\mbox{\scriptsize ex}}||} \;,
\end{equation}
 
where $\bX(\Delta t)$ is the tensor for a time step size $\Delta t$, $N_{\mbox{\scriptsize el}}$ is the number of elements in
the domain and $N_{\mbox{\scriptsize gauss}}$ is the number of Gauss-points per element which exhibit plastic deformation
in the current time step. 

In the following, the global relative error is displayed as a function of time step size $\Delta t$ in
double logarithmic scaling. For uniform convergence, the mean order of convergence is calculated by means 
of linear regression.

Table \ref{tab:MethodsAbbreviations1} lists the methods and their variants and introduces abbreviations which will be used in
the following.

\begin{Table}[htbp]
\center
\renewcommand{\arraystretch}{1.3}
%\begin{tabular}{llcc}
\begin{tabular}{ll}
\hline
method                                                                          & abbreviation       \\
\hline
$\bullet$ Backward Euler (i.e. Radau IIa, $s=1$)                                & BE                 \\
$\bullet$ Radau IIa, $s=2, 3$, with linear/quadratic strain interpolation       & RIIa-l/q           \\
\hspace*{1.6cm} as ''RIIa-l/q'', with SP detection by lin./quad. interpolation  & RIIa-l/q-SP        \\
\hspace*{1.6cm} as ''RIIa-l/q'', with SP detection by lin. extrapolation        & RIIa-l/q-exSP      \\
%\hspace*{1.6cm} as ''RIIa-l/q'', with initial step                              & RIIa-l/q-i         \\
\hline
\end{tabular}
\caption{Variants of Radau IIa with $s$ stages of theoretical order $p=2s-1$ of the differential variables
and their abbreviations.\label{tab:MethodsAbbreviations1}}

\end{Table}

%\vfill
%\newpage

%----------------------------------------------------------------------------
\subsection{Biaxial stretch}
\label{ElPl-BiaxialStretch}
%----------------------------------------------------------------------------

We study the elasto-plastic deformation of a cube subject to biaxial stretch
 which is discretized with one hexahedron element of side length $L=1$\,mm.
Material parameters are given in Tab.~\ref{tab:ElastoPlasticMaterialParameters}.
Boundary conditions at $Y=L/2$ and at $X=-L/2$ are chosen to avoid rigid body motions and to
ensure a homogeneous, plane-stress deformation state ($S_{zz}=S_{xz}=S_{yz}=0$).
The time-dependent loading is prescribed by displacement control, at $X=L/2$ by $u_{x}(t)=0.0005\, t$
and at $Y=-L/2$ by $u_{y}(t)=-0.001\,t$.

\begin{Figure}[htbp]
\begin{minipage}{14.5cm}
\begin{center}
  \includegraphics[width=.25\textwidth, angle=0, clip=]{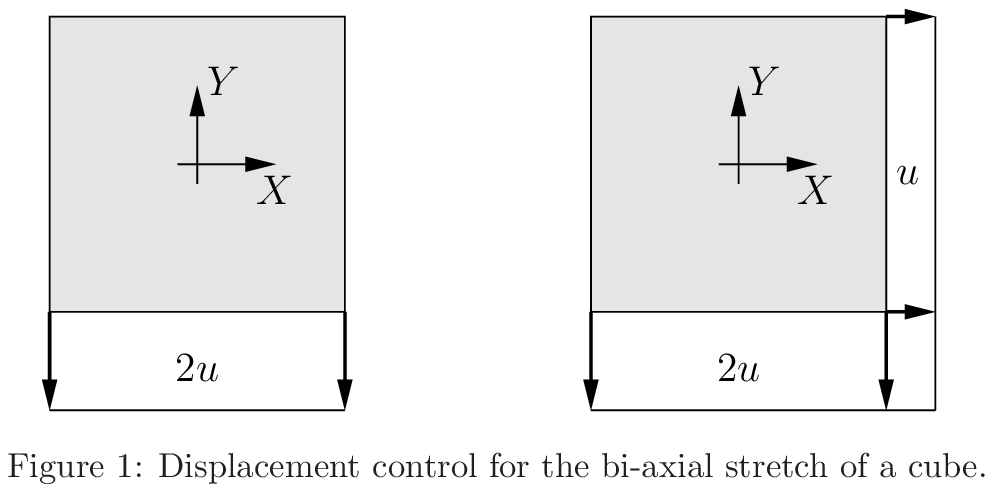}
\end{center}
\end{minipage}
\caption{Bi-axial stretch, mediated by displacement control.}
\label{fig:DisplacementsInBiaxialStretch}
\end{Figure}

\begin{Table}[htbp]
\center
\renewcommand{\arraystretch}{1.4}
%\begin{tabular}{llcc}
\begin{tabular}{cccccc}
\hline
 $E$ $({\text{N}}/{\text{mm}^2})$                            &
 $\nu$ $(1)$                                                 &
 $\sigma_Y$ $({\text{N}}/{\text{mm}^2})$                     &
 $\sigma_\infty - \sigma_Y$ $({\text{N}}/{\text{mm}^2})$     &
 $H$ $({\text{N}}/{\text{mm}^2})$                            &
 $\delta$  $(1)$                                            \\
\hline
$700\,000$   &   $0.2$   &   $875$   &   $211$   &   $1500$   &  $300$\\
 \hline
\end{tabular}
\caption{Biaxial stretch: elasto-plastic material parameters for von-Mises plasticity with nonlinear
isotropic hardening.\label{tab:ElastoPlasticMaterialParameters}}
\end{Table}

The stress state is calculated at $t=1,2,5,10$ using different time step sizes $\Delta t$, the SP is at $t_{\mbox{\scriptsize SP}} = 0.6575$. 
The reference solutions are calculated by Radau IIa employing quadratic interpolation of strain along with extrapolation for 
SP detection (RIIaq-exSP) and with $\Delta t=1.0$E$-04$.

Table~\ref{tab:ConvergenceBiaxialStretch-3-5thorder} and Fig.~\ref{fig:ConvergenceBiaxialStretch-3-5thorder} summarize the results of the convergence study.
For the application of Radau IIa with two stages, the order of convergence is restricted to order two for linear 
interpolation (RIIa-l and RIIa-l-SP) 
and SP detection does not cure order reduction. More accurate results are obtained for quadratic interpolation. However, the two variants, RIIa-q as well as RIIa-q-SP 
show no strictly uniform convergence behavior at $t=1$ and at $t=2$. The improved SP detection via extrapolation along with quadratic interpolation (RIIa-q-exSP) yields 
a more uniform convergence behavior and achieves the theoretical convergence order of 3.

For the 5th order Radau IIa version with $s=3$ stages, the same example is considered. Figure~\ref{fig:ConvergenceBiaxialStretch-3-5thorder} and 
Tab.~\ref{tab:ConvergenceBiaxialStretch-3-5thorder} display the result of the convergence test at various evaluation times. 
The results are similar to the case of the 3rd order Radau IIa version. The observed convergence order is 3 at maximum for RIIa-q-SP and 
RIIa-q-exSP thus falling back behind the theoretical order of $5$. The reason is, that the interpolation error for total strain 
interpolation is of order $\mathcal{O}(\Delta t^3)$, which limits the convergence order of time integration. 
The results underpin again the validity of \eqref{eq:min-pq} for the prediction of the convergence order for the differential 
variables; for $p=5, q=3$ the order is $\mathcal{O}(\Delta t^{\mbox{\scriptsize min}\{p,q\}})= \mathcal{O}(\Delta t^{3})$.
 
\begin{Table}[htbp]
\center
\renewcommand{\arraystretch}{1.4}
\begin{tabular}{rcccccccc}
   \hline
             &  time $t$ \,\,      &   \qquad  & RIIa-l   & RIIa-l-SP & \qquad & RIIa-q & \cellcolor{hellgrau} RIIa-q-SP  & \cellcolor{hellgrau} RIIa-q-exSP   \\
   \hline
$s=2$, $p=3$ &  $\phantom{1}1$\, : &           &  2.05    & 2.07      &        & 2.31   & \cellcolor{hellgrau} 2.19    & \cellcolor{hellgrau} 3.02       \\
             &  $\phantom{1}2$\, : &           &  1.89    & 1.90      &        & 2.20   & \cellcolor{hellgrau} 2.47    & \cellcolor{hellgrau} 3.00       \\
             &  $\phantom{1}5$\, : &           &  2.22    & 2.23      &        & 3.19   & \cellcolor{hellgrau} 3.15    & \cellcolor{hellgrau} 3.24       \\
             &            $10$\, : &           &  2.02    & 2.07      &        & 2.63   & \cellcolor{hellgrau} 2.87    & \cellcolor{hellgrau} 3.01       \\
 \hline
$s=3$, $p=5$ & $\phantom{1}1$\, :  &           &  2.05    &  2.08     &        & 2.21  &  \cellcolor{hellgrau} 2.18    &  \cellcolor{hellgrau} 2.96      \\
             & $\phantom{1}2$\, :  &           &  2.04    &  2.05     &        & 2.02  &  \cellcolor{hellgrau} 2.42    &  \cellcolor{hellgrau} 3.51      \\
             & $\phantom{1}5$\, :  &           &  2.30    &  2.31     &        & 3.12  &  \cellcolor{hellgrau} 3.23    &  \cellcolor{hellgrau} 3.42      \\
             & $          10$\, :  &           &  2.09    &  2.15     &        & 2.60  &  \cellcolor{hellgrau} 2.99    &  \cellcolor{hellgrau} 3.23      \\
 \hline
\end{tabular}
\caption{Biaxial stretch: convergence order of $e(\bS)$ for variants of the 3rd order (upper part) and the 5th order (lower part) Radau IIa method.}
\label{tab:ConvergenceBiaxialStretch-3-5thorder}
\end{Table}

\begin{Figure}[htbp]
   \begin{minipage}{15.5cm}
%   { 
     \centering%
     \tabcolsep1mm%
  \begin{tabular}[]{rcc}
  $e(\bS)$  &  $t=1$  &  $t=10$  \\ 
\raisebox{3.2cm}{$p=3$, $s=2$}
  &
   \includegraphics[width=.40\textwidth, angle=0, clip=]{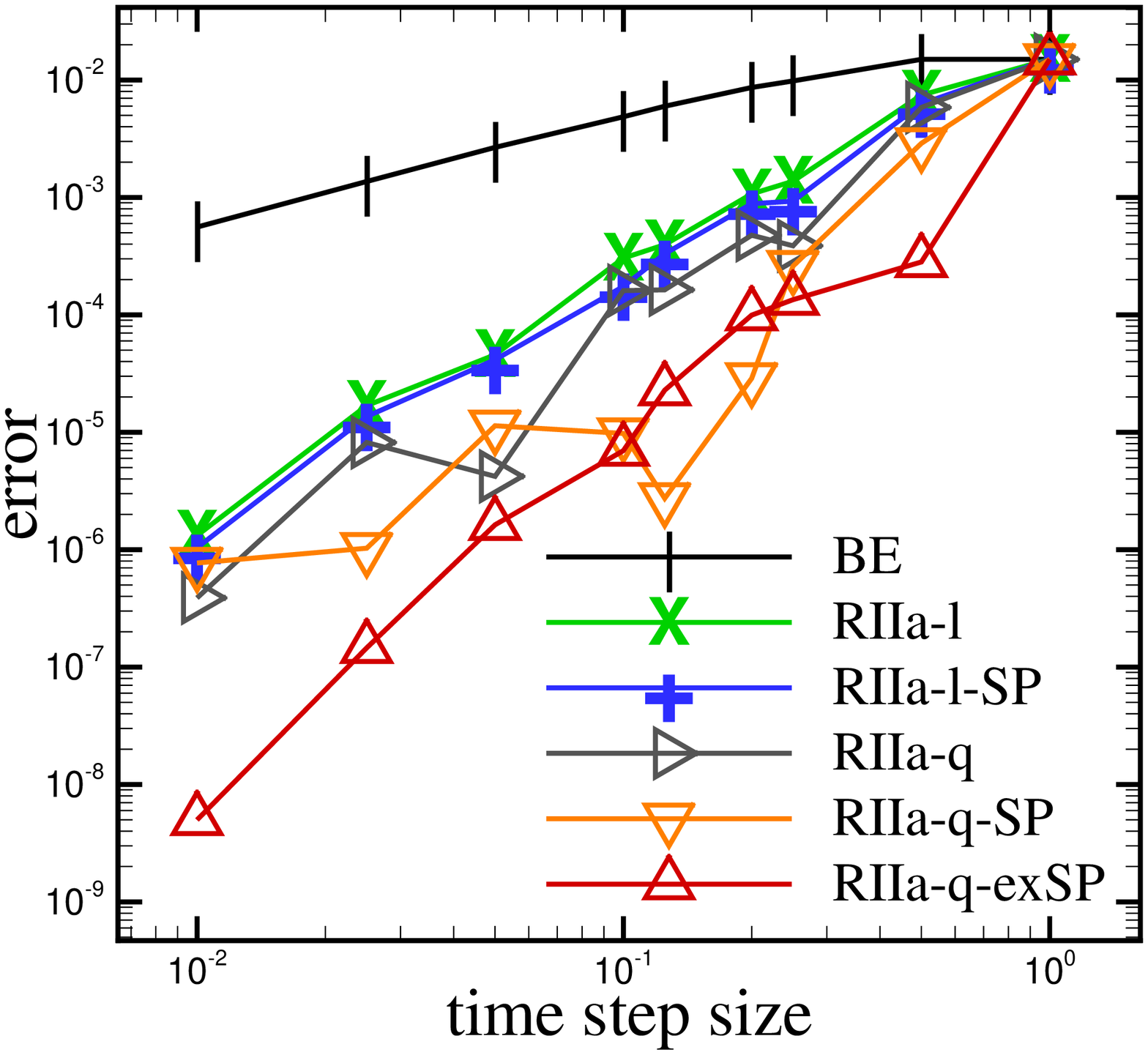} \hspace*{4mm}
  &
   \includegraphics[width=.40\textwidth, angle=0, clip=]{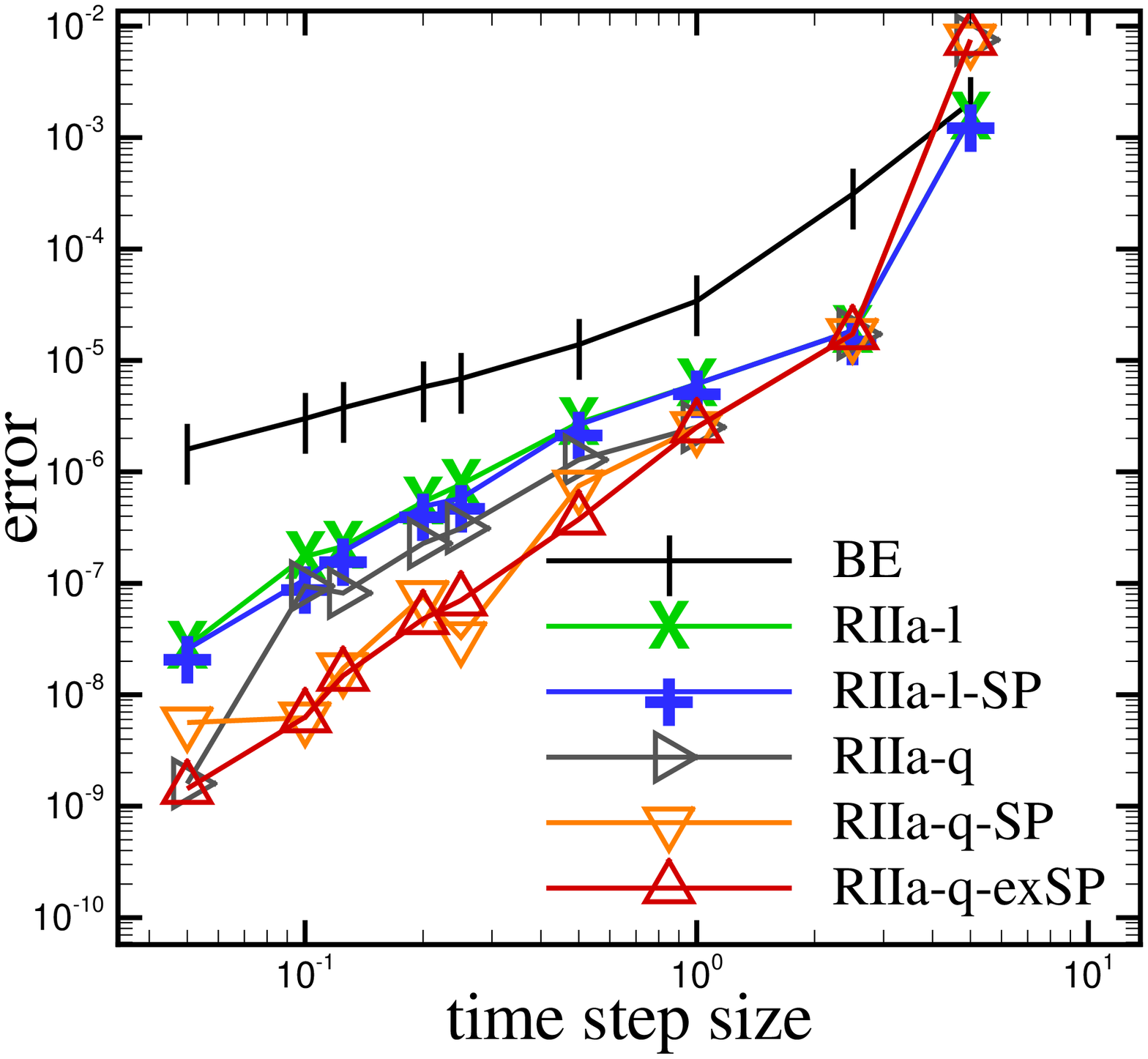} 
   \\
  \raisebox{3.2cm}{$p=5$, $s=3$}
  &
   \includegraphics[width=.40\textwidth, angle=0, clip=]{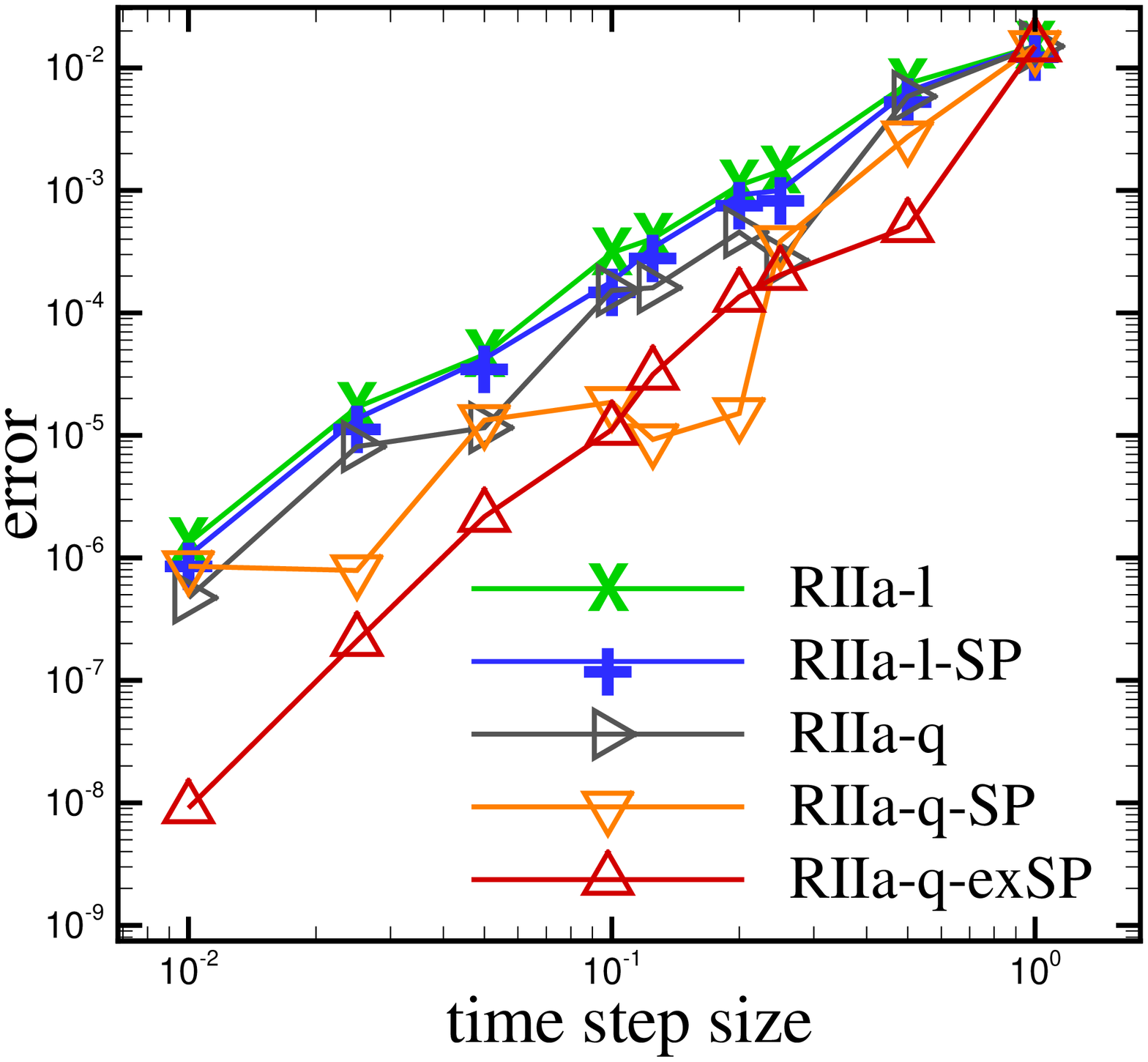}  \hspace*{4mm}
  & 
   \includegraphics[width=.40\textwidth, angle=0, clip=]{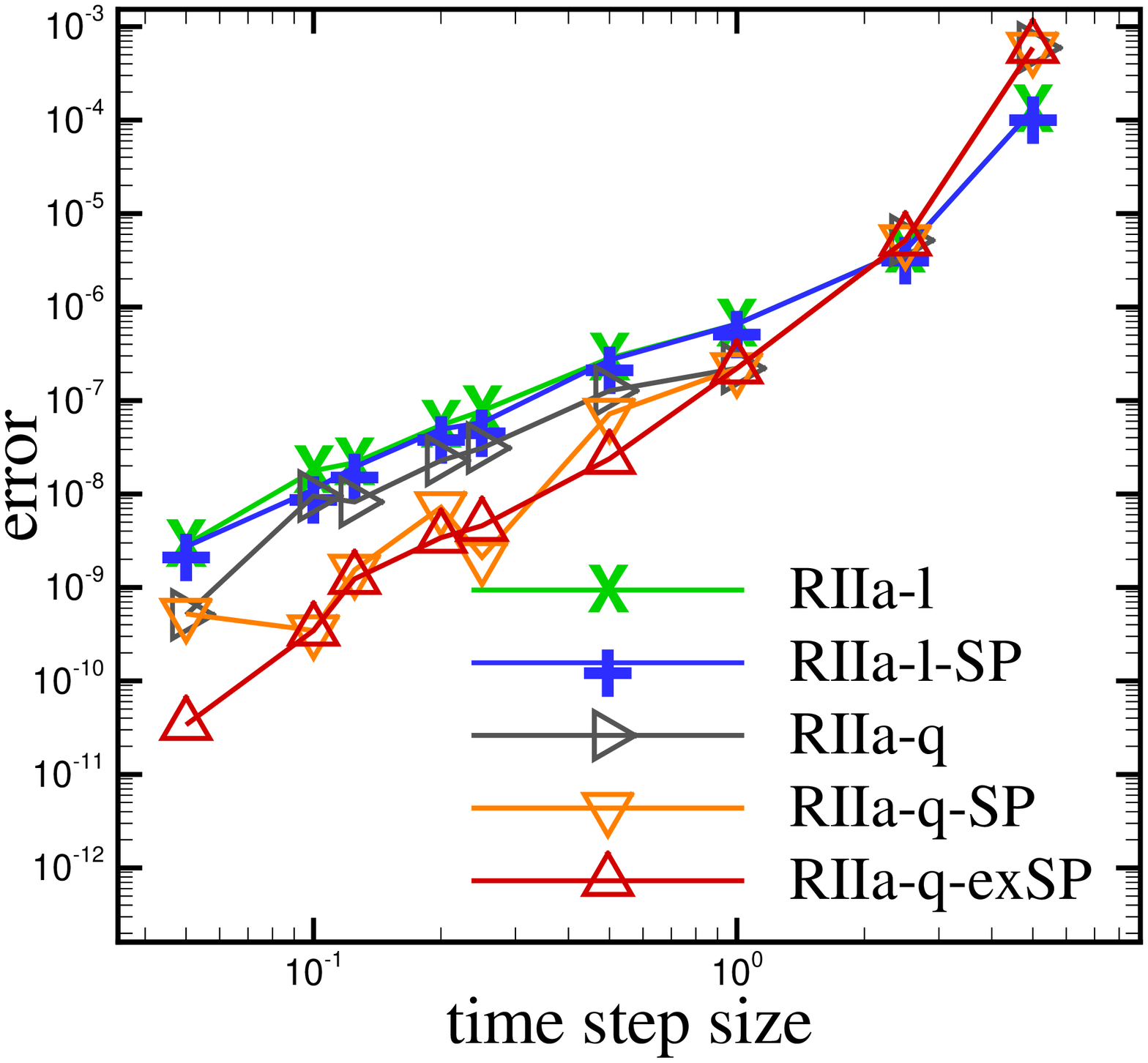}\\
 \end{tabular}
%  } 
  \end{minipage}
  \caption{Biaxial stretch: convergence of $e(\bS)$ for different versions of the 2-stage Radau IIa method
           at $t=1$ (left) and $t=10$ (right). First row: $s=2, p=3$, second row: $s=3, p=5$.
           Note that the two versions show a very similar error pattern, which can be explained by the
           same way of strain interpolation. \label{fig:ConvergenceBiaxialStretch-3-5thorder}}
\end{Figure}

%%%%%%%%%%%%%%%%%%%%%%%%%%%%%%%%%%%%%%%%%%%%%%%%%%%%%%%%%%%%%%%%%%%%%%%%%%%%%%%%%%%%%%%%%%%%%%%%%%%%%%%%%%%%%%%%%%%%%%
 
%\vfill
%\newpage

%---------------------------------------------------------------------------------------------------------------

\subsection{Intermediate considerations: construction of two idealized cases}
In the following, we aim to neatly separate the influence of the two sources for order reduction.
Therefore, two idealized cases are constructed:

\subsubsection{Case I: Discarding switching points by setting $\Bsigma_Y = \mathbf{0}$.}\quad 
                   For zero initial yield stress $\sigma_Y$, the yield surface degenerates into a point, and each material point starts to yield at $t=0$. 
                   Hence, the elastic branch is discarded and no elastic-plastic switching occurs in deformation histories of monotonously increasing load. 
                   This implies that initial conditions are given exactly and do not need to be approximately calculated by SP detection. The 
                   second benefit of discarding the SP (and the necessity to detect it) is that the strain path will be inherently smooth for 
                   monotonously increasing loading. Without SP there is no sudden loss of stiffness manifesting in a kink in the strain path 
                   as in Fig.~\ref{fig:VergleichSP}. Then, quadratic interpolation is expected to become effective and 3rd order convergence shall be obtained. 
                   This ''idealized'' DAE-case of elasto-plasticity therefore comes close to the smooth ODE-case of viscoelasticity, for an analysis of 
                   that latter case see \cite{EidelKuhn2010}. 

For this particular case we consider the two-stage Radau IIa along with linear and quadratic interpolation of the strain path in time.

\begin{Table}[htbp]
\center
\renewcommand{\arraystretch}{1.4}
%\begin{tabular}{llcc}
\begin{tabular}{cccccc}
\hline
 $E$ $({\text{N}}/{\text{mm}^2})$                            &
 $\nu$ $(1)$                                                 &
 $\sigma_Y$ $({\text{N}}/{\text{mm}^2})$                     &
 $\sigma_\infty - \sigma_Y$ $({\text{N}}/{\text{mm}^2})$     &
 $H$ $({\text{N}}/{\text{mm}^2})$                            &
 $\delta$  $(1)$                                            \\
\hline
$700\,000$   &   $0.2$   &   $0$   &   $211$   &   $1500$   &  $300$\\
 \hline
\end{tabular}
\caption{Biaxial stretch, $\sigma_Y=0$: elasto-plastic material parameters for von-Mises plasticity with nonlinear
isotropic hardening.}
\label{tab:ElPlMatParam-Biaxial-sig0=0}
\end{Table}

 The results of the convergence analysis are shown in Tab.~\ref{tab:ConvergenceBiaxialStretch-sigma0=0}
 and in Fig.~\ref{fig:KonvergenzZug2y0}. Linear interpolation results in order $2$ and full convergence order $3$
 is achieved for quadratic interpolation. These observations are in line with the results of previous simulations
 with nonzero yield stress. In contrast to the aforementioned case, convergence is strictly uniform for $\sigma_Y=0$.
 
\begin{Table}[htbp]
\center
\renewcommand{\arraystretch}{1.4}
%\begin{tabular}{llcc}
\begin{tabular}{lccccccc}
\hline
 time $t$   & \hspace*{3mm}  &  BE     &  \hspace*{3mm}  &  RIIa-l    & \hspace*{3mm} &  \cellcolor{hellgrau} RIIa-q    \\
\hline
 $1.5$\, :  &                &  1.05   &                 &  1.72   &               &  \cellcolor{hellgrau} 2.84     \\
 $3.0$\, :  &                &  1.08   &                 &  1.74   &               &  \cellcolor{hellgrau} 2.86     \\
 \hline
\end{tabular}
\caption{Case I: biaxial stretch with $\sigma_Y = 0$, order of convergence for different methods.}
\label{tab:ConvergenceBiaxialStretch-sigma0=0}
\end{Table}

\begin{Figure}[htbp]
\begin{minipage}{15.0cm}
{\hspace*{4.5cm}
  \includegraphics[width=.42\textwidth, angle=0, clip=]{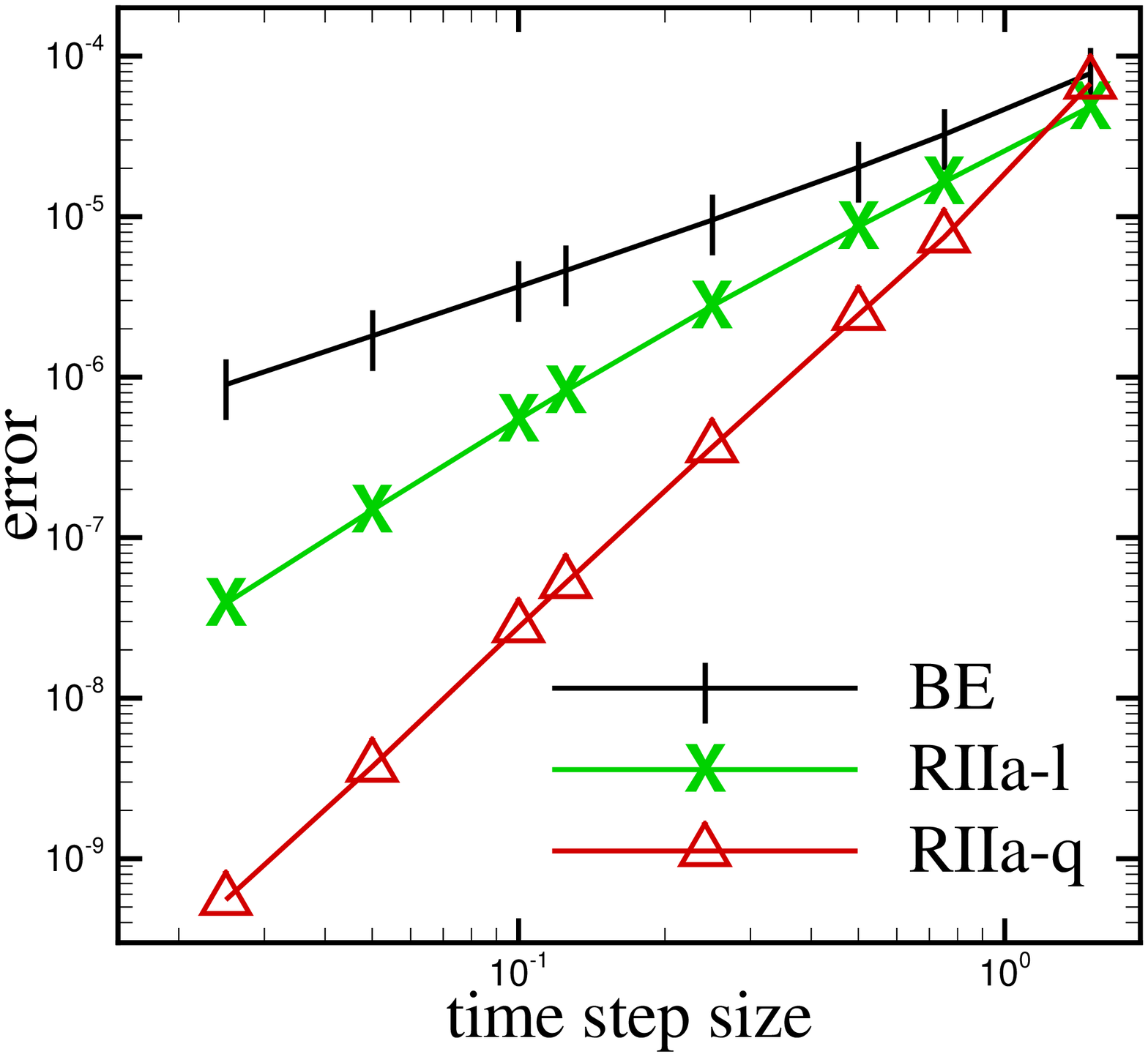}
}
\end{minipage}
 \caption{Case I: biaxial stretch for $\sigma_Y = 0$, comparison of convergence order for $e(\bS)$ evaluated at $t=3$.\label{fig:KonvergenzZug2y0}}
\end{Figure}
  
\subsubsection{Case II: Annihilating the strain interpolation error by direct discretization of the evolution equations.} \quad
                   Calculations so far have been carried out using time integration based on the partitioned ansatz. 
                   Within the partitioned ansatz, stage values of total strains have to be calculated via interpolation and therefore, 
                   a corresponding approximation error is inevitable, see Sec.~\ref{sec:Reasons4OrderReduction}.
                   In the following, we consider the special case, where discretization in time is directly applied to the evolution equations
                   without finite element discretization in space. Of course, this is restricted to homogeneous deformation states with prescribed 
                   displacement boundary conditions. This way, the strain interpolation error is annihilated and the effect of SP detection 
                   onto the convergence order can be analyzed separately. 

\begin{Figure}[htbp]
  \begin{minipage}{15.5cm}
     \centering%
     \tabcolsep1mm%
  \begin{tabular}[]{rcc}
    &  $t=1$  &  $t=10$  \\ 
\raisebox{3.0cm}{$e(E^p_{zz})$}
  &
   \includegraphics[width=.38\textwidth, angle=0, clip=]{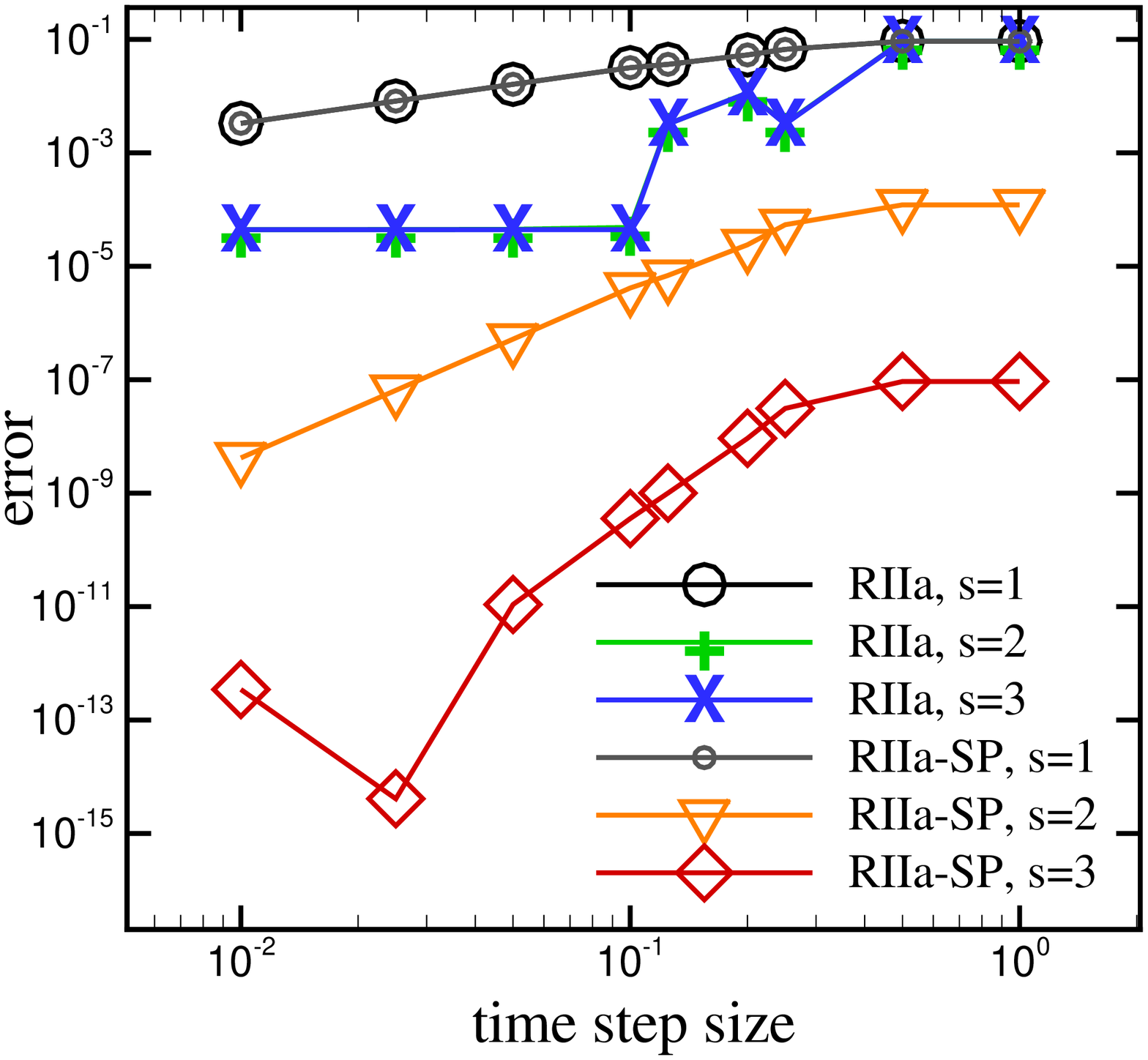} \hspace*{4mm}
  &
   \includegraphics[width=.38\textwidth, angle=0, clip=]{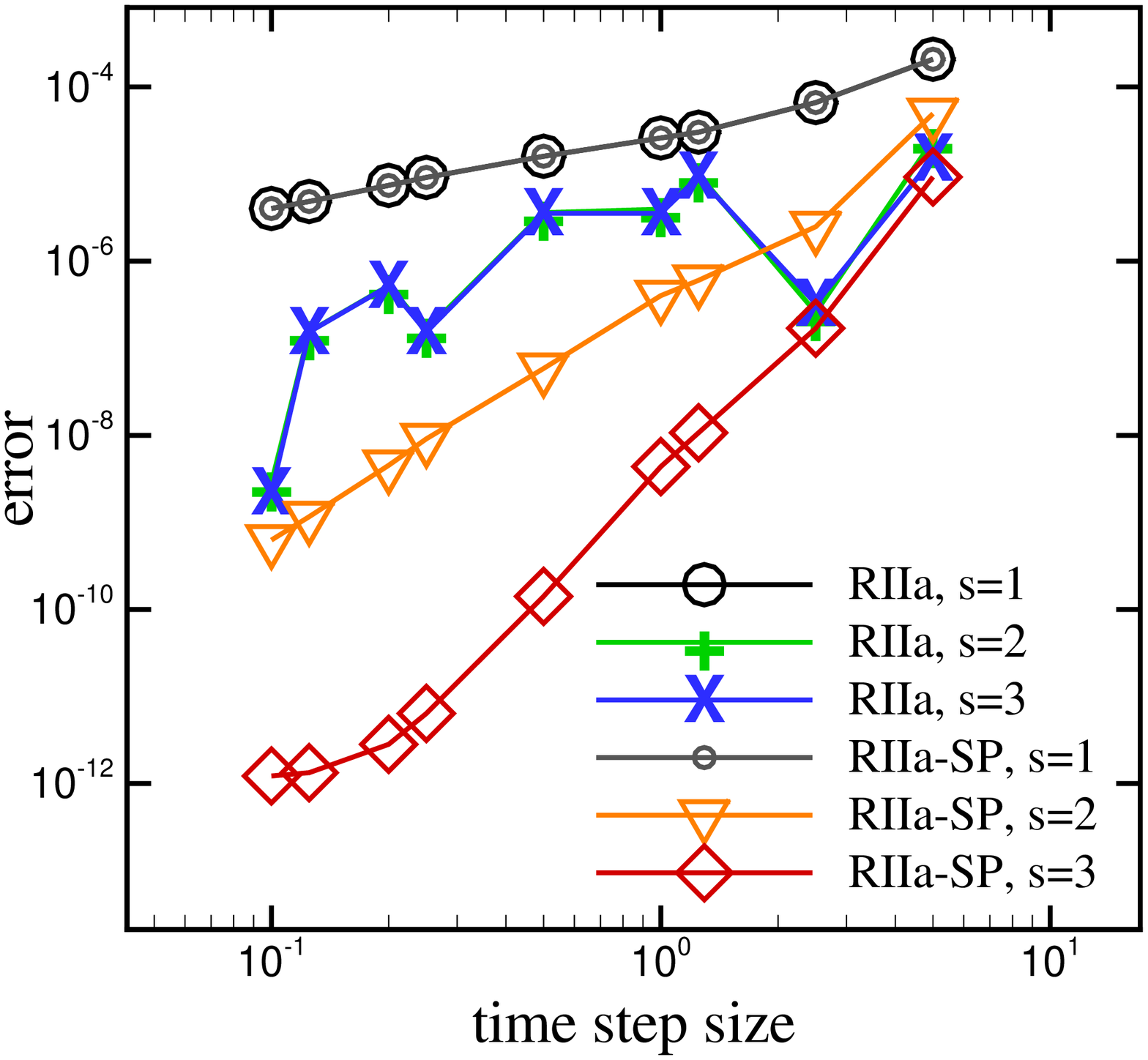}
 \end{tabular}
  \end{minipage}
  \caption{Case II: biaxial stretch using direct discretization of the evolution equations of plastic flow. Error in $E^{\mbox{\scriptsize p}}_{zz}$ versus time step size evaluated 
   at (left) $t=1$ and at (right) $t=10$. Full convergence order (5 for $s=3$ and 3 for $s=2$) is obtained by SP detection.\label{fig:KonvergenzZug2direkt}}
\end{Figure}

\begin{Table}[htbp]
\center
%\begin{tabular}{c|lr|lr|lr|}
\renewcommand{\arraystretch}{1.4}
\begin{tabular}{ccrrcrrcrr}
\hline
          time $t$          &  & $s = 1$:                           & $p=1$ & \,\,  & $s = 2$:                        & \cellcolor{hellgrau} $p=3$  & \,\, & $s = 3$:                                    & \cellcolor{hellgrau} $p=5$   \\
\hline $\phantom{1}1$\, :   &  & {\scriptsize $\Delta t\leq 0.25$:} & 0.98  &       & \scriptsize $\Delta t\leq 0.25$:& \cellcolor{hellgrau} 2.86   &      & {\scriptsize $0.025\leq\Delta t\leq 0.25$:} & \cellcolor{hellgrau} 5.22    \\
       $\phantom{1}2$\, :   &  &                                    & 0.94  &       & \scriptsize $\Delta t\leq 0.25$:& \cellcolor{hellgrau} 2.86   &      & {\scriptsize $\Delta t\geq 0.025$:}         & \cellcolor{hellgrau} 5.08    \\
       $\phantom{1}5$\, :   &  &                                    & 1.36  &       &                                 & \cellcolor{hellgrau} 2.94   &      &                                             & \cellcolor{hellgrau} 5.25    \\
       $10$\, :             &  &                                    & 0.68  &       &                                 & \cellcolor{hellgrau} 2.69   &      &                                             & \cellcolor{hellgrau} 4.65    \\
\hline
\end{tabular}
\caption{Case II: biaxial stretch using direct discretization of the evolution equations. Convergence order of the $s$-stage Radau IIa method {\bf with} switching 
point (SP) detection. The theoretical order of convergence is $p=2s-1$ is obtained for $s=1,2,3$.}\label{tab:RadauKonvergenz}
\end{Table}

We choose elastic Poisson's ratio $\nu=0$ and prescribe strains in $x$-direction and $y$-direction via $E_{xx}(t)$ and
$E_{yy}(t)$. Hence, we may write $\bE = \bE^{\mbox{\scriptsize e}} + \bE^{\mbox{\scriptsize p}} = (E_{xx}(t)\, , \, E_{yy}(t) \, , \, E^{p}_{zz})^T$ 
%\begin{equation}\label{specialeps}
%\bE = \bE^{\mbox{\scriptsize e}} + \bE^{\mbox{\scriptsize p}} =
%\begin{pmatrix}
%E_{xx}(t) \\
%E_{yy}(t) \\
%E^{p}_{zz}
%\end{pmatrix}\,
%\end{equation}
for the reduced set of nonzero strains, since $E_{xy}$, $E_{xz}$ and $E_{yz}$ identically vanish. The DAE-system from the
 evolution equations and the yield condition can be solved %using \eqref{specialeps} 
without solving the BVP.
Hence, the partitioned ansatz is no longer necessary and the rate-equations of plastic flow are directly accessible for
discretization in time. \\\\
For the deviator of the reduced strain tensor the relation holds
$$
\bE^{\mbox{\scriptsize D}} =
\begin{pmatrix}
E_{xx} \\
E_{yy} \\
E^{p}_{zz}
\end{pmatrix} - \dfrac{1}{3} (E_{xx}+E_{yy}+E^{p}_{zz})
\begin{pmatrix}
1\\
1 \\
1
\end{pmatrix} = \frac{1}{3}
\begin{pmatrix}
2\,E_{xx}-\,E_{yy}-E^{p}_{zz} \\
-E_{xx}+2\,E_{yy}-E^{p}_{zz}  \\
-E_{xx}-E_{yy}+2\,E^{p}_{zz}
\end{pmatrix}.
$$
For $\bE^{\mbox{\scriptsize p}}\equiv \mathbf{0}$ before the onset of plastic yielding, the SP can be calculated exactly from the yield condition. 
Here the time of the SP is $t_{\mbox{\scriptsize SP}}=0.693375$.\\\\
In the considered case, stage values of $\bE^{\mbox{\scriptsize D}}$ can be expressed in terms of the stage values $\bE_{ni}^{\mbox{\scriptsize p}}$, 
$E_{xx}(t_{ni})$ and $E_{yy}(t_{ni})$:
$$
\bE^{\mbox{\scriptsize D}}_{ni} = \frac{1}{3}
\begin{pmatrix}
2\,E_{xx}(t_{ni})-E_{yy}(t_{ni})-E_{zz\,ni}^{\mbox{\scriptsize p}} \\
  -E_{xx}(t_{ni})+2\,E_{yy}(t_{ni})-E_{zz\,ni}^{\mbox{\scriptsize p}} \\
  -E_{xx}(t_{ni})-E_{yy}(t_{ni})+2E_{zz\,ni}^{\mbox{\scriptsize p}}
\end{pmatrix}\,.
$$

For the convergence study, we choose the material parameters of Tab.~\ref{tab:ElastoPlasticMaterialParameters} except of $\nu=0$. Strains
in $x$- and $y$-direction are prescribed by $E_{xx}(t) = 0.0005\,t$ and  $E_{yy}(t) = 0.002\,t$.
In the present test we assess different variants of Radau IIa for $s=\{1,2,3\}$ stages, with and without SP detection.
Recall, that the theoretical convergence order $p$ for an $s$-stage method of this RK-family is $p=2s-1$ for the differential variables.

The results for different evaluation times are displayed in Fig. \ref{fig:KonvergenzZug2direkt} and in Tab.~\ref{tab:RadauKonvergenz}.
%The dashed lines in Fig.\ref{fig:KonvergenzZug2direkt} represent convergence orders 1, 3 and 5, respectively. 
For time step sizes, where the simulation
results come close to machine precision, and therefore no longer show uniform convergence, we reduce the data set for linear regression to calculate the mean
convergence order.
\\[2mm]
As can be seen in Fig.~\ref{fig:KonvergenzZug2direkt} and in Tab.~\ref{tab:RadauKonvergenz}, full convergence order is obtained (3 and 5, respectively) 
if SP detection is used.  Without SP detection, Radau IIa suffers from order reduction and shows a strongly non-uniform, even 
non-monotonous convergence; then the 2-stage solution coincides with the 3-stage solution.

\bigskip

{\bf Summary.} \quad 
The results of the cases I and II underpin the hypothesis, that order reduction in elasto-plasticity within a partitioned finite element ansatz 
is due to two reasons, (i) inaccurate strain values at RK-stages, which are obtained by linear interpolation bounding the convergence order to 
order 2 and (ii) a missing or inaccurate SP calculation giving inaccurate initial data. 
Complementary, quadratic interpolation of total strain implying approximation order 3 enables convergence order 3 in time integration, SP 
detection provides consistent initial data and thereby enables full convergence order. 

\bigskip

Even without the highly idealized settings of case I and case II, the examples so far have been homogeneous deformation states and are therefore 
of limited use for reliable conclusions concerning the quality of the time integration algorithms. This is true, since in that case each and 
every Gauss-point passes the elasto-plastic SP at the same time. As a consequence, no stress-redistribution will occur during 
plastic deformations, which is, however, a very characteristic of elasto-plasticity. Its reason is the sudden loss of stiffness at a material 
point when yielding commences and the \emph{''flexibility''} of the material to instantaneously redistribute the loading. How this \emph{''structural''} 
effect of stress-redistribution affects the characteristics of the strain path and the SP detection shall be studied in the next example
showing non-homogeneous deformation states.

%----------------------------------------------------------------------------
\subsection{Radial contraction of an annulus}
\label{ElPl-Annulus}
%----------------------------------------------------------------------------

\begin{Figure}[htbp]
\begin{minipage}{15.5cm}
   \centering
 %     \scalebox{.55}{\input{bild1202c.pstex_t}}
 %     \scalebox{.55}{\input{bild1202c.pstex_t}}
       \includegraphics[width=.42\textwidth, angle=0, clip=]{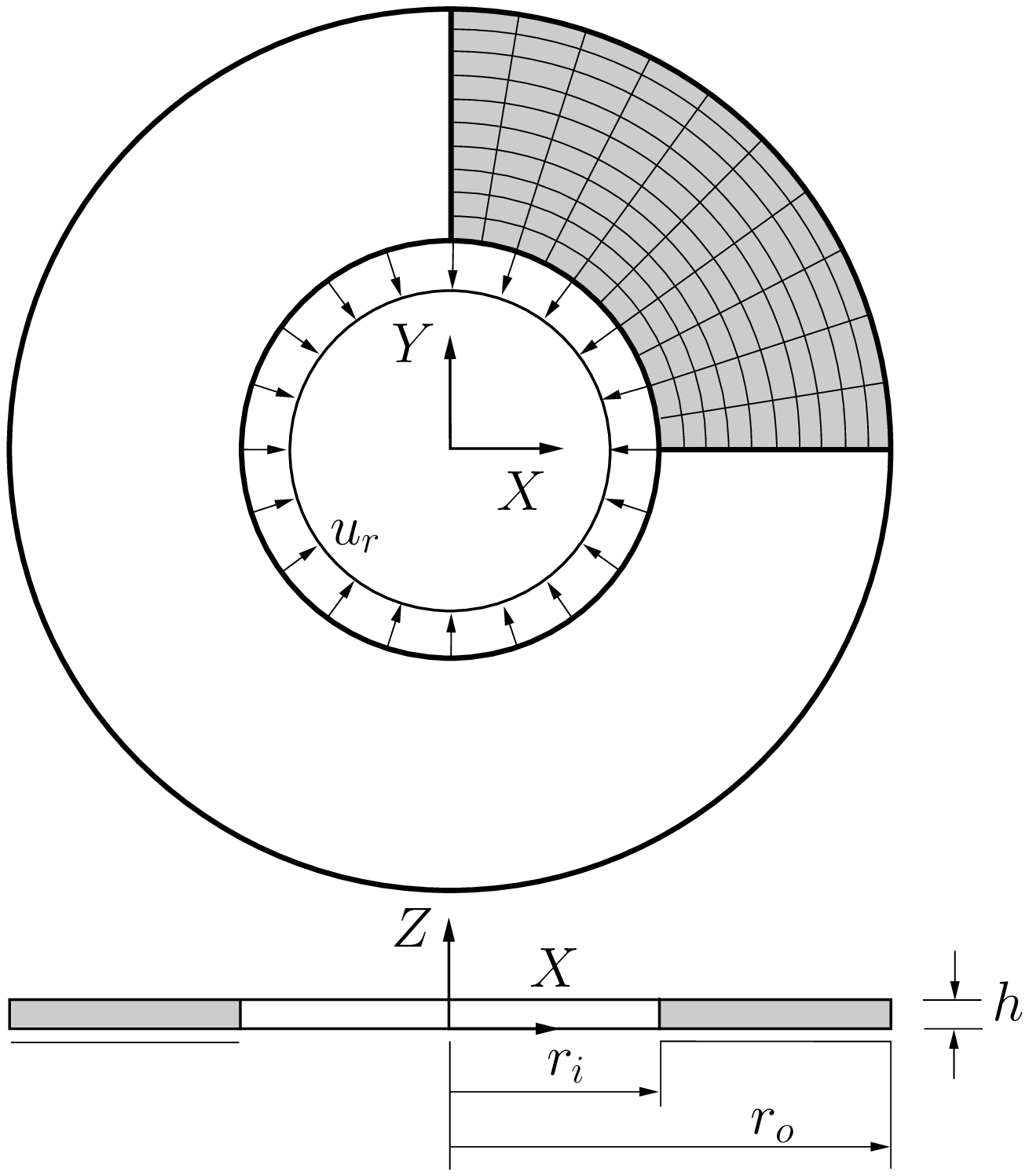}
 
\end{minipage}
\caption{Radial contraction of an annulus by displacement control applied to the inner rim.
Geometry, loading and finite element discretization of one quarter.}
\label{fig:Annulus-pic}
\end{Figure}

This test set is the radial contraction of an annulus as displayed in Fig.~\ref{fig:Annulus-pic}. The annulus 
exhibits inner and outer radii $r_i = 20\,$mm, $r_o = 40\,$mm and thickness $h = 1$\,mm.  
{\color{black}It is supported at $z=0\,$mm in $z$-direction}, and two symmetry planes are exploited, such that the simulation 
is carried out at a quarter system. It is discretized with 10 elements in circumferential direction, 10 elements in radial 
direction and one element over the thickness.
The inner rim is continuously pulled in radial direction by displacement control $u_r(t)=$ 1.0 mm$\,\cdot\,t$.
As a consequence, plasticity spreads out radially from the inner rim to the outer rim.
 
\begin{Figure}[htbp]
  \begin{minipage}{15.5cm}
      \centering%
      \tabcolsep1mm%
  \begin{tabular}[]{rcc}
     &
     \includegraphics[width=0.40\textwidth, angle=0, clip=]{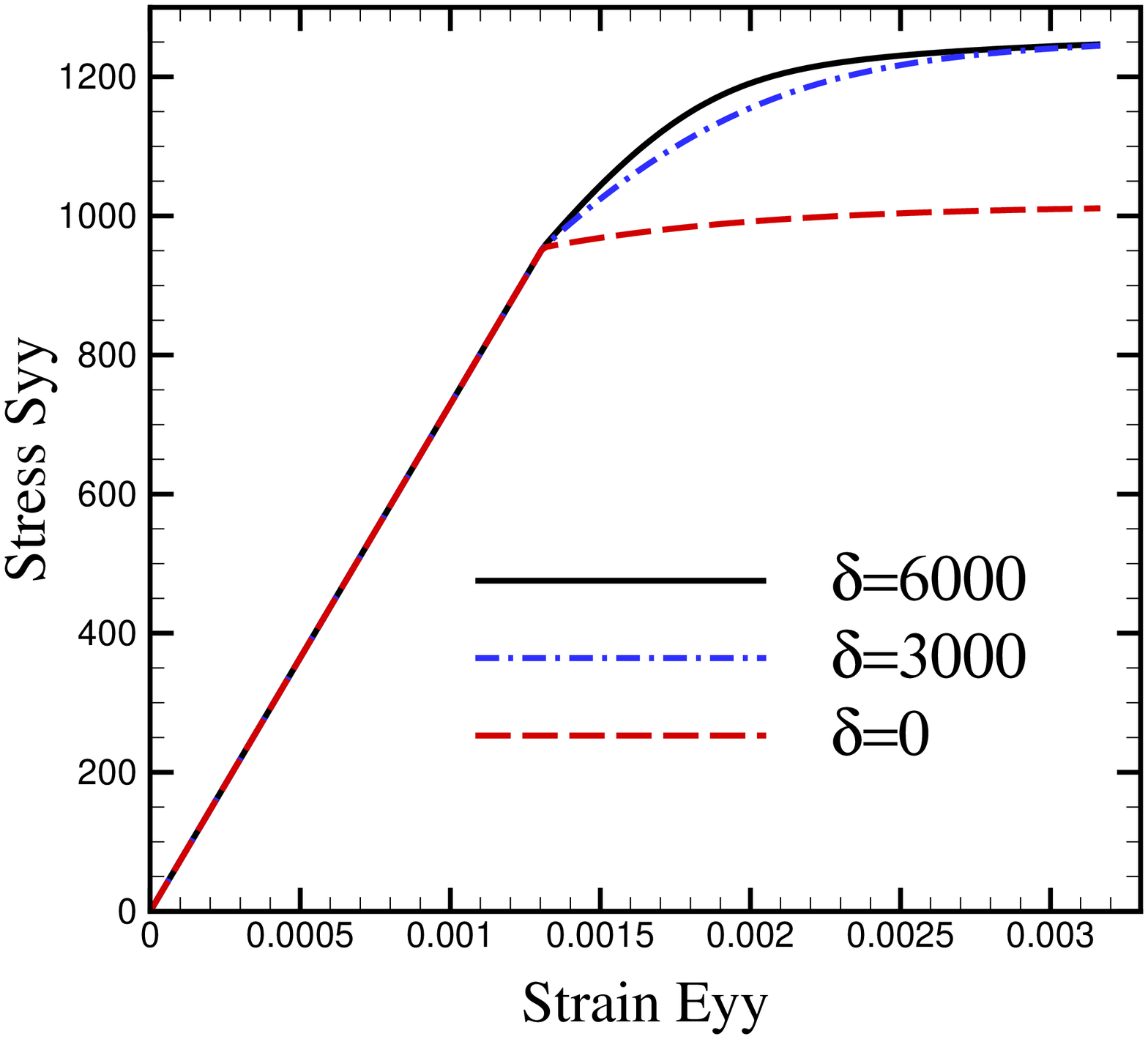} \hspace{5mm}
     & 
     \includegraphics[width=0.40\textwidth, angle=0, clip=]{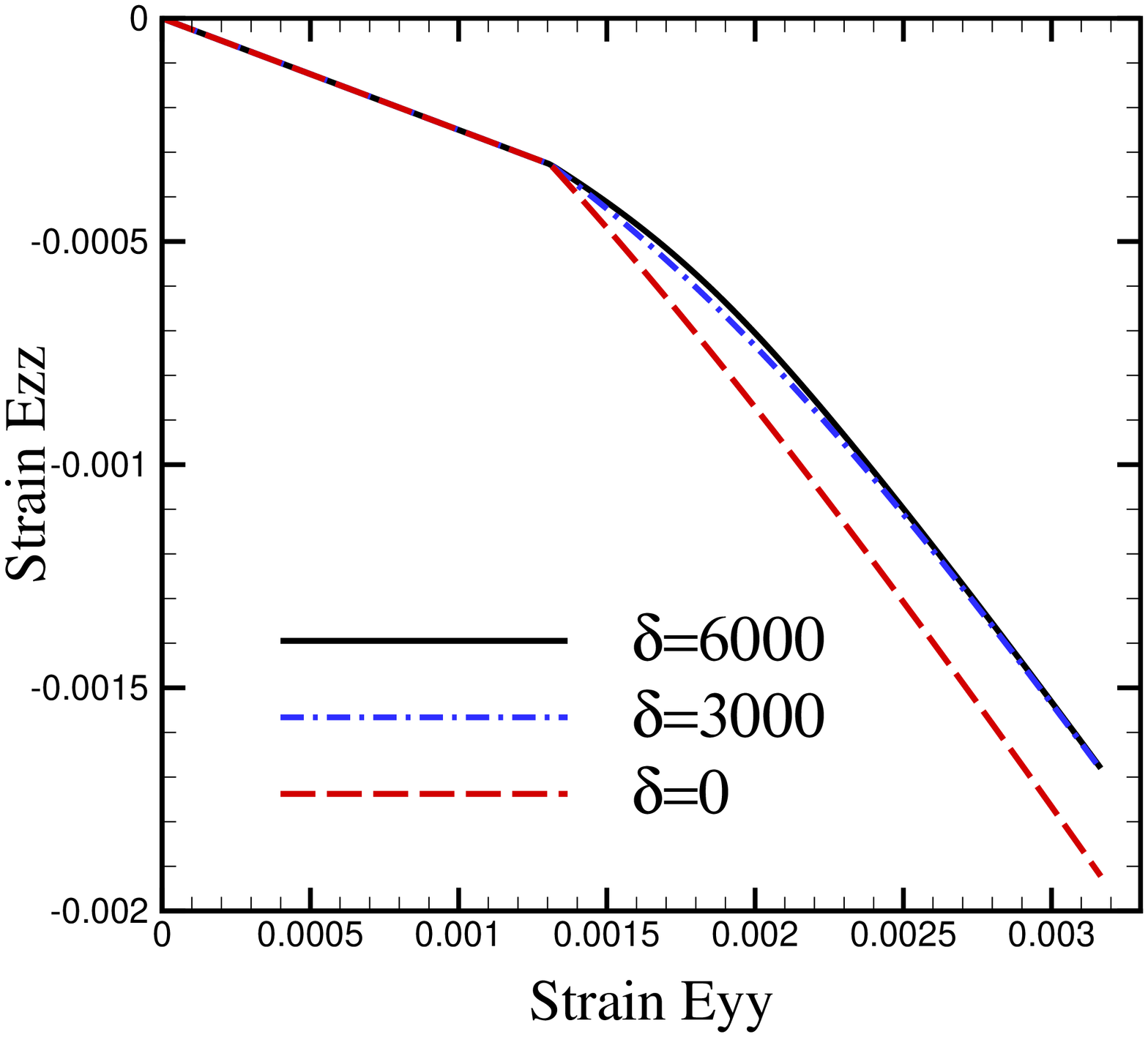}
  \end{tabular} 
 \end{minipage}
\caption{Radial contraction of an elasto-plastic annulus, design of a smooth transition from elastic-to-plastic (left) in 
         the stress-strain curve by a proper choice of the hardening variable $\delta$. A smooth stress-strain curve directly 
         propagates to a corresponding (right) smooth transition in the strain-path.\label{fig:smooth-glatt}}
\end{Figure}

%\noindent
In the following, 4 distinct cases are analyzed, which differ in their plastic material parameters:
\begin{itemize}
\item {\bf A0\,:}            \quad $\sigma_Y=0$, linear hardening          \\[-6mm]
\item {\bf B0\,:}            \quad $\sigma_Y=0$, exponential hardening     \\[-6mm]
\item {\bf A\phantom{0}\,:}  \quad $\sigma_Y>0$, linear hardening          \\[-6mm]
\item {\bf B\phantom{0}\,:}  \quad $\sigma_Y>0$, exponential hardening for a smooth elastic-plastic transition 
\end{itemize}

In the left of Fig.~\ref{fig:smooth-glatt} the stress-strain curve of a (homogeneous) bi-axial tension-compression test 
is shown for various hardening parameters $\delta$ in the exponential/saturation part of isotropic hardening 
$\hat{K'}(\alpha) = H\alpha + (\sigma_\infty - \sigma_Y)(1-e^{-\delta\alpha})$, where $H=600$\, N/mm$^2$. 
For case B, we choose $\delta =5\,000$ for the calculations, since it ensures not only a truly seamless transition from 
elastic-to-plastic state in the stress-strain curve but also a corresponding smooth transition in the strain components.

%\vfill
%\newpage

\subsubsection{Case A0: $\Bsigma_Y=\mathbf{0}$, linear hardening.}
%-----------------------------------------------------------------------------------------------------------------------------
\begin{Table}[htbp]
\center
\renewcommand{\arraystretch}{1.4}
%\begin{tabular}{llcc}
\begin{tabular}{ccccccc}
\hline
 case                                                        &
 $E$ $({\text{N}}/{\text{mm}^2})$                            &
 $\nu$ $(1)$                                                 &
 $\sigma_Y$ $({\text{N}}/{\text{mm}^2})$                     &
 $\sigma_\infty - \sigma_Y$ $({\text{N}}/{\text{mm}^2})$     &
 $H$ $({\text{N}}/{\text{mm}^2})$                            &
 $\delta$  $(1)$                                            \\
\hline
A0 & $68\,900$ &  $0.33$  &  $0.0$   &  $0.0$    &  $10\,000$  &  $0.0$    \\ 
% \hline
B0 &           &          &  $0.0$   &  $200$    &  $3\,000$   &  $5\,000$ \\ 
% \hline
 A  &          &          &  $300$   &   $0.0$   &  $10\,000$  &  $0.0$    \\
% \hline
 B  &          &          &  $300$   &   $200$   &  $3000$     &  $5\,000$ \\
 \hline 
\end{tabular}
\caption{Radial contraction of an annulus, elasto-plastic material parameters for different cases.}
\label{tab:KonvergenzKreisringElastoPlasto}
\end{Table}

The reference solution for the error analysis is calculated employing Radau IIa along with quadratic interpolation for
a time step size of $\Delta t=1.0$E$-05$.

\begin{Table}[htbp]
\center
{%\small
\renewcommand{\arraystretch}{1.4}
\begin{tabular}{lrcclcrcclcrccl}
\hline
error     &  \multicolumn{4}{ c }{$e(\bE)$}  & \,\,  & \multicolumn{4}{ c }{$e(\bE^{\mbox{\scriptsize p}})$} & \,\,  &  \multicolumn{4}{ c }{$e(\bS)$}   \\
 time $t$ &    0.10  &  0.25   &  0.50   &   &      &  0.10 &  0.25  &  0.50   &      &       &  0.10  &  0.25   &  0.50   &      \\
 \hline
 BE        &    1.00  &  1.00   &  1.00  &   &       &  1.00 &  1.00  &  1.00   &      &       &  1.00 &   1.00  &  1.00  &        \\
 RIIa-l      &    1.61  &  1.65   &  1.69  &   &       &  1.61 &  1.65  &  1.69   &      &       &  1.63 &   1.67  &  1.70  &        \\[2mm]
 \cellcolor{hellgrau} RIIa-q      &    2.05  &  2.06   &  2.05  &   &       &  2.05 &  2.06  &  2.04   &      &       &  2.07 &   2.05  &  2.02  &    \\ 
% \cellcolor{hellgrau} RIIa-q-i    &    0.97  &  1.20   &  1.40  &   &       &  0.98 &  1.21  &  1.41   &      &       &  1.04 &   1.27  &  1.46  &        \\
 \hline
\end{tabular}
%\newline
\caption{Radial contraction of an elasto-plastic annulus, {\bf case A0}: order of convergence for different methods.}
\label{tab:ConvergenceAnnulus-ElPl-y0=0-LinHard}
}
\end{Table}

%\bigskip
  
\begin{Figure}[htbp]
  \begin{minipage}{15.5cm}
     \centering%
     \tabcolsep1mm%
  \begin{tabular}[]{rcc}
    &  $t=0.1$  &  $t=0.5$  \\ 
\raisebox{3.2cm}{$e(\bS)$}
  &
   \includegraphics[width=.40\textwidth, angle=0, clip=]{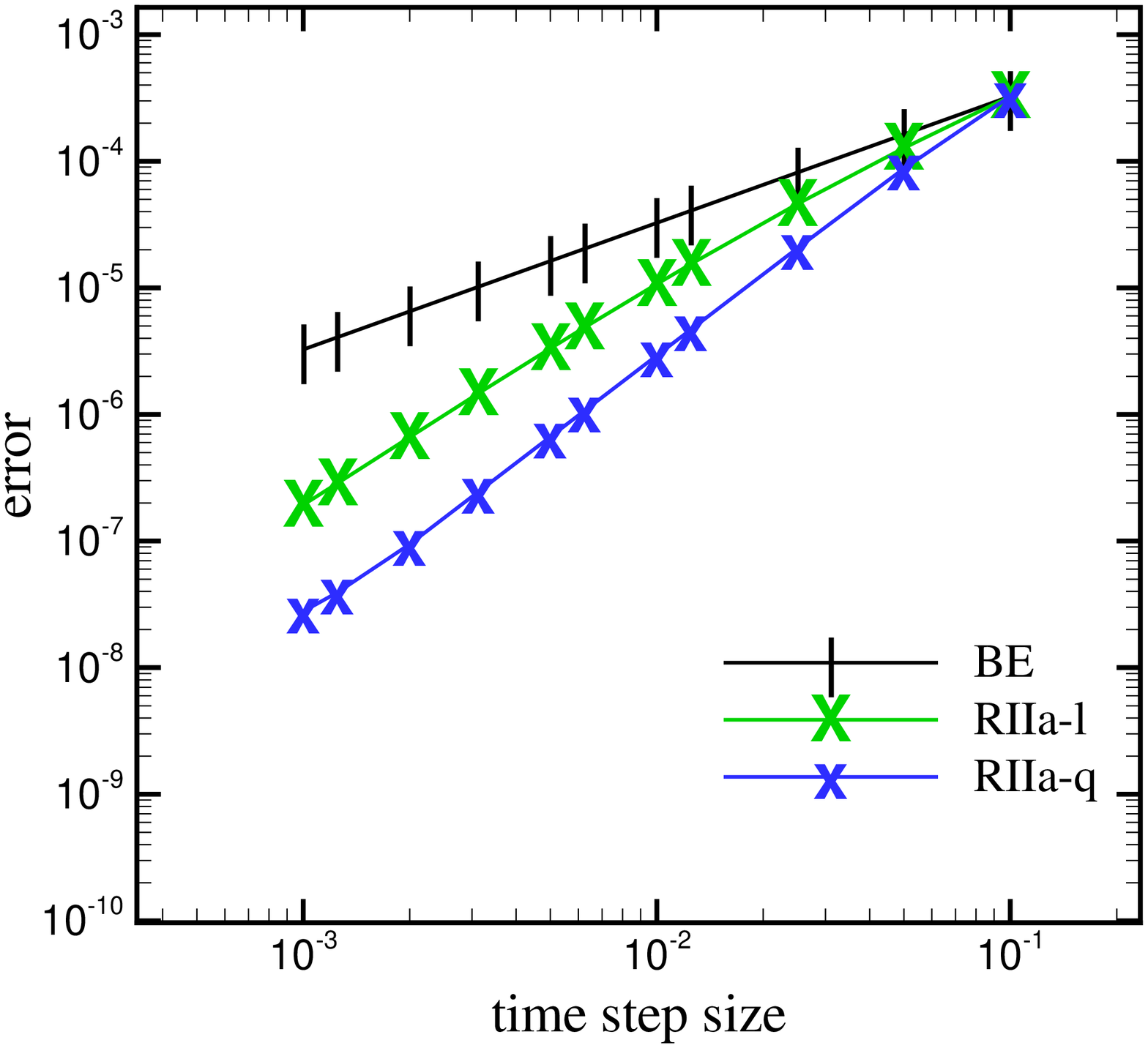} \hspace{5mm}
  &
   \includegraphics[width=.40\textwidth, angle=0, clip=]{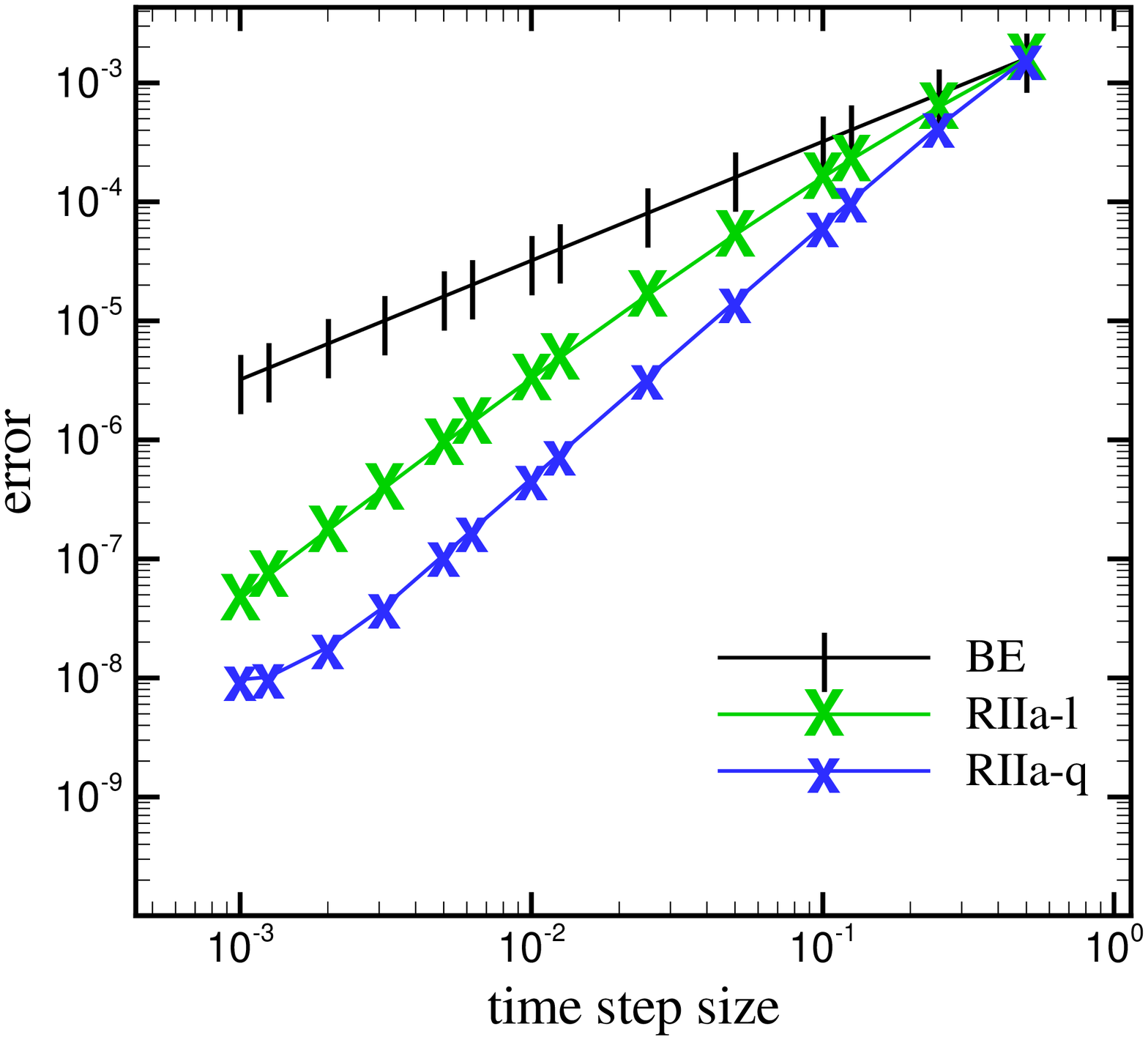}
 \end{tabular}
  \end{minipage}
\caption{Radial contraction of an annulus, {\bf case A0}: convergence of $e(\bS)$ 
at $t=0.1\,$ (left) and at $t=0.5\,$ (right).\label{fig:RadialContraction-Convergence-y0=0-LinHard}}
\end{Figure}

The simulation results are shown in Tab.~\ref{tab:ConvergenceAnnulus-ElPl-y0=0-LinHard} and in 
Fig.~\ref{fig:RadialContraction-Convergence-y0=0-LinHard}. 
Radau IIa along with linear interpolation (RIIa-l) exhibits a reduced order in the range of $1.6-1.7$.
Radau IIa with quadratic interpolation (RIIa-q) is slightly better but hardly above 2nd order.
The reason is that plasticity starts at $t=0$ and therefore, interpolation in the first plastic
time step is necessarily linear, since $t_{n-1}$ data are missing. 
%If, however, an initial step
%is chosen, a quadratic polynomial can be constructed and an improvement in convergence can be achieved.
%This can be seen in Fig.~\ref{fig:RadialContraction-Convergence-y0=0-LinHard} as well as in Table \ref{tab:ConvergenceAnnulus-ElPl-y0=0-LinHard}.

A comparison of the absolute errors for different methods shows, that third order Radau IIa is already more accurate for the maximum
time step of $\Delta t=0.125$ than Backward-Euler for the minimum time step size $\Delta t=0.010$, which even applies for Radau IIa
suffering from order reduction.

%\vfill
%\newpage

\subsubsection{Case B0: $\Bsigma_Y=\mathbf{0}$, exponential hardening.}
%-----------------------------------------------------------------------------------------------------------------------------

The reference solution for the error analysis is calculated employing Radau IIa along with quadratic interpolation for
a time step size of $\Delta t=1.0$E$-05$.

\begin{Table}[htbp]
\center
{%\small
\renewcommand{\arraystretch}{1.4}
\begin{tabular}{lrcclcrcclcrccl}
\hline
error     &  \multicolumn{4}{ c }{$e(\bE)$}  & \,\,  & \multicolumn{4}{ c }{$e(\bE^{\mbox{\scriptsize p}})$} & \,\,  &  \multicolumn{4}{ c }{$e(\bS)$}   \\
 time $t$ &    0.10  &  0.25   &  0.50   &   &      &  0.10 &  0.25  &  0.50   &      &       &  0.10  &  0.25   &  0.50   &      \\
 \hline
 BE        &    0.91  &  0.86   &  0.77  &   &       &  0.91 &  1.00  &  0.77   &      &       &  0.91 &   0.85  &  0.93  &        \\
 RIIa-l      &    2.02  &  2.06   &  2.04  &   &       &  2.00 &  1.65  &  2.04   &      &       &  2.04 &   2.09  &  2.03  &        \\[2mm]
\rowcolor{hellgrau}
 RIIa-q      &    3.10  &  2.97   &  2.69  &   &       &  2.98 &  2.06  &  2.69   &      &       &  3.08 &   2.96  &  2.71  &        \\[2mm]
% RIIa-q-i    &    2.95  &  2.58   &  2.62  &   &       &  2.92 &  1.21  &  2.61   &      &       &  2.92 &   2.78  &  2.66  &        \\
 \hline
\end{tabular}
%\newline
\caption{Radial contraction of an annulus, {\bf case B0}: order of convergence for different methods.}
\label{tab:ConvergenceAnnulus-ElPl-y0=0-expHard}
}
\end{Table}
 
\begin{Figure}[htbp]
  \begin{minipage}{15.5cm}
     \centering%
     \tabcolsep1mm%
  \begin{tabular}[]{rcc}
    &  $t=0.1$  &  $t=0.5$  \\
\raisebox{3.0cm}{$e(\bS)$}
  &
   \includegraphics[width=.36\textwidth, angle=0, clip=]{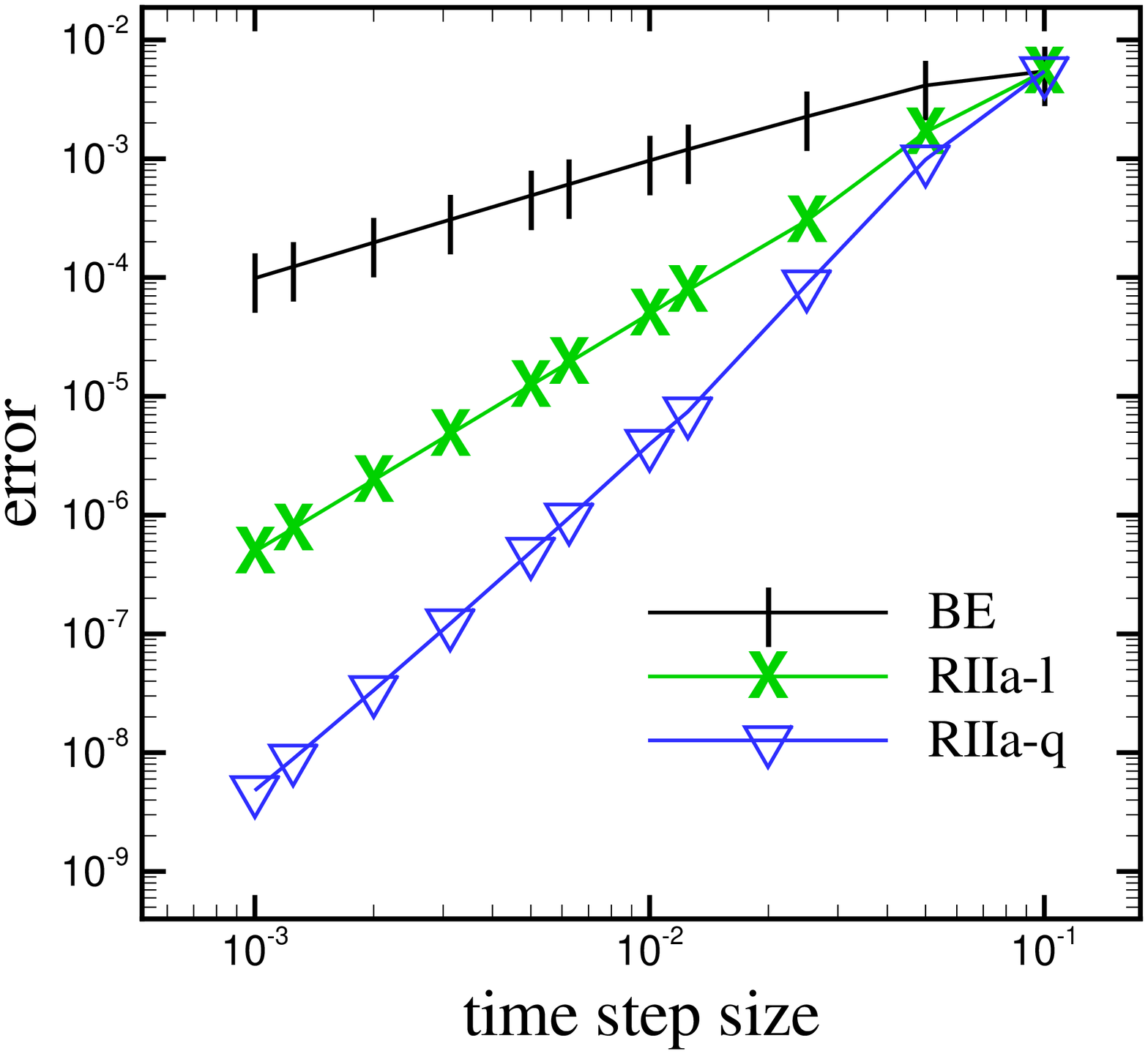} \hspace{5mm}
  &
   \includegraphics[width=.36\textwidth, angle=0, clip=]{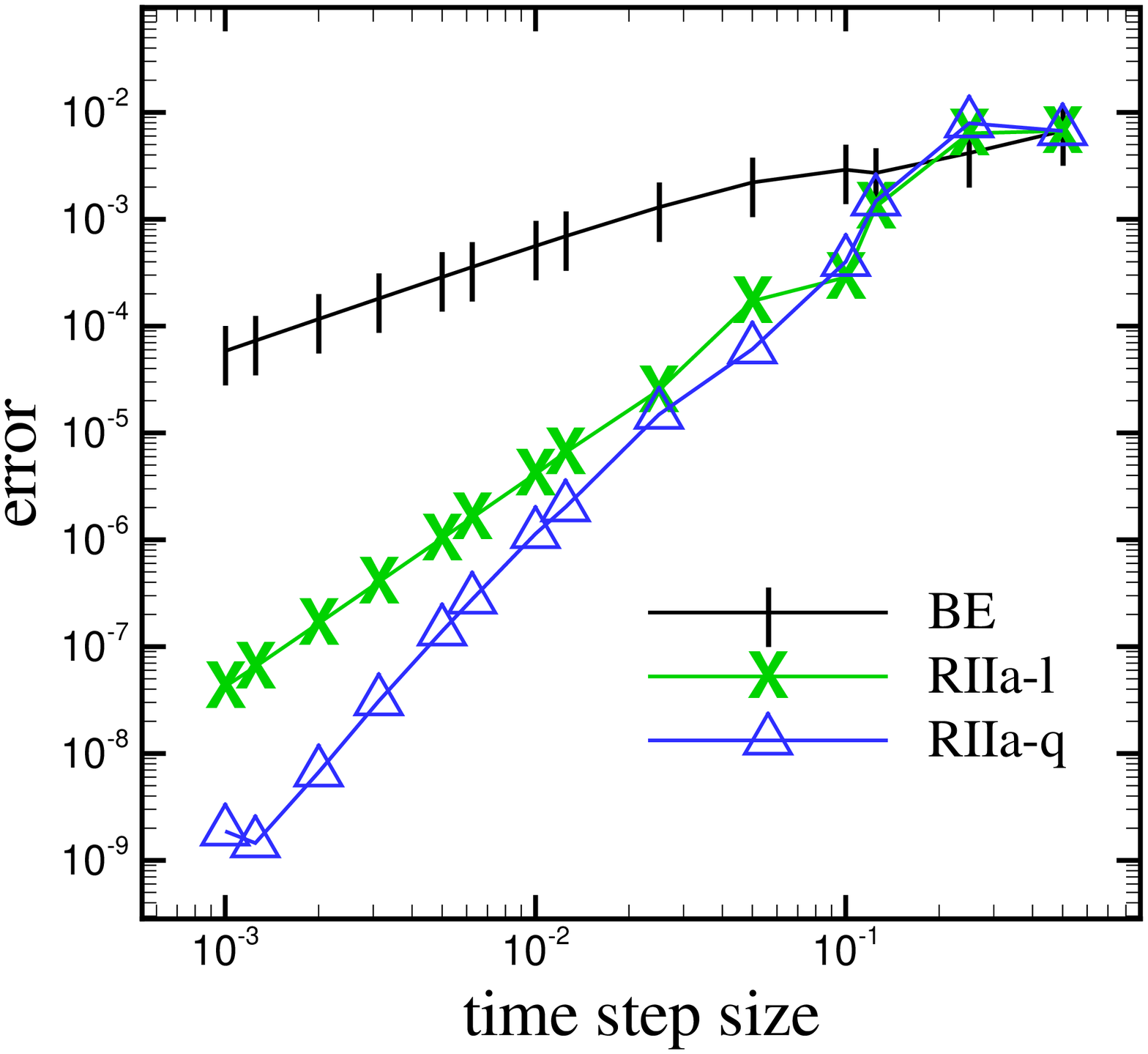}
 \end{tabular}
  \end{minipage}
\caption{Radial contraction of an elasto-plastic annulus, {\bf case B0}: convergence of the error $e(\bS)$ at $t=0.1\,$ (left)
and at $t=0.5\,$ (right).\label{fig:RadialContraction-Convergence-y0=0-expHard} }
\end{Figure}
 
The simulation results are shown in Tab.~\ref{tab:ConvergenceAnnulus-ElPl-y0=0-expHard} and in 
Fig.~\ref{fig:RadialContraction-Convergence-y0=0-expHard}. 
Radau IIa along with linear interpolation (RIIa-l) exhibits order reduction resulting in quadratic convergence.
Radau IIa with quadratic interpolation (RIIa-q) shows full convergence order 3. 
  
A comparison of the errors for different methods shows, that third order Radau IIa is already more accurate for the maximum
time step of $\Delta t=0.125$ than Backward-Euler for the minimum time step size $\Delta t=0.010$, which even applies for Radau IIa
suffering from order reduction.

%\vfill
%\newpage

\subsubsection{Case A: $\Bsigma_Y>\mathbf{0}$, linear hardening.}
%-----------------------------------------------------------------------------------------------------------------------------

\begin{Table}[htbp]
\center
{%\small %\scriptsize
\renewcommand{\arraystretch}{1.2}
\begin{tabular}{lrclcrclcrcl}
\hline
error     &  \multicolumn{3}{ c }{$e(\bE)$}  & \,\,  & \multicolumn{3}{ c }{$e(\bE^{\mbox{\scriptsize p}})$} & \quad &  \multicolumn{3}{ c }{$e(\bS)$}    \\
 time $t$ &    0.10  &  0.25  &            0.50 &       &  0.10 & 0.25 &             0.50    &       &  0.10  &  0.25  &     0.50     \\
 \hline
 BE        &    0.92  &  0.94   &           0.85 &       &  0.89 &  0.94  &            0.85    &       &  0.96  &   0.94  &     0.93     \\[2mm]
 RIIa-l      &    2.28  &  1.94   &           1.99 &       &  2.15 &  1.94  &            2.00    &       &  2.34  &   2.15  &     2.09     \\
 RIIa-l-SP   &    2.17  &  1.97   &           2.10 &       &  2.12 &  1.95  &            2.09    &       &  2.18  &   2.12  &     2.12     \\
 RIIa-q      &    2.73  &  2.29   &           1.99 &       &  2.36 &  2.20  &            2.00    &       &  2.65  &   2.39  &     2.13     \\
 \cellcolor{hellgrau} RIIa-q-SP   &    2.17  &  2.36   &           2.34 &       &  2.24 &  2.34  &            2.32    &       &  2.33  &   2.19  &     2.18     \\
 \cellcolor{hellgrau} RIIa-q-exSP &    2.28  &  2.22   &           2.29 &       &  2.58 &  2.33  &            2.29    &       &  2.54  &   2.15  &     2.13     \\
 \hline
\end{tabular}
}
%\newline
\caption{Radial contraction of an annulus, {\bf case A}: order of convergence for different methods.}
\label{tab:StripWithHole-Convergence-A}
\end{Table}
  
\begin{Figure}[htbp]
  \begin{minipage}{15.5cm}
     \centering%
     \tabcolsep1mm%
  \begin{tabular}[]{rcc}
    &  $t=0.1$  &  $t=0.5$  \\ 
       \raisebox{3.0cm}{$e(\bS)$}
    &
       \includegraphics[width=.38\textwidth, angle=0, clip=]{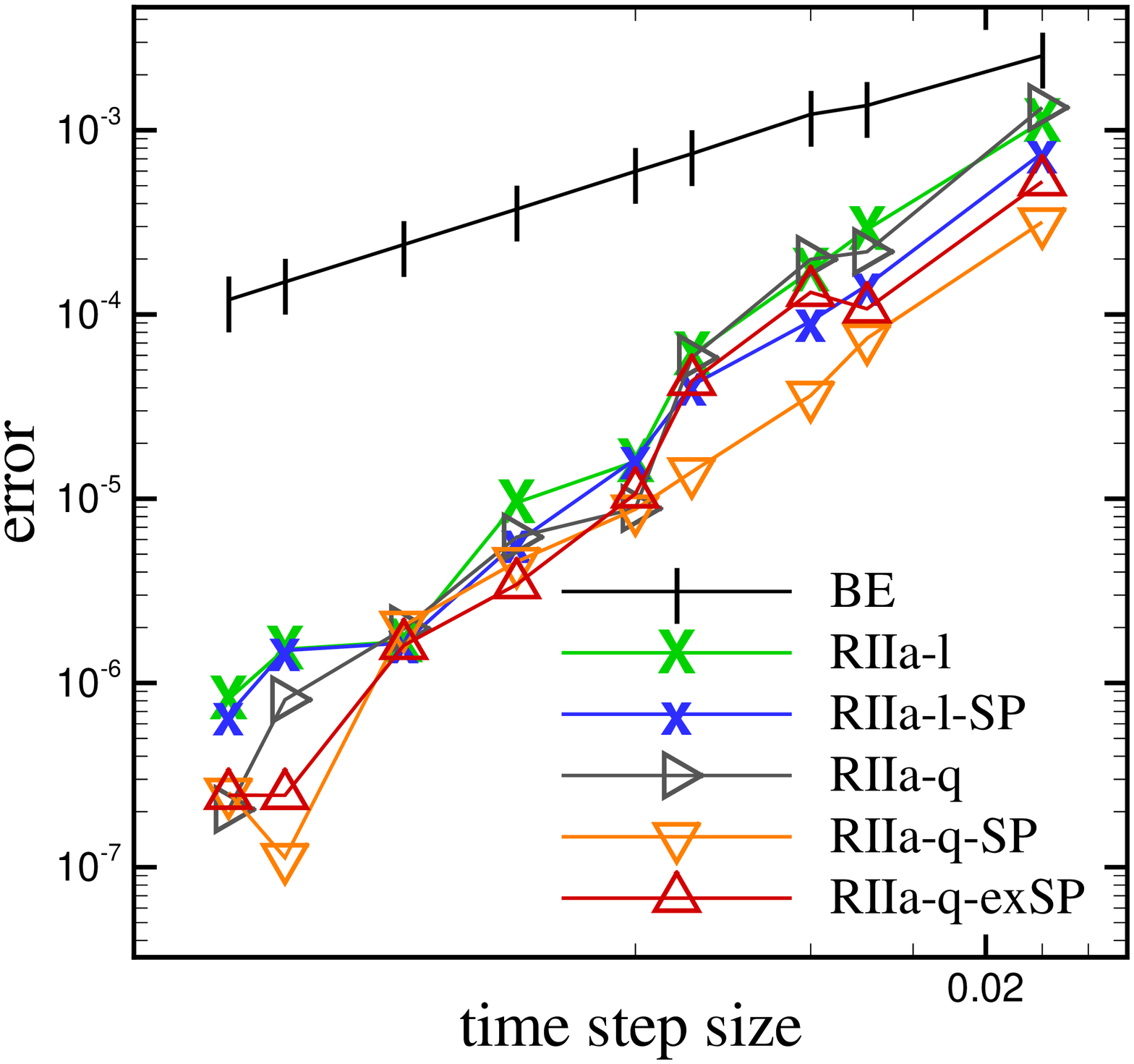} \hspace{5mm}
    &
       \includegraphics[width=.38\textwidth, angle=0, clip=]{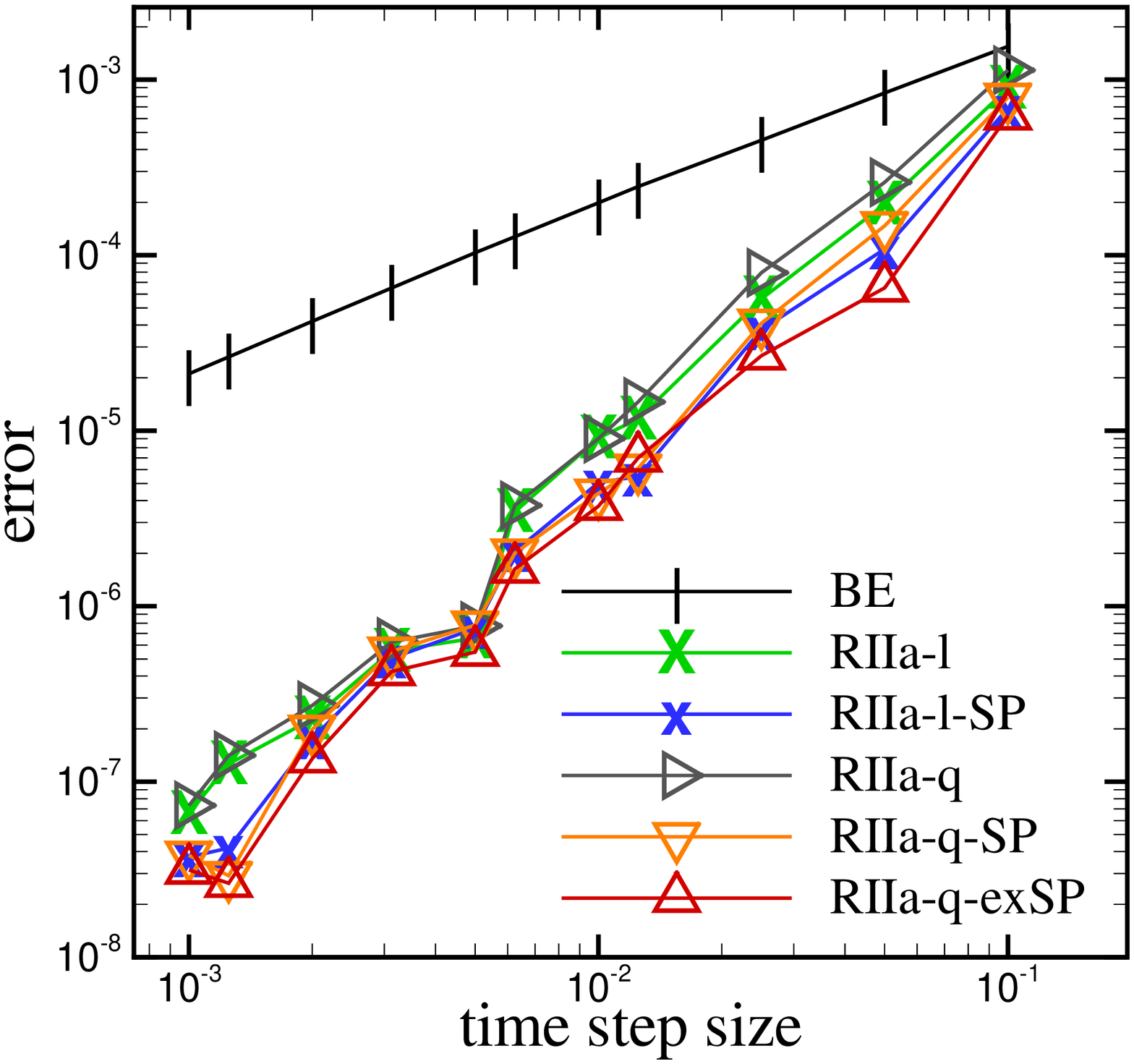}
  \end{tabular}
  \end{minipage}
\caption{Radial contraction of an elasto-plastic annulus, {\bf case A}: convergence of the error $e(\bS)$ at $t=0.1\,$ (left)
and at $t=0.5\,$ (right).\label{fig:Annulus-CaseA-Errors}}
\end{Figure}

Opposed to the uniform convergence of Backward-Euler, all variants of Radau IIa do not exhibit strictly uniform convergence behavior, see
Fig.\ref{fig:Annulus-CaseA-Errors}.

Radau IIa with linear interpolation, with or without SP detection shows order reduction to second order.
Radau IIa with quadratic interpolation slightly improves convergence, but toggles in between order 2 and 2.7. SP 
detection does not lead to a clear improvement. This is true for all considered quantities, for $\bE$, $\bE^{\mbox{\scriptsize p}}$ and for $\bS$,
all quantities show a similar error pattern, see Tab.\ref{tab:StripWithHole-Convergence-A}.

%\vfill
%\newpage

\subsubsection{Case B: $\Bsigma_Y>\mathbf{0}$, exponential hardening for a smooth elasto-plastic transition.}
  
\begin{Table}[htbp]
\center
{%\small %\scriptsize
\renewcommand{\arraystretch}{1.4}
\begin{tabular}{lcrclcrclcrcl}
\hline
error      & \hspace*{2mm}  & \multicolumn{3}{ c }{$e(\bE)$}  & \hspace*{3mm}  & \multicolumn{3}{ c }{$e(\bE^{\mbox{\scriptsize p}})$} & \hspace*{3mm}  &  \multicolumn{3}{ c }{$e(\bS)$}   \\
time $t$   &                & 0.10  &  0.25     &  0.50 &       &  0.10 &  0.25   &    0.50    &       &  0.10 &   0.25   &    0.50    \\
 \hline
BE          &                & 0.91  &  0.95     &  0.95 &       &  0.86 &  0.95  &    0.95    &       &  0.95 &   0.96  &    0.94     \\[2mm]
RIIa-l      &                & 2.05  &  2.03     &  2.05 &       &  2.02 &  2.02  &    2.06    &       &  2.09 &   2.02  &    2.17     \\
RIIa-l-SP   &                & 1.99  &  2.03     &  2.06 &       &  1.96 &  2.02  &    2.06    &       &  2.00 &   2.02  &    2.17     \\[2mm]
\cellcolor{hellgrau}RIIa-q   &  & 2.26  &   \cellcolor{hellgrau} 2.91 &  \cellcolor{hellgrau}2.76 &       &  2.22 &  \cellcolor{hellgrau}2.61  &  \cellcolor{hellgrau}2.77  &  &  2.24 & \cellcolor{hellgrau}2.94  &  \cellcolor{hellgrau}2.95  \\
\cellcolor{hellgrau}RIIa-q-SP   &  & 2.12  &   \cellcolor{hellgrau} 2.82 &   \cellcolor{hellgrau}2.70 &       &  2.13 &   \cellcolor{hellgrau}2.74  & \cellcolor{hellgrau}2.70 &  &  2.14 & \cellcolor{hellgrau}2.92  &  \cellcolor{hellgrau}2.95   \\
\cellcolor{hellgrau}RIIa-q-exSP &  & 2.25  &   \cellcolor{hellgrau} 2.84 &   \cellcolor{hellgrau}2.63 &       &  2.40 &   \cellcolor{hellgrau}2.75  & \cellcolor{hellgrau}2.64 &  &  2.47 & \cellcolor{hellgrau}2.93  &  \cellcolor{hellgrau}2.94   \\
 \hline
\end{tabular}
%\newline
}
\caption{Radial contraction of an elasto-plastic annulus, {\bf case B}: order of convergence for different methods.}
\label{tab:StripWithHole-Convergence-B}
\end{Table}

\begin{Figure}[htbp]
  \begin{minipage}{15.5cm}
     \centering%
     \tabcolsep1mm%
  \begin{tabular}[]{rcc}
    &  $t=0.1$  &  $t=0.25$  \\ 
       \raisebox{3.2cm}{$e(\bS)$}
    &
       \includegraphics[width=.42\textwidth, angle=0, clip=]{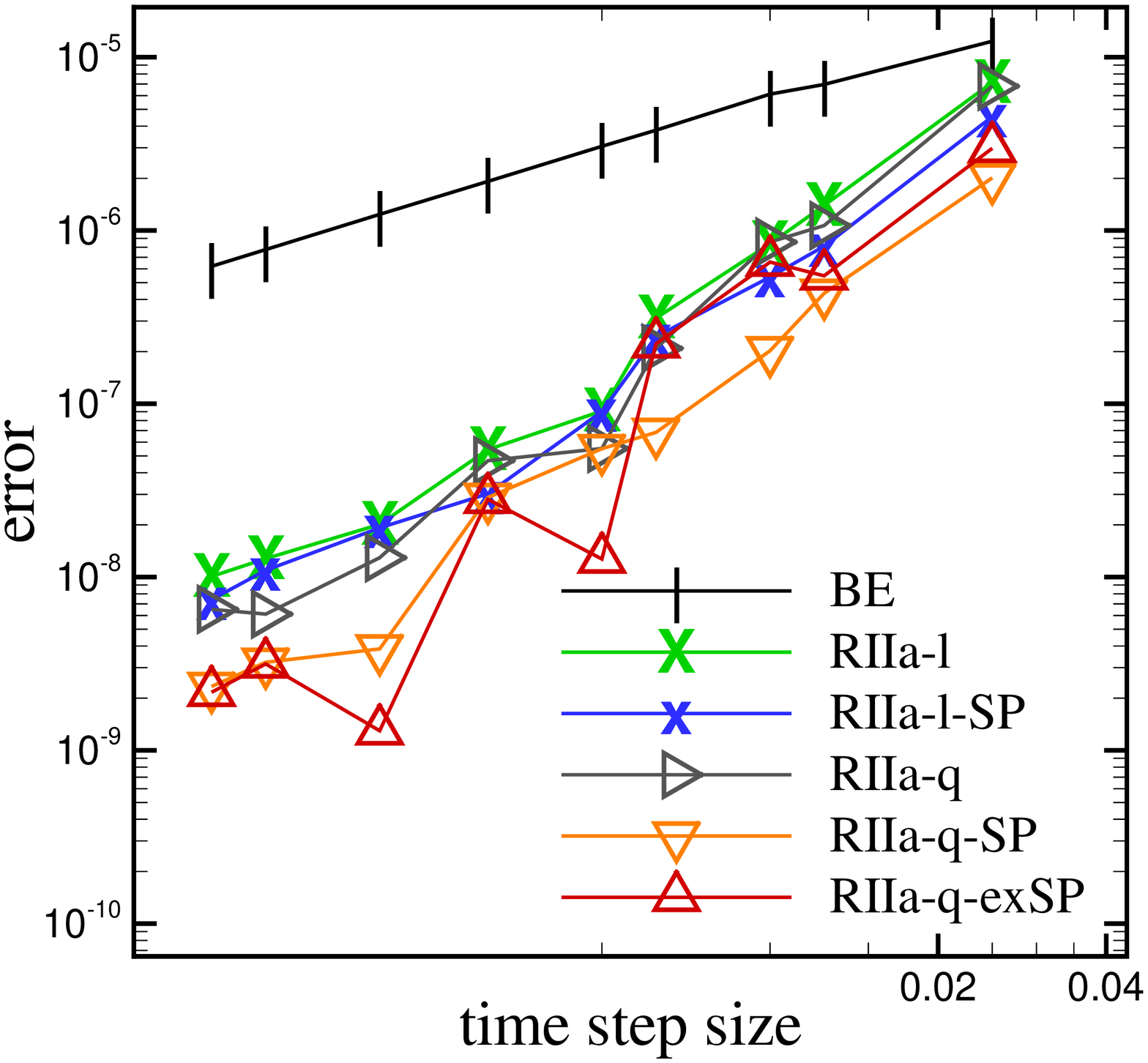} \hspace{5mm}
    &
       \includegraphics[width=.42\textwidth, angle=0, clip=]{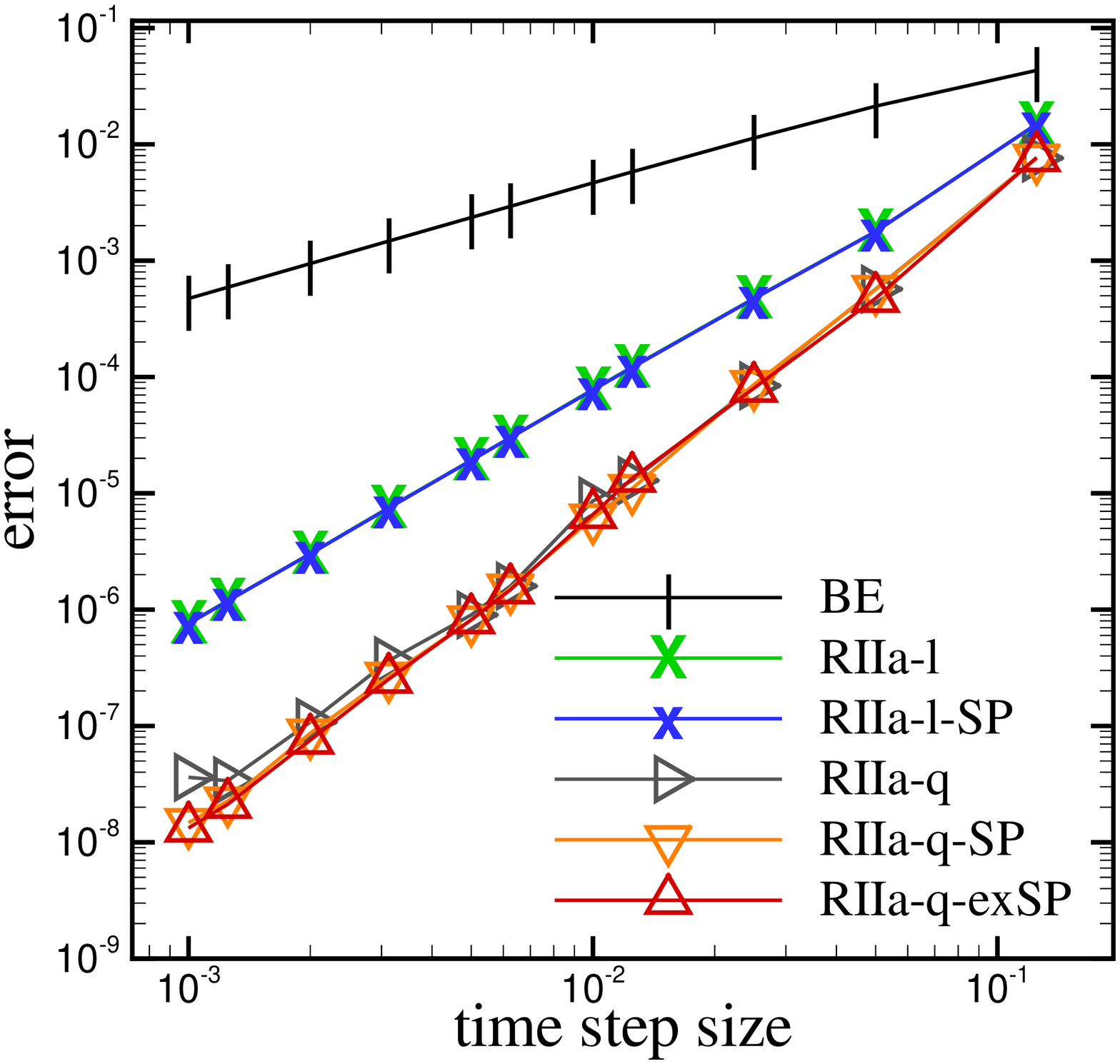}
  \end{tabular}
  \end{minipage}
\caption{Radial contraction of an elasto-plastic annulus, {\bf case B}: convergence of error $e(\bS)$ for different methods at (left) time $t=0.1$ and 
(right) $t=0.25$.\label{fig:Annulus-CaseB-Errors}}
\end{Figure}

For the present smooth case, Radau IIa along with quadratic strain interpolation shows full convergence order 3, see 
Tab.~\ref{tab:StripWithHole-Convergence-B} and Fig.~\ref{fig:Annulus-CaseB-Errors}. Remarkably, the additional switching 
point detection does not improve the convergence. Full convergence order is observed for {\color{black}$\bE$}, which indicates that 
the total strain path in time is equally smooth as the stress-strain curve. Note, that this smoothness was realized by 
the choice of the hardening parameters. Then, the consistent quadratic interpolation is effective and enables the full convergence order in time integration as 
reflected in the order of $\bE^{\mbox{\scriptsize p}}$ as well as $\bS$. 

\begin{Figure}[htbp]
   \begin{minipage}{15.0cm}
     \begin{center}
           \includegraphics[scale=0.35, clip=]{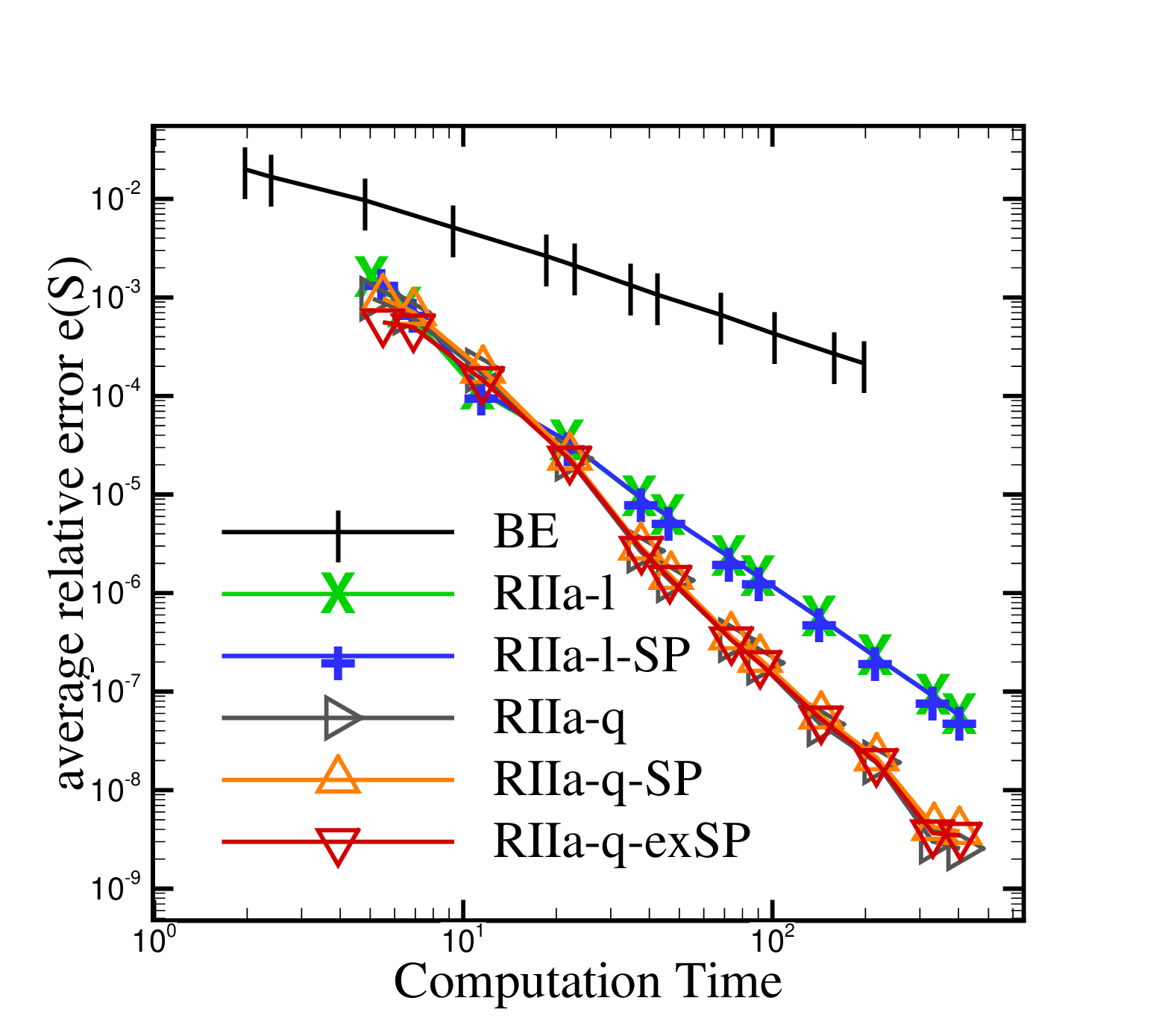}
     \end{center}
   \end{minipage}
\caption{Radial contraction of an elasto-plastic annulus, {\bf case B}: efficiency in terms of the error $e(\bS^{\mbox{\scriptsize}})$ 
versus the overall computation time (s) at time $t=0.5$.\label{fig:Efficiency-Annulus-CaseB}}

%\begin{table}[htbp]
\end{Figure}

\begin{Table}[htbp]
\center
\renewcommand{\arraystretch}{1.1}
\begin{tabular}{lcrrrr}
\hline
error tol. $e(\bS^{\mbox{\scriptsize}})$ & &    1.0E-03     &    &    1.0E-04        &                   \\
                &          &       speed-up    &    &       speed-up &                   \\
 \hline
 BE             &          &           1.0     &    &           1.0  &                   \\[2mm]
 RIIa-l         &          &           7.2     &    &          34.2  &                   \\[2mm]
 RIIa-l-SP      &          &           7.7     &    &          33.3  &                   \\[2mm]
 RIIa-q         &          &           8.1     &    &          29.9  &                   \\
 RIIa-q-SP      &          &           8.7     &    &          29.3  &                   \\[2mm]
 RIIa-q-exSP    &          &          12.2     &    &          31.1  &                   \\
 \hline
\end{tabular}
%\newline
\caption{Radial contraction of an elasto-plastic annulus, {\bf case B}: Speed-up factors of different methods compared with Backward-Euler for 
different error tolerances. Error calculations at time $t=0.25$.\label{tab:SpeedupAnnulus-CaseB}}
\end{Table}

Full convergence order 3 and the corresponding high accuracy of Radau IIa is a nice result, but the method is due to
its 2-stage, fully implicit characteristics more expensive than Backward-Euler. The true improvement of the 3rd order
method must be quantified in terms of the realized speed-up compared with Backward-Euler. Figure~\ref{fig:Efficiency-Annulus-CaseB}
and Tab.~\ref{tab:SpeedupAnnulus-CaseB} show, that for an error tolerance of 1.0E-03 Radau IIa already exhibits 
a speed-up of 7-12$\times$. If a higher accuracy shall be achieved in terms of a smaller error tolerance of 
1.0E-04, the speed-up is considerably increased to factors of approx. 30$\times$. 

Figure~\ref{fig:Annulus-CaseC-Stress-EqPlStr} shows for Case B the annulus distribution of shear stresses and equivalent plastic
strain on the increasingly deformed structure. 
   
\begin{Figure}[htbp]
   \begin{minipage}{15.5cm}
     \centering
          \begin{minipage}{77.5mm}
           \includegraphics[scale=0.35, clip=]{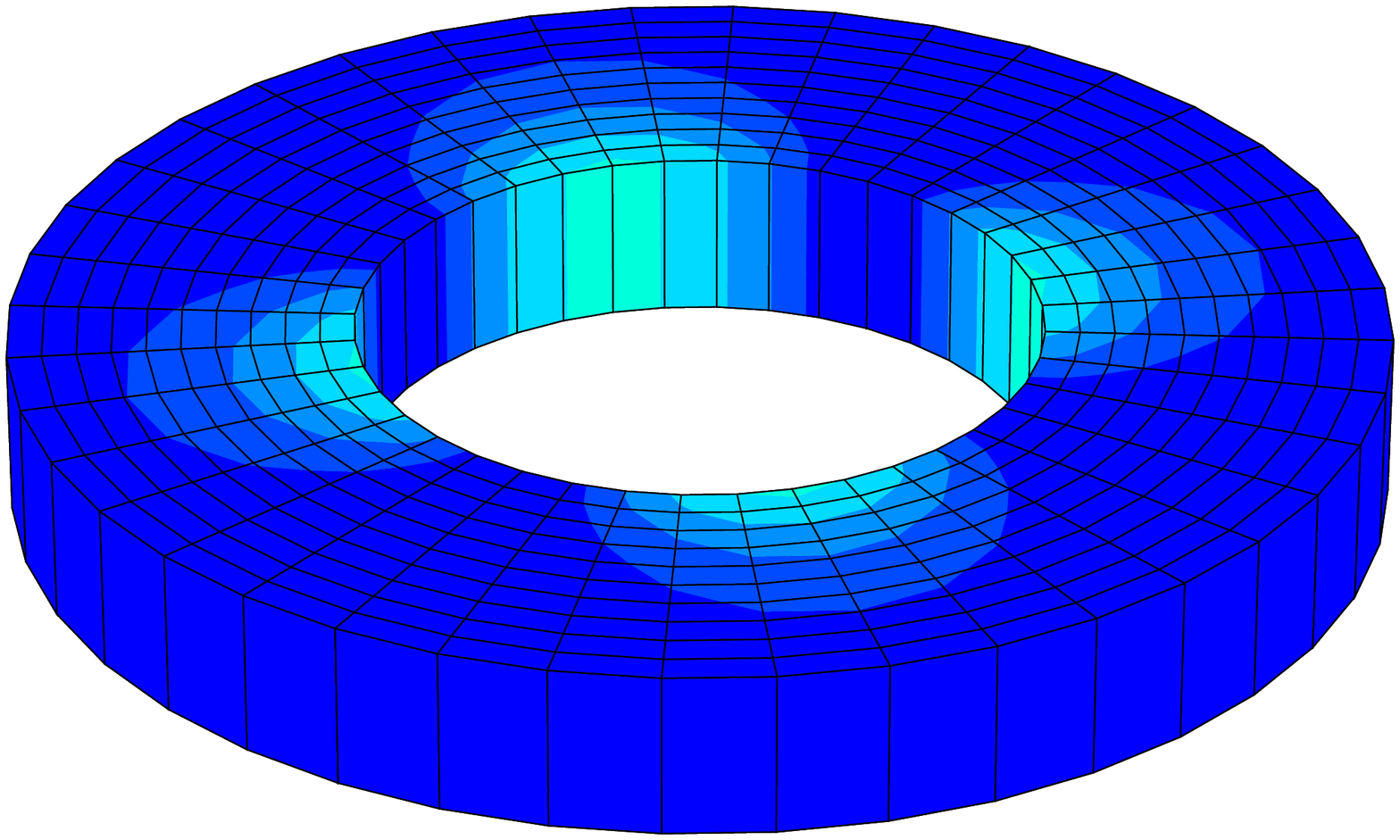}
          \end{minipage}\hfil
          \begin{minipage}{77.5mm}
            \includegraphics[scale=0.35, clip=]{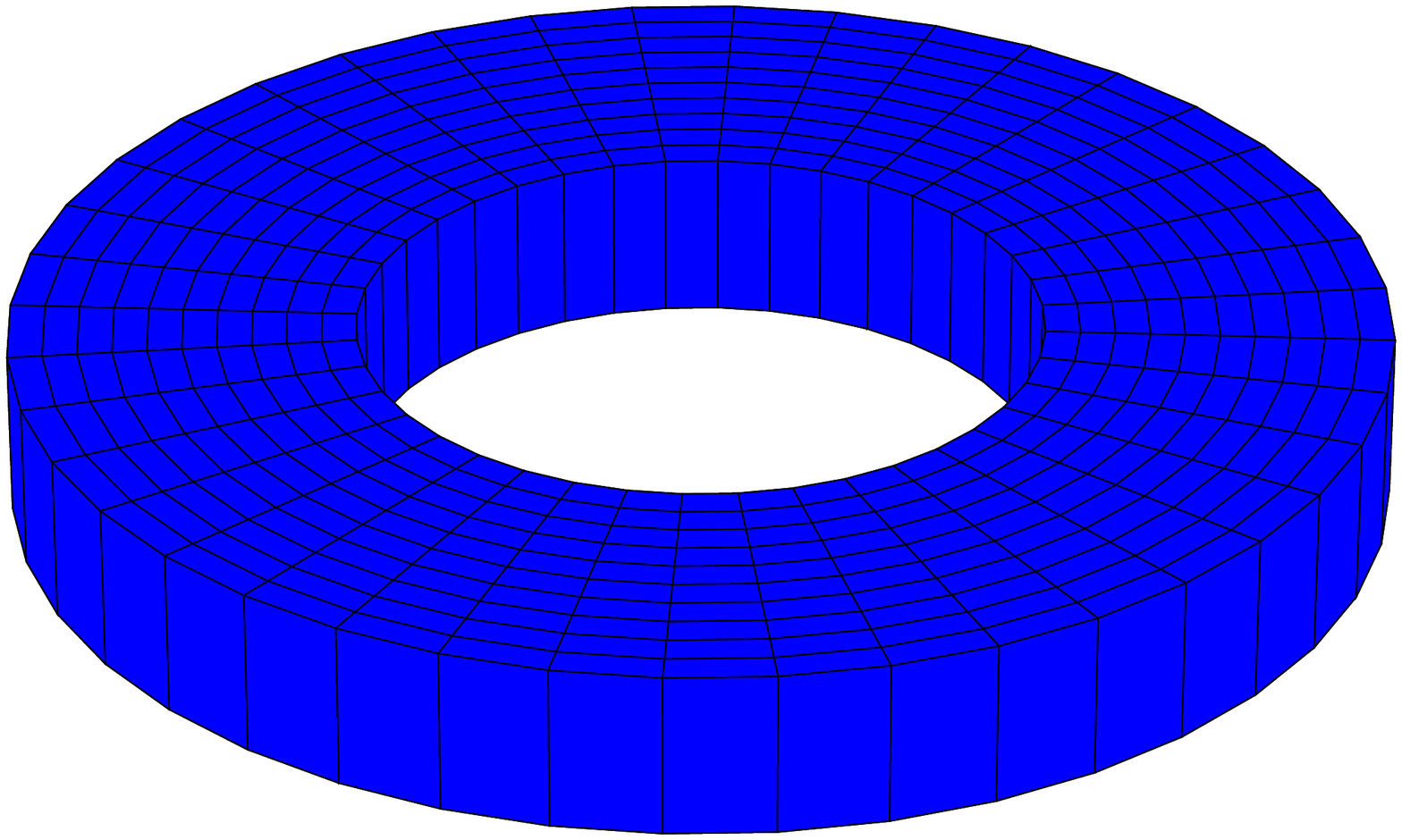}
          \end{minipage}
          \\[0mm]
          \begin{minipage}{77.5mm}
           \includegraphics[scale=0.35, clip=]{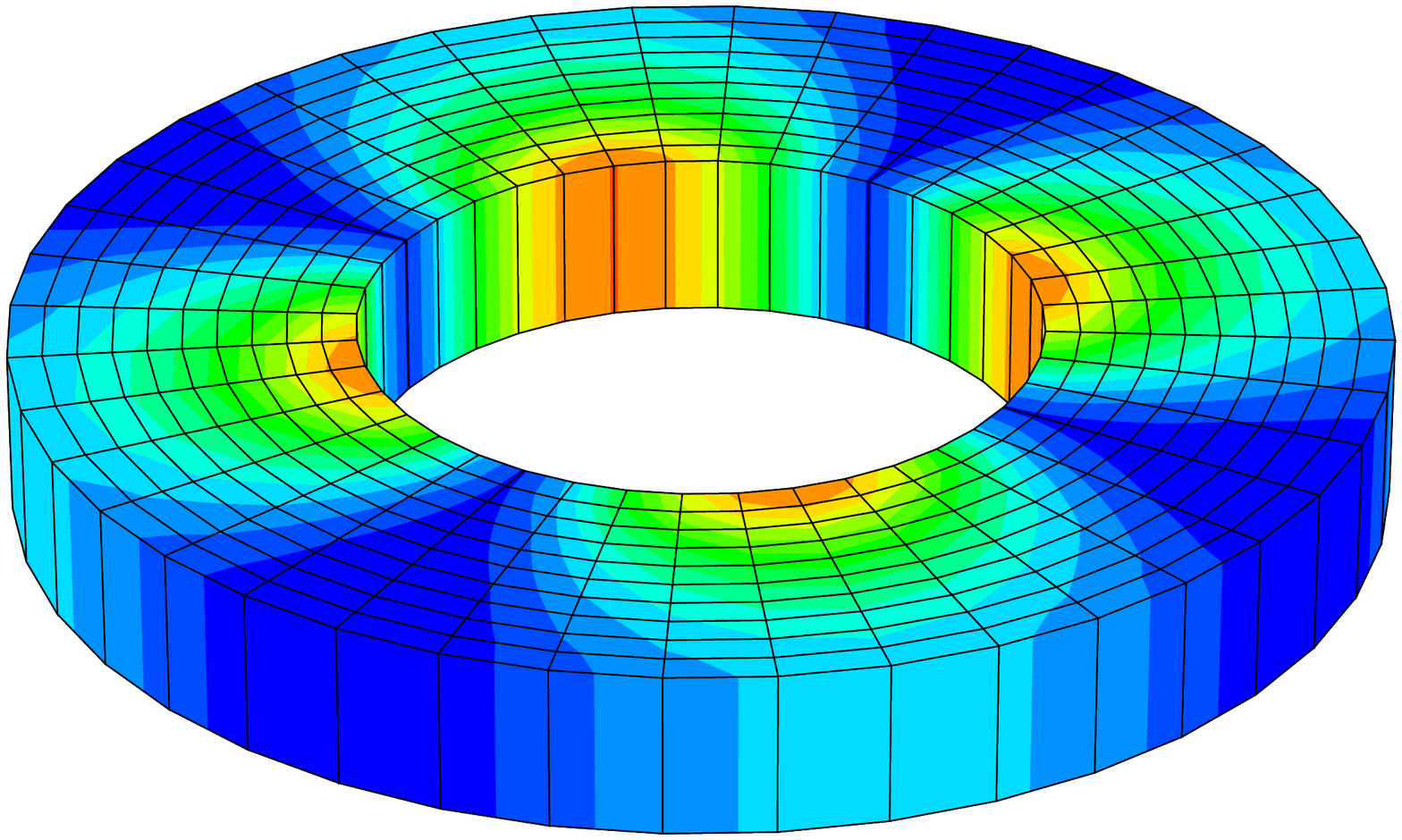}
          \end{minipage}\hfil
          \begin{minipage}{77.5mm}
            \includegraphics[scale=0.35, clip=]{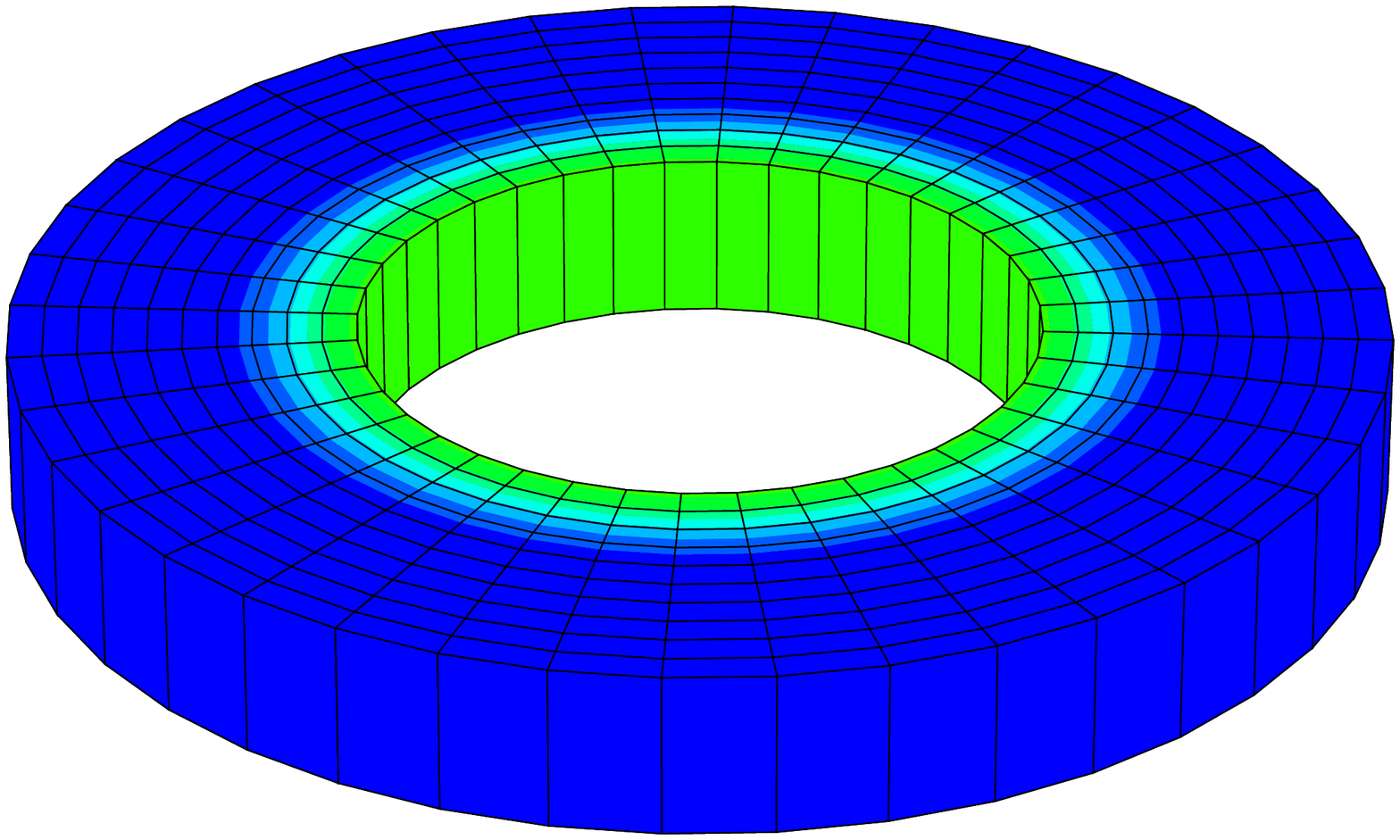}
          \end{minipage}
           \\[0mm]
          \begin{minipage}{77.5mm}
           \includegraphics[scale=0.35, clip=]{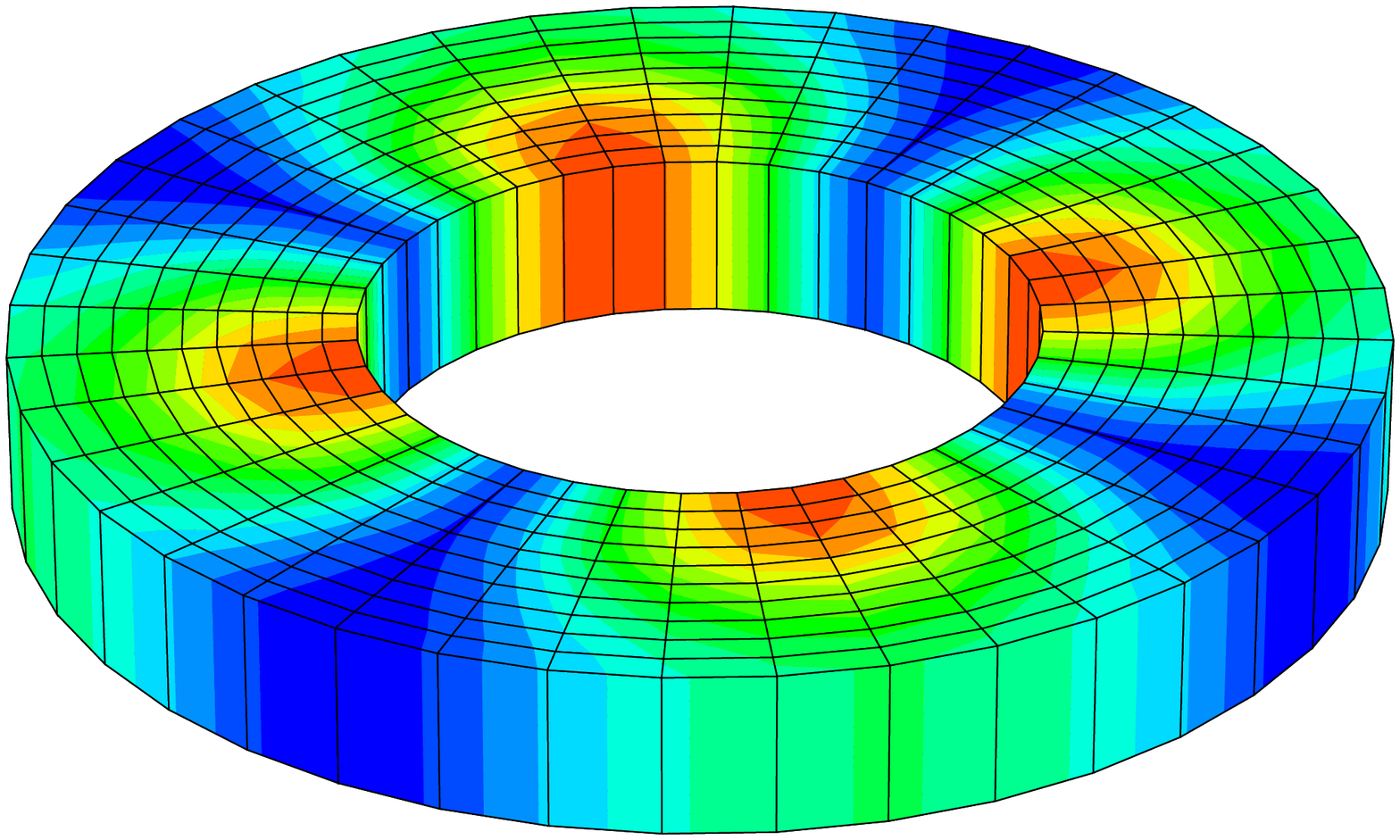}
          \end{minipage}\hfil
          \begin{minipage}{77.5mm}
            \includegraphics[scale=0.35, clip=]{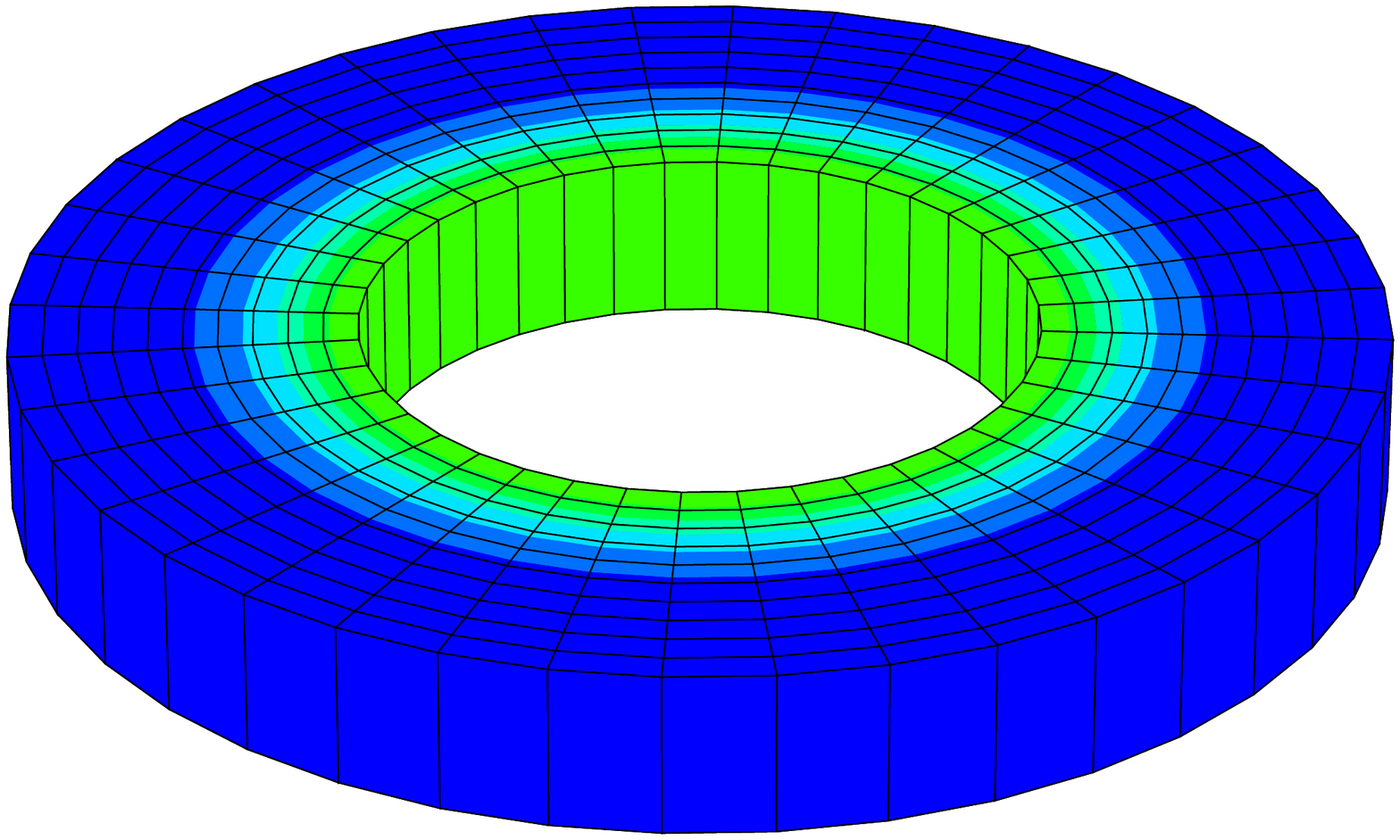}
          \end{minipage}
           \\[0mm]
          \begin{minipage}{77.5mm}
           \includegraphics[scale=0.35, clip=]{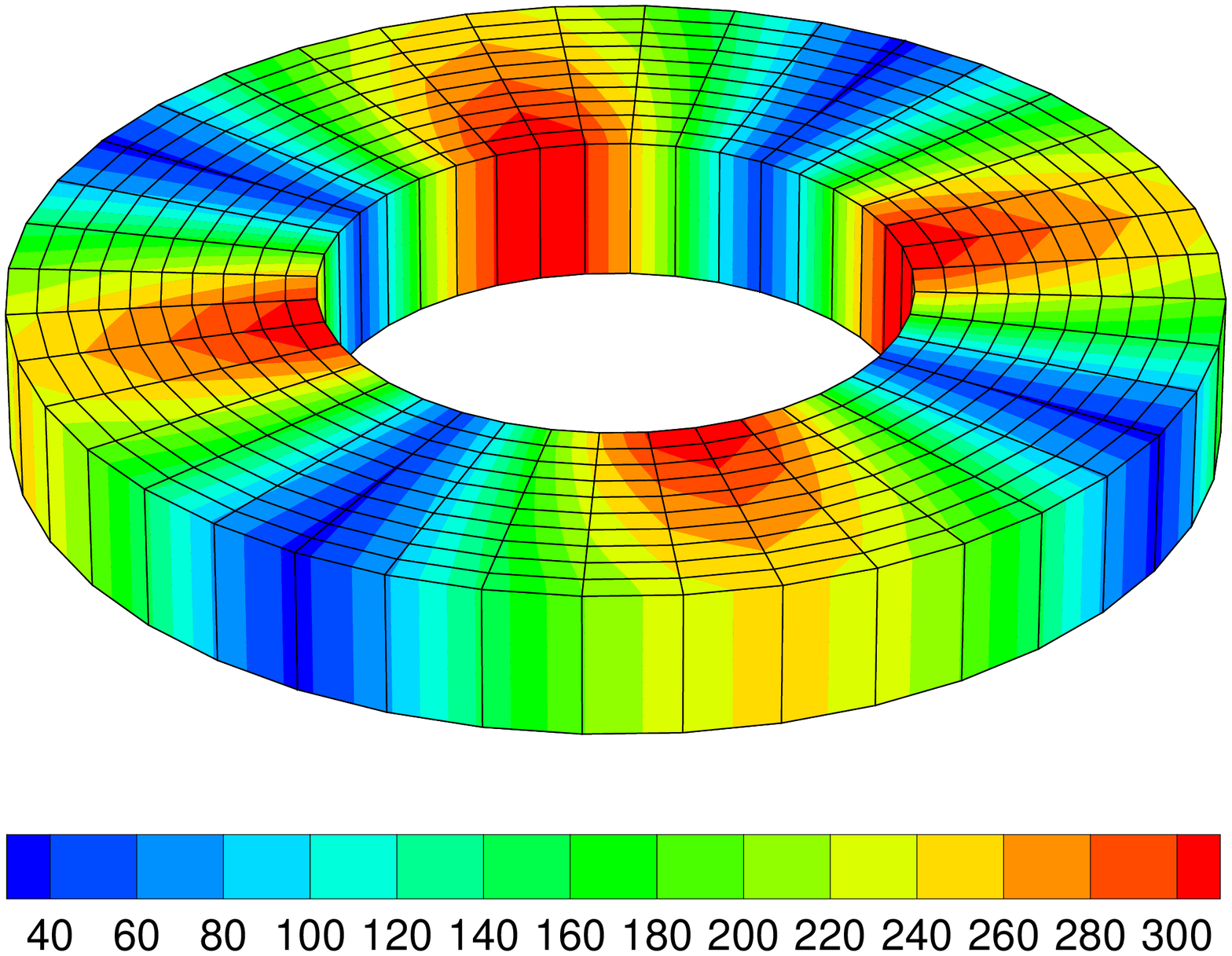}
          \end{minipage}\hfil
          \begin{minipage}{77.5mm}
            \includegraphics[scale=0.35, clip=]{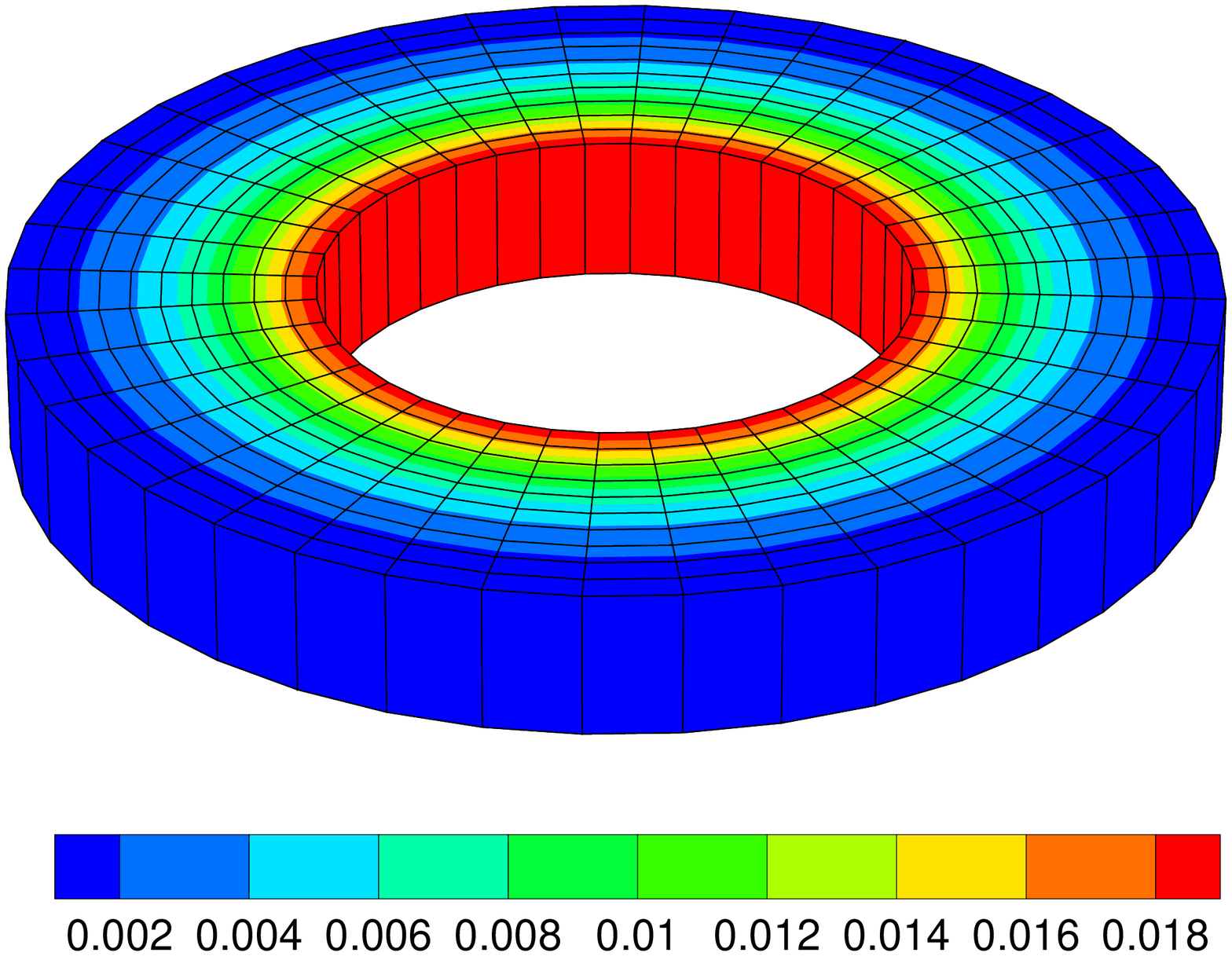}
          \end{minipage}
  \end{minipage}
\caption{Radial contraction of an elasto-plastic annulus for {\bf case B}: (left) shear-stress $S_{xy}$ $({\text{N}}/{\text{mm}^2})$, 
(right) equivalent plastic strain $\alpha$ for linearly increasing contraction of the annulus from top to bottom.
\label{fig:Annulus-CaseC-Stress-EqPlStr}}
\end{Figure}

\vfill

\subsubsection{Discussion: The signature of elasto-plasticity.}

The results of the simulations show, that it is not purely the hardening rate, which is the decisive agency for the convergence order.
Instead, it is the smoothness of the elasto-plastic stiffness that directly transfers to a corresponding smoothness
of the strain path in time. 
 
\begin{Figure}[htbp]
\center
      \includegraphics[width=9.2cm, angle=0, clip=]{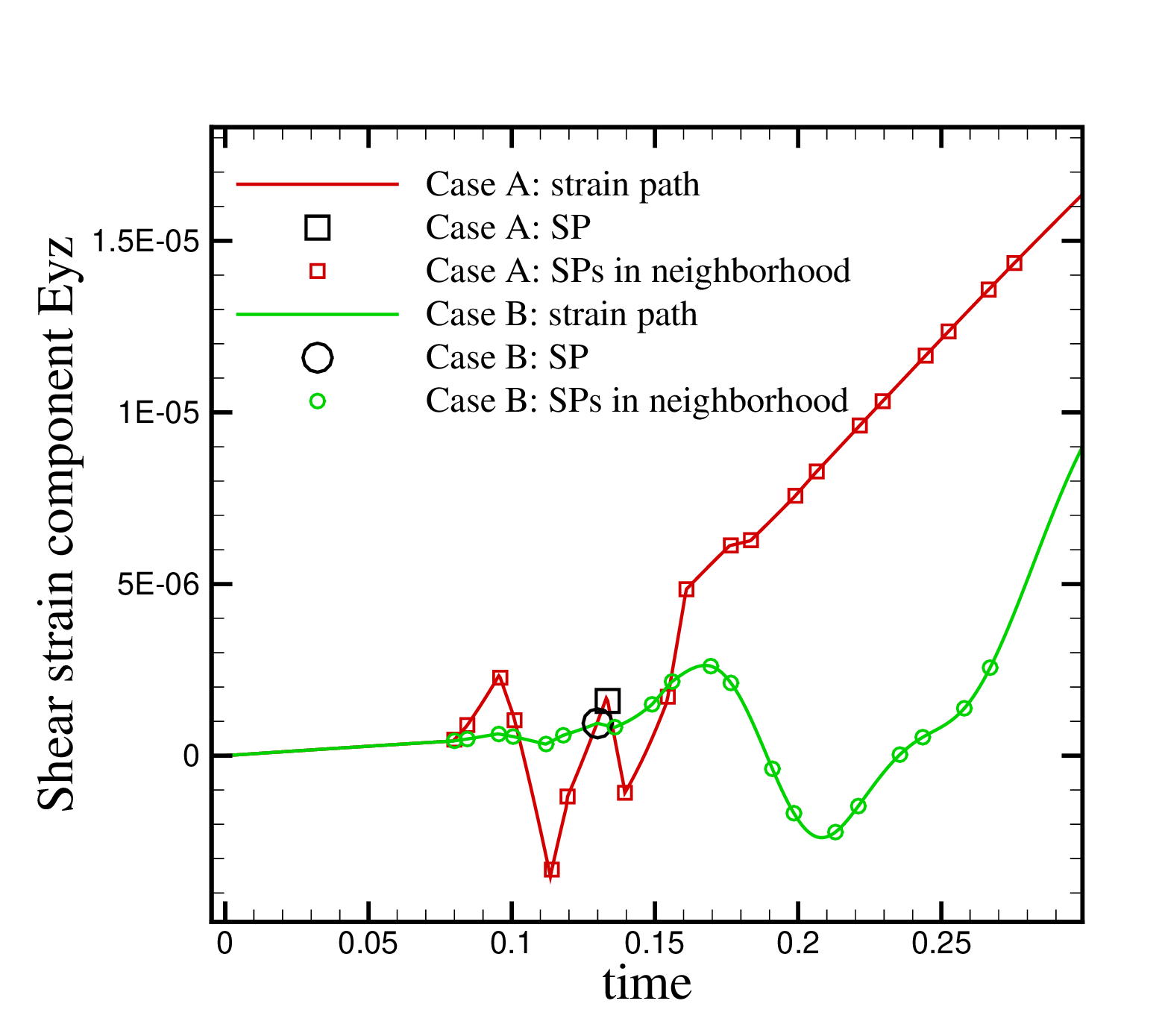}
   %   \hspace*{2mm} \scalebox{1.0}{\input{bild1202f.pstex_t}}
   %   \hspace*{2mm} \scalebox{1.0}{\input{bild1202f.pstex_t}}
      \includegraphics[width=.36\textwidth, angle=0, clip=]{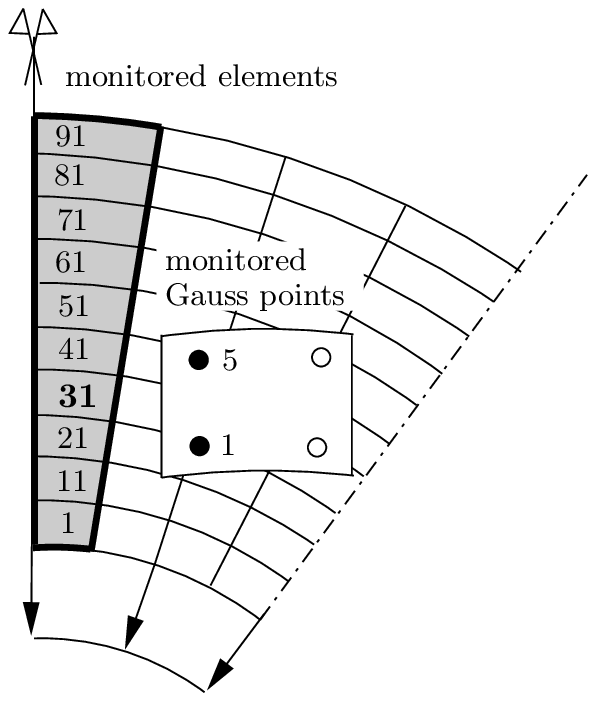}
\caption{Radial contraction of an elasto-plastic annulus, case A vs. case B. {\bf Left:} strain component $E_{yz}$ in {\bf Gauss point 1 of 
element 31} as a function of time. The symbols $\Box$ (for case A) and $\bigcirc$ (for case B) each indicate the event of a elastic-plastic 
switching occurring in the neighborhood of the considered Gauss point (31/1). {\bf Right:} Among the points in the neighborhood we monitor 
the Gauss points (1, 5) in the elements (1, 11, 21, \ldots, 91). Each of these SPs are displayed in the strain paths in the 
left diagram.
\label{fig:Kreisring_epsSP}}
\end{Figure}

The fundamental difference in the convergence between case A and case B shall be underpinned by data. To this end, the evolution of the strain 
component $E_{yz}$ is instrumental. Here, we choose Gauss point 1 in element 31, see the right of Fig.~\ref{fig:Kreisring_epsSP}, for which 
$E_{yz}$ is monitored during the contraction of the annulus from time $t=0$ to $t=0.5$. 
% (coordinates $x= ....$, $y= ....$, $z= ....$,).
A comparison of case A with case B reveals remarkable differences. Right after the beginning at $t=0$, the strain paths are identical, since the 
deformation of the entire structure is still purely elastic and materials A and B have identical elasticity laws. At time $t=0.08$ plastic yielding 
commences at the inner rim of the annulus (Gauss point 1 of the first inner element row) and then, continuously spreads out to the outer rim. For 
case A, we observe several kinks in the strain path which each coincide with the event, that a particular Gauss point (the $\Box$ symbol in 
Fig.~\ref{fig:Kreisring_epsSP}) exhibits a SP. 

Remarkably, such a kink in the strain path of the monitored material point (Gauss point 1 of element 31) does not only occur once, namely for the 
elastic-plastic transition for this particular point itself. Instead, suchlike kinks equally occur \emph{before} and \emph{after} the plastic 
frontier has been passing that particular material point, see Fig.~\ref{fig:Kreisring_epsSP} (left). Hence, we can observe in a particular Gauss-point,
when other material points in its neighborhood undergo elastic-plastic transition, which is typically accompanied by an abrupt loss of material 
stiffness.
 
The observed mutual interference of different quadrature points in terms of kinks in the strain path is due to stress redistribution which is a 
signature for structures undergoing elasto-plastic deformations. Opposed to homogeneous deformations as in the considered test sets of biaxial 
stretch, plastic flow commences in this elasto-plastic structure subsequently and not simultaneously. By the particular choice of case A0 to set 
$\sigma_Y=0$ plastic flow starts in all points at the very beginning and no kinks in strain path are present.

The kinks in the strain-path of Fig.~\ref{fig:Kreisring_epsSP} visualize what is referred to as \emph{''low regularity''} of elasto-plasticity. 
The fact that strains are not continuously differentiable in time is not a deficiency or artefact of the time-discrete algorithmic treatment of 
elasto-plasticity, and it is not an outcome of the spatially discretized problem. Quite in contrast, it is an inherent property of the very 
continuous problem of elasto-plasticity. As a consequence, the application of nonlinear strain interpolation by higher-order polynomials as 
proposed in Sec.~\ref{subsec:DehnInterpol} is of limited use. 
For non-smooth strain paths as in the DAE-case of elasto-plasticity, however, this concept can hardly cure order-reduction in all cases. 

For case B, which is smooth in its elasto-plastic stress-strain curve, we observe -- in contrast to case A -- an equally smooth strain 
path in terms of strain component $E_{yz}$ of Gauss point 1 of element 31; here, the event of elastic-plastic switching (symbol $\bigcirc$ in 
Fig.~\ref{fig:Kreisring_epsSP}) in neighboring points does not induce a kink. For that case, which comes close to the viscoelastic ODE-case, 
see \cite{EidelKuhn2010}, the consistent, higher-order strain interpolation is effective and overcomes the order-reduction.   

While some modeling approaches as e.g. in \cite{HollensteinJabareenRubin2013} assume for algorithmic convenience a smooth elastic-inelastic
transition, the smooth transition in the stress-strain curve as in case B is typical for a class of metals and metallic alloys. 
A look into timely experimental data underpins, that this class of metals is very broad, covering not only aluminum and its alloys, 
\cite{KhanMeredith2010}, but also high-performance dual-phase steels, \cite{Khan-etal-HighStrengthSteels2010}, as well as e.g. high-purity, 
polycrystalline $\alpha$-titanium, \cite{NixonCazacuLebensohn2010}, to mention but a few.

%----------------------------------------------------------------------------------------------------------------------------------------
\section{Summary and conclusions}
\label{sec:ElPl-Conclusions}
In this paper we have identified and analyzed the reasons for order reduction in computational elasto-plasticity, where the evolution 
equations for plastic flow exhibit the format of DAEs. For a consistent space-time coupling we have 
established novel consistency conditions and have proposed simple and effective remedies to fulfil them. The main findings shall be 
summarized:
\\
\begin{enumerate}
\item[(I)] {\bf Reasons.} The reasons for order reduction in computational elasto-plasticity from predicted order $p\geq3$ to order 2 is twofold. It is (i) an
           \emph{inaccurate} approximation of the displacement/strain path in time for the calculation of RK-stage values introducing a local interpolation
           error that propagates to the convergence order of stress-update algorithms. If the strain path is non-linear but smooth, then its approximation by
           quadratic interpolation leads to superior convergence compared with linear interpolation. Suchlike ''smooth'' elasto-plasticity comes close to the
           viscoelastic case. If, instead, the strain path exhibits kinks, the remedy of high-order polynomial approximation is of course a blunt knife which
           cannot improve convergence. Suchlike kinks, however, typically occur, when a material point switches from elastic to plastic state and hence passes the SP.
           Secondly (ii), it is a missing or inaccurate SP detection between elastic phase and plastic phase which leads to an error in the initial
           values for time integration and therefore to order reduction.

\item[(II)] {\bf Remedies.} For a third order, two-stage Runge-Kutta method of Radau IIa class, we have shown by numerical tests that a necessary condition 
            for full convergence
            order 3 in the case of a smooth strain path is a quadratic interpolation of the displacements/strains in time. Here we have chosen $t_n$, $t_{n+1}$ and, 
            additionally, $t_{n-1}$ data for the construction of a quadratic interpolation polynomial. Furthermore, we have shown that the smoothness of the strain path 
            requires a sufficiently smooth transition from elastic to plastic state in the stress-strain curve. Hence, it is the hardening behavior after 
            the SP which avoids the abrupt loss of stiffness at the SP. Complementary, a constant or linear strain approximation 
            inevitably leads to an order reduction to second order. Moreover we have demonstrated, that SP detection is a necessary condition 
            for full convergence order.

\item[(III)] {\bf Benefits.} A main advantage of the present approach is, that it solves a problem. Furthermore the measurements for overcoming order reduction 
             are easy to implement and require minor changes of existing implementations based on Backward-Euler as time integrator. The necessity to provide
             total strain data of the previous time-step for quadratic strain interpolation yields a moderately increased overhead of history data. But in view
             of the observed speed-up of at least a factor 8$\times$, the 3rd RK methods are most efficient and promising for further developments.
\end{enumerate}

%The common denominator of both measurements undertaken to overcome order reduction refers to an improved
%approximation of the strain path in time. This applies to the nonlinear (here: quadratic) interpolation of
%strains or a bi-linear approximation in case of a kink due to a switching point.
%On the other hand, switching point detection, which defines proper initial conditions for time integration,
%implies to accurately decompose the current, total strain increment into its true elastic part and plastic part, respectively.
%\\[2mm]
For materials, which exhibit a sharp kink in the stress strain curve indicating an abrupt loss of material stiffness at the SP, it is very 
likely that this property directly transfers into a corresponding non-smooth strain path in time. This non-smooth characteristics of the 
strain paths substantiate the frequently used, but not very sharp attribute ''\emph{low regularity}'' of elasto-plasticity. Here, the present 
work has elucidated this mathematical, but still vague notion by linking it to the mechanics of elasto-plastic solids.
\\[2mm]
Great many metals exhibit a smooth elastic-to-plastic transition in the stress-strain curve. For these materials a convergence order of 3 
can be obtained, if the algorithmic means of the present work are used and the strain paths are smooth. Then, RK methods of order 3 can 
be most efficient for time integration of elasto-plastic constitutive models as we have underpinned by the obtained speed-up. 
\\[3mm] 
{\bf Outlook.} In the future, the presented high-order time integration methods for elasto-plasticity should be tested for more complex
deformation histories with more than one SP like for unloading and reloading. Moreover, in the simulation of
a large number of repeated cyclic loading the accumulation of the integration errors should be studied. This will reveal
the long term behavior of the improved RK integration algorithms, which is necessary for applications like fatigue simulations.

%------------------------------------------------------gg------------------------------------------------------------------------------

\begin{appendix}

\addcontentsline{toc}{section}{Appendix}
\renewcommand{\thesubsection}{\Alph{subsection}}
\renewcommand{\theequation}{\Alph{subsection}.\arabic{equation}}
\renewcommand{\thefigure}{\Alph{subsection}.\arabic{figure}}
\renewcommand{\thetable}{\Alph{subsection}.\arabic{table}}
\newcommand {\ssectapp}{
                        \setcounter{equation}{0}
                        \setcounter{figure}{0}
                        \setcounter{table}{0}
		                \subsection
                        }

\vfill
\newpage

%\begin{center}
%{APPENDIX A}
%\end{center}

%-------------------------------------------------------------------------------------------------------------------------------
\ssectapp{Butcher tableaus}
\label{sec:RKSolutionStepsForDAEs}
%==============================================================================================

The RK methods, which are used in the present work, are displayed in the so-called \textit{Butcher-array}, 
in Tab.~\ref{tab:Butcher-Array}. For $s$-stage methods, $b,c\in\Real^s$ and $A\in\Real^{s\times s}$. 
%Between the lines $a_{i.}$ of the RK matrix and the stages $c_i$ the following relationship holds $c_i = \Sum_{j=1}^s a_{ij}\,.$
{\color{black} For the family of Radau IIa methods, where it holds for the order $p=2s-1$, the Butcher 
schemes are shown for $p=1,2,3$ in the second row of Tab.\ref{tab:Butcher-Array}}.

\begin{Table}[htbp]
\center
\begin{minipage}[b]{0.10\textwidth}
\center
\begin{tabular}{c|c}
     $c$ &   $A$  \\
\hline   &   $b$  \\
\end{tabular}
\end{minipage}
%\hspace*{4mm} 
\begin{minipage}[b]{0.35\textwidth}
\renewcommand{\arraystretch}{1.5}
\center
\begin{tabular}{c|cccc}
   $c_1$   &  $a_{11}$  & $a_{12}$  &  $\ldots$  &  $a_{1s}$ \\
   $c_2$   &  $a_{21}$  & $a_{22}$  &  $\ldots$  &  $a_{2s}$ \\
 $\vdots$  &  $\vdots$  &           &  $\ddots$  &  $\vdots$ \\
   $c_s$   &  $a_{s1}$  & $a_{s2}$  &  $\ldots$  &  $a_{ss}$ \\
                                                     \hline
           &   $b_1$    &  $b_2$    &  \ldots    &  $b_s$    \\
\end{tabular}
\end{minipage}
%\hspace*{4mm} 
\begin{minipage}[b]{0.35\textwidth}
\renewcommand{\arraystretch}{1.5}
\center
\begin{tabular}{c|cccc}
   $c_1$   &  $a_{11}$  & $a_{12}$   &  $\ldots$  &   $a_{1s}$  \\
   $c_2$   &  $a_{21}$  & $a_{22}$  &  $\ldots$  &    $a_{2s}$  \\
 $\vdots$  &  $\vdots$  &           &  $\ddots$  &    $\vdots$  \\
   $1$     &   $b_1$    &   $b_2$   &  $\ldots$  &   $b_s$  \\
                                                     \hline
           &   $b_1$    &   $b_2$   &  $\ldots$  &  $b_s$    \\
\end{tabular}
\end{minipage}

\vspace{5mm}
 
\begin{minipage}[b]{0.10\textwidth}
\center
\begin{tabular}{c|c}
     $1$ &   $1$  \\
\hline   &   $1$  \\
\end{tabular}
\end{minipage}
\begin{minipage}[b]{0.30\textwidth}
\center
\renewcommand{\arraystretch}{1.5}
\begin{tabular}{c|cc}
%               &                &                 \\
 $\dfrac{1}{3}$ & $\dfrac{5}{12}$ & $-\dfrac{1}{12}$ \\[2mm]
           $1$ & $\dfrac{3}{4}$  & $\dfrac{1}{4}$   \\[2mm]
                                             \hline \\[-6mm]
               & $\dfrac{3}{4}$  & $\dfrac{1}{4}$   \\
\end{tabular}
\end{minipage}
\begin{minipage}[b]{0.40\textwidth}
\center
\renewcommand{\arraystretch}{1.5}
\begin{tabular}{c|ccc}
 $\frac{4- \sqrt{6}}{10}$ & $\frac{88-7\sqrt{6}}{360}$ & $\frac{296-169\sqrt{6}}{1800}$ & $\frac{-2+3\sqrt{6}}{225}$ \\
 $\frac{4+\sqrt{6}}{10}$ & $\frac{296+169\sqrt{6}}{1800}$ & $\frac{88+7\sqrt{6}}{360}$ & $\frac{-2-3\sqrt{6}}{225}$ \\
 $1$ & $\frac{16-\sqrt{6}}{36}$ & $\frac{16+\sqrt{6}}{36}$ & $\frac{1}{9}$ \\
\hline  & $\frac{16-\sqrt{6}}{36}$ & $\frac{16+\sqrt{6}}{36}$ & $\frac{1}{9}$ \\
\end{tabular}
\end{minipage}
\caption{Butcher arrays for (1st line, from left to right) the general case, for implicit Runge-Kutta (IRK) methods and for 
stiffly accurate implicit RK methods, and for (2nd line, from left to right) Radau IIa-schemes for $s=1$, {\color{black} $p=1$},
i.e. Backward-Euler; for $s=2$, {\color{black} $p=3$}; and for $s=3$, {\color{black} $p=5$}, respectively.
\label{tab:Butcher-Array}}
\end{Table}

A nice survey on the historical development of the Radau IIa methods is given in \cite{HaWaRadau}. It has a focus on stiff problems including 
extended stability notions for stiff problems.

%\bigskip
%
%\begin{center}
%{APPENDIX B}
%\end{center}

%-------------------------------------------------------------------------------------------------------------------------------
\ssectapp{Algorithmically consistent tangent 
%consistent to the Radau IIa scheme
% using linear/quadratic strain interpolation
}
\label{sec:RKSolutionStepsForDAEs}
%==============================================================================================

Introducing the fourth-order deviatoric tensor $\mathbb{P}$ with $\mathbb{P} := \mathbb{I}\text{d} - \frac{1}{3}\1\otimes\1$ where 
$\mathbb{I}\text{d}$ is the 4th order unity tensor with components $\delta_{ij}\delta_{kl}$ such that $\mathbb{P}:\bE =\bE^D$, the stress tensor 
$\bS_{n+1}$ can be written as
\begin{eqnarray}
\label{Stress-strain}
\bS_{n+1} &=& \kappa\, \tr(\bE_{n+1})\1 + \bS^D_{n+1}                     \nonumber \\
&=& \kappa (\1\otimes\1):\bE_{n+1} + 2\mu(\mathbb{P}:\bE_{n+1} - \bE^{\mbox{\scriptsize p}}_{n+1})           \\
&=& \mathbb{C}:\bE_{n+1} - 2\mu\bE^{\mbox{\scriptsize p}}_{n+1}\,.                                              \nonumber
\end{eqnarray}
Based on representation \eqref{Stress-strain}$_3$, the consistent tangent can be written as
\begin{equation}
\mathbb{C}_{n+1}^{\mbox{\scriptsize ep}} = \parziell{\bS_{n+1}}{\bE_{n+1}}= \mathbb{C}- 2\mu\parziell{\bE^{\mbox{\scriptsize p}}_{n+1}}{\bE_{n+1}} \, .
\label{eq:Elpl-Tangent-General}
\end{equation}
Differentiation of \eqref{Solution}$_1$ and \eqref{Solution}$_3$ with respect to $\bE_{n+1}$
yields a linear set of equations for calculating
$\dfrac{\partial\bE^{\mbox{\scriptsize p}}_{n+1}}{\partial\bE_{n+1}} = \dfrac{\partial\bE^{\mbox{\scriptsize p}}_{ns}}{\partial\bE_{n+1}}$.
%With representation \eqref{Stress-strain}$_3$, the consistent tangent can be written as
%\begin{equation}
%\mathbb{C}^{\mbox{\scriptsize ep}} = \parziell{\bS_{n+1}}{\bE_{n+1}}= \mathbb{C}- 2\mu\parziell{\bE^{\mbox{\scriptsize p}}_{n+1}}{\bE_{n+1}} \, .
%%\label{eq:Elpl-Tangent-General}
%\end{equation}
%Differentiation of \eqref{Solution}$_1$ and \eqref{Solution}$_3$ with respect to $\bE_{n+1}$
%yields a linear set of equations for calculating
%$\dfrac{\partial\bE^{\mbox{\scriptsize p}}_{n+1}}{\partial\bE_{n+1}} = \dfrac{\partial\bE^{\mbox{\scriptsize p}}_{ns}}{\partial\bE_{n+1}}$.

Introducing $\bN_{nj} := \dfrac{\bE^D_{nj} - \bE^{\mbox{\scriptsize p}}_{nj}}{||\bE^D_{nj} -\bE^{\mbox{\scriptsize p}}_{nj} ||}$, differentiation of \eqref{Solution}$_1$ yields
\begin{eqnarray}
\parziell{\bE^{\mbox{\scriptsize p}}_{ni}}{\bE_{n+1}} &=& \Sum_{j=1}^s a_{ij} \parziell{ }{\bE_{n+1}}(\Delta\Gamma_{nj} \bN_{nj})\nonumber\\
&=& \Sum_{j=1}^s a_{ij} \biggl( \Delta\Gamma_{nj} \parziell{\bN_{nj}}{\bE_{n+1}} + \bN_{nj}\otimes\parziell{\Delta\Gamma_{nj}}{\bE_{n+1}}\biggl)\,.\label{AblepTeil}
\end{eqnarray}
With $\bx_{nj} := \bE^D_{nj} - \bE^{\mbox{\scriptsize p}}_{nj}$ and the chain rule the first derivative in \eqref{AblepTeil}$_2$ is readily obtained  
\begin{eqnarray}
\parziell{\bN_{nj}}{\bE_{n+1}}&=& \parziell{\bN_{nj}}{\bx_{nj}}:\parziell{\bx_{nj}}{\bE_{n+1}} \nonumber \\
&=& \biggl[\frac{\mathbb{I}\text{d} - \bN_{nj}\otimes\bN_{nj}}{{||\bx_{nj} ||}}\biggl]:\biggl[\parziell{\bE^D_{nj}}{\bE_{n+1}}-\parziell{\bE^{\mbox{\scriptsize p}}_{j}}{\bE_{n+1}}\biggl]\,.
\end{eqnarray}
Inserting this result into \eqref{AblepTeil} yields
\begin{eqnarray*}
\parziell{\bE^{\mbox{\scriptsize p}}_{ni}}{\bE_{n+1}} &=&\Sum_{j=1}^s a_{ij} \biggl( \Delta\Gamma_{nj} \biggl[\frac{\mathbb{I}\text{d} - \bN_{nj}\otimes\bN_{nj}}{||\bE^D_{nj} - \bE^{\mbox{\scriptsize p}}_{j} ||}\biggl]:\biggl[\parziell{\bE^D_{nj}}{\bE_{n+1}}-\parziell{\bE^{\mbox{\scriptsize p}}_{nj}}{\bE_{n+1}}\biggl] + \bN_{nj}\otimes\parziell{\Delta\Gamma_{nj}}{\bE_{n+1}}\biggl)\\
% &\Longleftrightarrow &
\end{eqnarray*}
\begin{equation*}
-\parziell{\bE^{\mbox{\scriptsize p}}_{ni}}{\bE_{n+1}}+\Sum_{j=1}^s a_{ij}\biggl( \bN_{nj}\otimes\parziell{\Delta\Gamma_{nj}}{\bE_{n+1}} - \Delta\Gamma_{nj} \biggl[\frac{\mathbb{I}\text{d} - \bN_{nj}\otimes\bN_{nj}}{||\bE^D_{nj} - \bE^{\mbox{\scriptsize p}}_{nj} ||}\biggl]:\parziell{\bE^{\mbox{\scriptsize p}}_{nj}}{\bE_{n+1}} \biggl)\hspace{2cm}
\end{equation*}
\begin{equation}
\hspace{5cm}= -\Sum_{j=1}^s a_{ij} \biggl( \Delta\Gamma_{nj} \biggl[\frac{\mathbb{I}\text{d} - \bN_{nj}\otimes\bN_{nj}}{||\bE^D_{nj} - \bE^{\mbox{\scriptsize p}}_{nj} ||}\biggl]:\parziell{\bE^D_{nj}}{\bE_{n+1}} \biggl) \,.
\label{ablre}
\end{equation}
If the stage values of strain are constructed by linear interpolation according to \eqref{linapproxE}, we obtain
\begin{eqnarray}
\parziell{\bE^D_{nj}}{\bE_{n+1}} &=& \parziell{ }{\bE_{n+1}}\left(\bE^D_{n} + c_{nj}(\bE^D_{n+1}-\bE^D_{n}) \right)  \nonumber \\
                                  &=& c_{nj}\,\mathbb{P}\,.
\end{eqnarray}
For quadratic interpolation it holds
\begin{eqnarray}
\parziell{\bE^D_{nj}}{\bE_{n+1}} &=& \parziell{ }{\bE_{n+1}}\biggl(\frac{c_j}{2}(c_j-1)\,\bE^D_{n-1}+(1-c_j^2)\,\bE^D_n + \frac{c_j}{2}(c_j+1)\,\bE^D_{n+1} \biggl) \nonumber \\
                                    &=& \frac{c_j}{2}(c_j+1)\,\mathbb{P}\,.
\end{eqnarray}
For notational convenience we avoid explicit case differentiation in the following by introducing a new stage variable $\bar{c}_j$.
\begin{equation}
\bar{c}_j:=\begin{cases} c_j & \text{, for linear interpolation}\\
\dfrac{c_j}{2}(c_j+1) & \text{, for quadratic interpolation} \end{cases}
\end{equation}
Inserting $\bar{c}_j$ into the right-hand side of \eqref{ablre} yields
\begin{eqnarray*}
\biggl[\frac{\mathbb{I}\text{d} - \bN_{nj}\otimes\bN_{nj}}{||\bE^D_{nj} - \bE^{\mbox{\scriptsize p}}_{nj} ||}\biggl]:\parziell{\bE^D_{nj}}{\bE_{n+1}} &=&\bar{c}_j \frac{\mathbb{P} - \bN_{nj}\otimes\bN_{nj}}{||\bE^D_{nj} - \bE^{\mbox{\scriptsize p}}_{nj} ||}\,,
\end{eqnarray*}
and finally \eqref{ablre} is obtained in the form

\fbox{ \parbox{15.2cm}{
\begin{eqnarray}
-\parziell{\bE^{\mbox{\scriptsize p}}_{ni}}{\bE_{n+1}}+\Sum_{j=1}^s a_{ij}\biggl( \bN_{nj}\otimes\parziell{\Delta\Gamma_{nj}}{\bE_{n+1}} - \Delta\Gamma_{nj} \biggl[\frac{\mathbb{I}\text{d} - \bN_{nj}\otimes\bN_{nj}}{||\bE^D_{nj} - \bE^{\mbox{\scriptsize p}}_{nj} ||}\biggl]:\parziell{\bE^{\mbox{\scriptsize p}}_{nj}}{\bE_{n+1}} \biggl)\nonumber\\
= -\Sum_{j=1}^s a_{ij} \frac{\bar{c_j}\Delta\Gamma_{nj}}{||\bE^D_{nj} - \bE^{\mbox{\scriptsize p}}_{nj} ||} (\mathbb{P} - \bN_{nj}\otimes\bN_{nj})\, .
\label{cepsyst1}
\end{eqnarray}
}}

Moreover, \eqref{Solution}$_3$ has to be differentiated with respect to $\bE_{n+1}$.
\begin{eqnarray}
0 &=& 2\mu\parziell{||\bE^D_{ni} - \bE^{\mbox{\scriptsize p}}_{ni} ||}{\bE_{n+1}}-\sqrt{\frac{2}{3}}\ddK(\Lambda_{ni})\parziell{\Lambda_{ni}}{\bE_{n+1}} \nonumber\\
0 &=& 2\mu\,\bN_{ni} :\parziell{(\bE^D_{ni} - \bE^{\mbox{\scriptsize p}}_{ni})}{\bE_{n+1}}-\frac{2}{3}\ddK(\Lambda_{ni})\Sum_{j=1}^s a_{ij} \parziell{\Delta\Gamma_{nj}}{\bE_{n+1}} \nonumber\\
-2\mu\,\bN_{ni}:(\bar{c}_{i}\,\mathbb{P}) &=&-2\mu\,\bN_{ni}:\parziell{\bE^{\mbox{\scriptsize p}}_{ni}}{\bE_{n+1}} -\frac{2}{3}\ddK(\Lambda_{ni})\Sum_{j=1}^s a_{ij} \parziell{\Delta\Gamma_{nj}}{\bE_{n+1}} \nonumber\\
\bar{c}_{i}\,\bN_{ni} &=&\bN_{ni}:\parziell{\bE^{\mbox{\scriptsize p}}_{ni}}{\bE_{n+1}}+\frac{\ddK(\Lambda_{ni})}{3\mu}\Sum_{j=1}^s a_{ij} \parziell{\Delta\Gamma_{nj}}{\bE_{n+1}} \nonumber 
%\bar{c_i}\frac{\bE^D_{ni} - \bE^{\mbox{\scriptsize p}}_{ni}}{||\bE^D_{ni} -\bE^{\mbox{\scriptsize p}}_{ni}||} &=&\frac{\bE^D_{ni} - \bE^{\mbox{\scriptsize p}}_{ni}}{||\bE^D_{ni} -\bE^{\mbox{\scriptsize p}}_{ni}||}:\parziell{\bE^{\mbox{\scriptsize p}}_{ni}}{\bE_{n+1}}+\frac{\ddK(\Lambda_{ni})}{3\mu}\Sum_{j=1}^s a_{ij} \parziell{\Delta\Gamma_{nj}}{\bE_{n+1}}\label{cepsyst2}
\end{eqnarray}

\fbox{ \parbox{15.2cm}{
\begin{eqnarray}
\bar{c_i}\frac{\bE^D_{ni} - \bE^{\mbox{\scriptsize p}}_{ni}}{||\bE^D_{ni} -\bE^{\mbox{\scriptsize p}}_{ni}||} &=&\frac{\bE^D_{ni} - \bE^{\mbox{\scriptsize p}}_{ni}}{||\bE^D_{ni} -\bE^{\mbox{\scriptsize p}}_{ni}||}:\parziell{\bE^{\mbox{\scriptsize p}}_{ni}}{\bE_{n+1}}+\frac{\ddK(\Lambda_{ni})}{3\mu}\Sum_{j=1}^s a_{ij} \parziell{\Delta\Gamma_{nj}}{\bE_{n+1}}\label{cepsyst2}
\end{eqnarray}
}}

\vspace*{2mm}

Equations \eqref{cepsyst1} and \eqref{cepsyst2} form the set of linear equations for the calculation of  $\dparziell{\bE^{\mbox{\scriptsize p}}_{ni}}{\bE_{n+1}}$ and
$\dparziell{\Delta\Gamma_{ni}}{\bE_{n+1}}$ such that $\dparziell{\bE^{\mbox{\scriptsize p}}_{n+1}}{\bE_{n+1}} = \dparziell{\bE^{\mbox{\scriptsize p}}_{ns}}{\bE_{n+1}}$ 
in \eqref{eq:Elpl-Tangent-General} can be calculated which completes the exercise.

%\vfill
%\newpage

%\input{./sec_appendix/Ten2Mat}  % Tensor representation by Matrices
%\cleardoublepage
%\input{./sec_appendix/DefConvexity}
%\cleardoublepage
%\input{./sec_appendix/FbarIncVars} % Incremental and variational quantities for the Fbar-approach

% \newpage

%\input{appendix_C}
%\input{appendix_D} 

\end{appendix}

%-----------------------------------------------------------------------------------------------------------

%=========================================================================================

% ------- layout-datei --------------
\bibliographystyle{plainnat}
%\bibliographystyle{plainnat}
%\bibliographystyle{plaindin}
%\bibliographystyle{plaindin_shortname2}
%\bibliographystyle{elsarticle-num-names}
% ------- bib-datei --------------
%\bibliography{js_master_2009}
%\bibliographystyle{model2-names}
%\bibliography{Plasticity}
%\begin{appendix}
%\include{appendix_A}
%\end{appendix}
\end{document}